\documentclass[12pt]{article}
\usepackage{amssymb}
\usepackage{graphicx}
\usepackage{epsfig}
\usepackage{color}
\include{psfig}

\newcommand{\eq}[1]{$(\ref{#1})$}

\newcommand{\ud}{\mathrm{d}}
\newtheorem{theorem}{Theorem}[section]
\newtheorem{proposition}{Proposition}[section]
\newtheorem{lemma}{Lemma}[section]
\newtheorem{corollary}{Corollary}[section]
\newtheorem{definition}{Definition}[section]
\newcommand{\bean}{\begin{eqnarray*}}
\newcommand{\eean}{\end{eqnarray*}}
\newcommand{\bea}{\begin{eqnarray}}
\newcommand{\eea}{\end{eqnarray}}
\newcommand{\tod}{\stackrel{{\cal D}}{\longrightarrow}}
\newcommand{\eqd}{\stackrel{{\cal D}}{=}}
\newcommand{\toP}{\stackrel{P}{\longrightarrow}}
\newcommand{\toas}{\stackrel{{\rm a.s.}}{\longrightarrow}}
\newcommand{\inL}{\stackrel{L^1}{\longrightarrow}}
\newcommand{\inLL}{\stackrel{L^2}{\longrightarrow}}
\newcommand{\rem}{\noindent \textbf{Remark. }}
\newcommand{\rems}{\noindent \textbf{Remarks. }}
\newcommand{\proof}{\noindent \textbf{Proof. }}
\def\Exp{E}
\def\Pr{P}
\def\Var{{\mathrm{Var}}}
\def\Cov{{\rm Cov}}
\def\R{{\bf R }}
\def\Z{{\bf Z }}
\def\1{{\bf 1 }}
\def\bx{{\bf x }}
\def\bz{{\bf z }}
\def\bw{{\bf w }}
\def\by{{\bf y }}
\def\N{{\bf N }}
\def\NN{{\cal N }}

\def\potp{\stackrel{\theta,\phi}{\preccurlyeq}}
\def\postar{\preccurlyeq^*}

\def\0{{\bf 0}}
\def\SS{{\cal S}}
\def\RR{{\cal R}}
\def\QQ{{\cal Q}}
\def\X{{\cal X}}

\def\FF{{\cal F}}
\def\GG{{\cal G}}

\def\VV{{\cal V}}
\def\bX{\mathbf{X}}

\def\tLL{\tilde{\cal L}}

\def\tD{\tilde{D}}
\def\tF{\tilde{F}}
\def\bV{\mathbf{V}}
\def\bU{\mathbf{U}}

\def\Po{{\cal P}}
\def\H{{\cal H}}
\def\LL{{\cal L}}

\def\UU{{\cal U}}
\def\card{{\rm  card}}
\def\dist{{\rm  dist}}
\def\diam{{\rm  diam}}

\def\eps{{\varepsilon}}

\def\dinf{{\Delta ( \infty )}}

\def\tDone{\tD_1}
\def\tFone{\tF_1}
\def\Falph{F_\alpha}
\def\Dalph{D_\alpha}
\def\tFalph{\tF_\alpha}
\def\tDalph{\tD_\alpha}
\def\Walph{W_\alpha}
\def\Yalph{Y_\alpha}
\def\tLone{\tilde \LL^1}
\def\tLalph{\tilde \LL^\alpha}
\def\Lalph{\LL^\alpha}                                                                                
\title{{\bf On the total length of the random minimal directed spanning tree}}
\author{Mathew D. Penrose$^{1}$ and Andrew R. Wade$^{2}$ 
}
                                                                               
\date{September 2004}
\begin{document}
                                                                               
\footnotetext[1]{ Department of Mathematical Sciences, University of Bath,
Bath, BA2 7AY, England: {\texttt m.d.penrose@bath.ac.uk} }
                                                                               
\footnotetext[2]
{ Department of Mathematical Sciences, University
of Durham, South Road, Durham DH1 3LE, England: {\texttt
a.r.wade@durham.ac.uk} }

\maketitle 
\abstract{ In Bhatt and Roy's minimal directed spanning tree (MDST)
 construction for a random partially ordered set of
 points in the unit square,
 all edges must respect the ``coordinatewise'' partial order
and there must be a directed path from each vertex to a minimal
element. We study the asymptotic behaviour of the total length of this
graph with power weighted edges. The limiting distribution is
given by the sum of a normal component
away from the boundary and a contribution
introduced by the boundary effects, which can be characterized
by a fixed point equation, and is reminiscent of limits
arising in the probabilistic analysis of
certain algorithms. As the exponent
of the power weighting
increases, the distribution undergoes a phase
transition from the normal contribution being
dominant to the boundary effects dominating. In the
critical case where the weight is simple Euclidean length,
both effects contribute significantly to the limit law.
We also give a law of large
numbers for the total weight of the graph. } \\


\section{Introduction}

Recent interest in graphs, generated
over random point sets 
consisting of independent uniform points in the unit square
 by connecting nearby points
according to some deterministic rule, has been
considerable. Such graphs include the geometric graph, the nearest
neighbour graph and the minimal-length spanning tree. Many aspects
of the large-sample asymptotic theory for such graphs,
when they are locally determined in a certain sense, are by now quite
well understood. See for example
\cite{KL,penbook,penyuk1,penyuk2,SY,steelebook,yukbook}.

\begin{figure}[h!]
\begin{center}
\includegraphics[angle=0, width=0.8\textwidth]{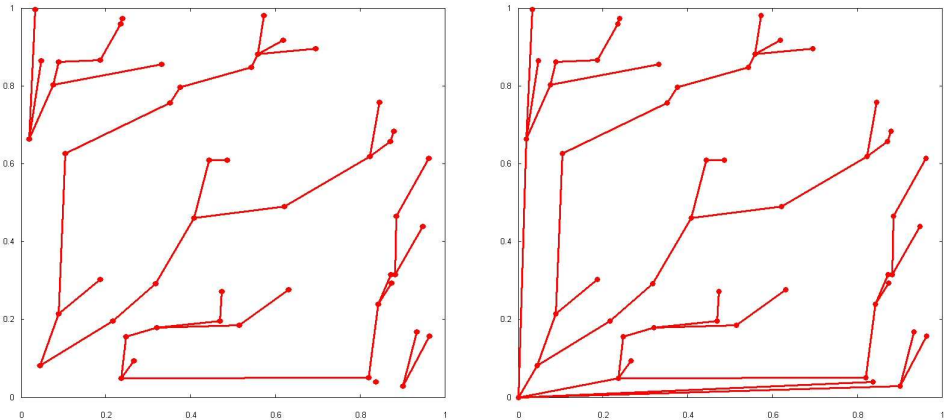}
\end{center}
\caption{Realizations of the MDSF (left) and MDST
 on 100 simulated random points in the unit square, under
 the partial ordering $\postar$.}
\label{mdstfig}
\end{figure}

One such graph is the \textit{minimal directed spanning tree} 
(or MDST for short), which was introduced
by Bhatt and Roy in~\cite{bhattroy2002}.
In the MDST, each point $\bx$ of a finite (random) subset $\SS$
of $(0,1]^2$ is connected by a directed edge to the nearest
$\by \in \SS \cup \{(0,0)\}$ such that $\by \neq \bx$ and 
$\by \preccurlyeq^* \bx$, where
$\by \preccurlyeq^* \bx$ means that
each component of $\bx -\by$ is nonnegative.
See Figure \ref{mdstfig} for a realisation of
 the MDST on simulated random points.

Motivation comes from the
modelling of communications or drainage networks
(see \cite{bhattroy2002,rooted,rodriguez}). For example,
consider the problem of designing a set of canals
to connect a set of hubs, so as to minimize their total length
subject to a constraint that all canals must flow downhill.
The mathematical formulation given above for this
 constraint can lead to significant boundary effects due
to the possibility  of long edges occurring near 
the lower and left boundaries of the unit square;
 these boundary effects distinguish
the MDST qualitatively from the standard minimal spanning tree and the
nearest neighbour graph for point sets in the plane.
Another difference is the fact that there is no
uniform upper bound on vertex degrees in the MDST.

In the present work, we consider the total length of the MDST on
random points in $(0,1]^2$, as the number of points becomes large.
We also consider the total length of the {\em minimal
directed spanning forest} (MDSF), which is the  MDST
with edges incident to the origin removed
(see Figure \ref{mdstfig} for an example). 
In \cite{bhattroy2002}, Bhatt and Roy mention that the total
length is an object of considerable interest, although they
restrict their analysis to the length of the edges joined to the origin
(subsequently also examined in~\cite{rooted}).
A first order result for the total length of the MDST
or MDSF is a law of large numbers; we 
derive this in Theorem \ref{llnthm} for
 a family of MDSFs indexed by partial orderings
 on $\R^2$, which include $\preccurlyeq^*$ as a special case.

This paper is mainly concerned with establishing second order results, i.e.,
weak convergence results for the distribution of the total length, suitably 
 centred and scaled.  For the length of edges from points
 in the region away from the boundary, we prove a central limit theorem.
The boundary effects are significant, and near the
boundary the MDST can be described in terms of a one-dimensional,
on-line version of the MDST  which we call the directed
 linear tree (DLT), and which we examine in Section 
\ref{secdlt}. In the DLT, each point in a sequence of
independent uniform random points in an interval is
joined to its nearest neighbour to the left, amongst
those points arriving earlier in the sequence.
This DLT 
is of separate interest in relation to, for example,
 network modelling and molecular fragmentation (see
 \cite{boll}, \cite{bertoin}, and references therein).

In Theorem \ref{dltthm} we establish that
 the limiting distribution of the centred total length of the DLT
 is characterized by a distributional fixed-point equation, which resembles
those encountered in the probabilistic analysis of algorithms such
as Quicksort~\cite{hoare}. Such fixed-point distributional
equalities, and the so-called `divide and conquer' or recursive algorithms
from which they arise, have received considerable attention recently; see, for
example,~\cite{hwang1998,neinrusch,rosler0,rosler}.

We consider power-weighted edges. Our weak convergence results
(Theorem \ref{mainth})
demonstrate that, depending on the value chosen for the weight
exponent of the edges, there are two regimes in which either the boundary
effects dominate or those edges away from the boundary are dominant, and that
there is a critical value (when we take simple
Euclidean length as the weight)
for which neither effect dominates.

In the related paper~\cite{rooted}, we give results dealing 
with the weight of the
edges joined to the origin, including weak convergence results,
in which the limiting distributions are given in terms of
some generalized Dickman distributions.
Subsequently, it
has been shown~\cite{blp} that this two dimensional case is rather
special -- in higher dimensions the corresponding limits 
are normally
distributed. \cite{rooted} also deals with the maximum edge
 length of the MDST (the maximum length of those edges incident
 to the origin was dealt with in~\cite{bhattroy2002}).

In the next section we give formal definitions of the MDST and MDSF,
 and state our main results (Theorems \ref{llnthm} and \ref{mainth}) 
on the total length of the MDST and MDSF. The 
 results on the DLT which we present in Section \ref{secdlt},
and  the
general central limit theorems  
which we present in Section \ref{secgeneral}, 
are of some independent interest.

\section{Definitions and main results} \label{con}
We work in the same framework as \cite{rooted}. Here we briefly
recall the relevant terminology. See \cite{rooted} for more detail.

Suppose $V$ is a finite set endowed with a partial ordering
$\preccurlyeq$. % (see e.g.~\cite{kolmogorovformin1975}).
A {\em minimal element}, or {\em sink}, of $V$ is a vertex 
$v_0 \in V$ for which there
exists no $v\in V \setminus \{v_0\}$ such that $v \preccurlyeq v_0$.
Let $V_0$ denote the set of all sinks of $V$.

The partial ordering induces a directed graph
$G=(V,E)$, with vertex set $V$ and
with edge set $E$ consisting of all ordered pairs
$(v,u)$ of distinct elements of $V$ such that $u \preccurlyeq v$.
A {\em directed spanning forest (DSF)} on $V$ 
is a subgraph $T=(V_T,E_T)$
of $(V,E)$ such that (i) $V_T=V$ and $E_T \subseteq E$, and  (ii)
 for each vertex $v \in V \setminus V_0$
there exists a unique directed path in
$T$ that starts at $v$ and ends at some sink
$u \in V_0$.  In the case where  $V_0$ consists of a
single sink,  we refer to any   DSF on $V$ as a
 \emph{directed spanning tree (DST)} on $V$.
If we ignore the orientation of edges then \cite{rooted}
 a DSF on $V$ is indeed a
forest and, if  there is just one sink, then
any DST on $V$ is a tree.

Suppose the directed graph $(V,E)$ 
carries a {\em weight function} on its edges, i.e.,
 a function $w:E
\to [0,\infty)$.  If $T$ is a DSF on $V$, we set $w(T):=
 \sum_{e \in E_T} w(e)$.
A {\em minimal directed spanning forest (MDSF)}
on $V$ (or, equivalently, on $G$), is a directed spanning forest $T$
on $V$ 
 such that $w(T) \leq w(T')$
for every DSF $T'$ on $V$.
If $V$ has a single sink, then a minimal
directed spanning forest on $V$ is called a
 {\em minimal directed spanning tree (MDST)} on $V$.

For $v\in V$, 
we say that $u\in V \setminus \{v\}$
 is a \emph{directed nearest neighbour} of $v$
if $u \preccurlyeq v$ and $w(v,u) \leq w(v,u') $ for all
$ u' \in V\setminus  \{ v \} $ such that $ u' \preccurlyeq v$.
For each $v\in V \setminus V_0$,
let $n_v$ denote  a directed nearest neighbour of $v$ (chosen
arbitrarily if $v$ has more than one directed nearest neighbour).
Then \cite{rooted} the subgraph
  $(V,E_M)$ of $(V,E)$, obtained by taking
$
E_M := \{ (v,n_v): v \in V \setminus V_0 \},
$
is a MDSF of $V$. Thus, if all edge-weights are distinct, the MDSF
is unique,  and is
 obtained by connecting each non-minimal vertex to its directed
nearest neighbour. \\

For what follows, we consider a general type of partial ordering
of $\R^2$, denoted
$\stackrel{\theta,\phi}{\preccurlyeq}$, specified by the angles
%$0 \leq \theta < \theta+\phi \leq \pi$.
$\theta \in [0 ,2 \pi)$ and $\phi \in (0,\pi ] \cup \{2\pi\}$.
For $\mathbf{x} \in
\R^2$, let $C_{\theta,\phi}(\mathbf{x})$ be the closed cone with
vertex $\mathbf{x}$ and boundaries given by the rays from
$\mathbf{x}$ at angles $\theta$ and $\theta+\phi$, measuring
anticlockwise from the upwards vertical. The partial order is such
that, for $\mathbf{x}_1, \mathbf{x}_2 \in \R^2$,
\bea
 \mathbf{x}_1 \stackrel{\theta,\phi}{\preccurlyeq} \mathbf{x}_2
\textrm{ iff } \mathbf{x}_1 \in C_{\theta,\phi} (\mathbf{x}_2) .
\label{0719}
\eea
We shall use $\preccurlyeq^*$ as shorthand
for the special case $\stackrel{\pi/2,\pi/2}{\preccurlyeq}$,
which is of particular interest, as in~\cite{bhattroy2002}.
In this case
$u \postar v$ for $u=(u_1,u_2),v=(v_1,v_2) \in E$ if and only if
$u_1 \leq v_1$ and $u_2 \leq v_2$.
The
symbol $\preccurlyeq$ will denote a general partial order on
$\R^2$. 

We do not permit here the case $\phi=0$,  which
would almost surely give us a disconnected point set.
Nor do
%Note that if
 we allow $\pi < \phi < 2\pi$, since in this case 
the
directional relation (\ref{0719})
is not a partial order, since the transitivity property
(if $u \preccurlyeq v$ and $v \preccurlyeq w$ then $u \preccurlyeq w$)
fails for $\pi < \phi < 2\pi$. 
We shall, however, allow the
case $\phi=2\pi$ which leads to the standard nearest
neighbour (directed) graph.

The weight function is given by power-weighted Euclidean distance,
i.e., for $(u,v) \in E$ we assign weight $w(u,v) = \|u-v\|^\alpha $
to the edge $(u,v)$, where $\|\cdot\|$ denotes the Euclidean norm
on $\R^2$, and $\alpha >0$ is an arbitrary fixed parameter.
Thus, when $\alpha =1$ the weight of an edge is simply 
its Euclidean length.
Moreover, we shall assume that $V \subset \R^2$ is given by $V=
\mathcal{S}$ or $V= \SS^0 := \SS \cup \{ \0 \} $, where
$\0$ is the origin in $\R^2$ and
 $ \SS$ is generated in a {\em random} manner.
The random point set $\SS$ will usually be either the set of
points given by a homogeneous Poisson point process $\Po_n$ of
intensity $n$ on the unit square $(0,1]^2$, or a binomial point
process $\X_n$ consisting
of $n$ independent uniformly distributed points on $(0,1]^2$.

Note that in
this random setting,
%with probability one $V_0 = \{\0\}$ and
each point of $\SS$ almost surely
 has  a unique directed nearest neighbour, so
that $V$ has a unique MDSF, which does not
depend on the choice of $\alpha$.  Denote by $\LL^\alpha(\SS)$ the total
weight of all the edges in the MDSF on $\SS$, and let
 $\tilde \LL^\alpha(\SS):= \LL^\alpha(\SS)-\Exp[\LL^\alpha(\SS)]$, 
the centred total weight. 

Our first result presents laws of large numbers for the
 total edge weight for the
general partial order $\potp$ and general $0<\alpha<2$.
We state the result for $n$ points uniformly distributed
on $(0,1]^2$, but the proof carries through to other distributions
(see the start of Section \ref{seclln}).

\begin{theorem} \label{llnthm} 
Suppose $0 < \alpha < 2$.
Under the general partial
order $\potp$, with $0 \leq \theta < 2\pi$ and
$0 < \phi \leq \pi$ or $\phi = 2\pi$, it is the case that
\bea
n^{(\alpha/2)-1} \LL^\alpha ( \X_n ) \inL 
(2/\phi)^{\alpha/2} \Gamma( 1+\alpha/2), 
~~~ {\rm as }~ n\to
\infty.
\label{0728e}
\eea
Also, 
when the partial order is $\postar$, (\ref{0728e}) remains true with the addition
of the origin, i.e.~with $\X_n$
replaced by $\X_n^0$.
\end{theorem}
\rem In the special case $\alpha=1$, the limit in
(\ref{0728e}) is $\sqrt{ \pi / (2 \phi)}$. This limit is 1 when
$\phi=\pi/2$. Also, for $\phi=2\pi$ we have the standard
nearest neighbour (directed) graph (that is, every point is joined to its
nearest neighbour by a directed edge), and this limit is then $1/2$. 
This result (for $\alpha =1, \phi=2\pi$) is stated without
proof (and attributed to Miles \cite{Miles})
in \cite{avrambertsimas1993}, but we have 
%Results of this type date back at least to Miles \cite{Miles};
%see also \cite{penyuk2}, and page 101 of \cite{yukbook}. However,
%we have 
not previously seen the limiting constant 
derived explicitly, either in \cite{Miles} or anywhere else. \\ 

Our main result (Theorem \ref{mainth})
presents convergence in distribution 
for the case where the partial order is $\postar$;
the limiting distributions are of a different type in 
the three cases $\alpha=1$ (the same
situation as \cite{bhattroy2002}),  $0<\alpha < 1$, and
$\alpha >1$.  We define these 
limiting distributions
 in Theorem \ref{mainth},
in terms of distributional fixed-point equations.
These fixed-point equations are of the form
\bea
\label{0701a}
X \eqd \sum_{r=1}^k A_r X^{\{r\}} + B ,
\eea
where $k \in \N$,
$X^{\{r\}}$, $r=1,\ldots,k$, are
independent copies of the random variable $X$,
 and $(A_1,\ldots,A_k,B)$
is a random vector, independent of
 $( X^{\{1\}}, \ldots, X^{\{k\}})$,
satisfying the  conditions
\bea
\label{0701b}
\Exp \sum_{r=1}^k | A_r |^2 < 1 , ~~~~\Exp [B] =0, ~~~~\Exp [B^2] < \infty.
\eea
%for some $s$ with $0<s \leq 3$.  
 Theorem 3 of R\"osler \cite{rosler0} 
(proved using the contraction mapping theorem;
see also \cite{neinrusch,rosler})
says that if (\ref{0701b}) holds, 
there is a unique square-integrable distribution
with mean zero satisfying the fixed-point equation (\ref{0701a}), and
this will guarantee uniqueness of solutions to all
 the distributional
fixed-point equalities
 considered in the sequel.

Define the random variable $\tDone$,
to have the distribution that is the unique solution
to the distributional fixed-point
equation
\bea
\label{0628a}
 \tDone \eqd U \tDone^{\{1\}} +(1-U)
\tDone^{\{2\}} +U \log{U} +(1-U) \log(1-U) +U ,
\eea
where
$U$ is uniform on $(0,1)$ and independent of the other variables on the right.
We shall see later (in Propositions \ref{dconverge} and \ref{1020d})
 that $\Exp [\tDone] =0$ and $\Var[\tDone] = 2 -\pi^2/6$;
higher order moments are given recursively by eqn (\ref{0709c}).

For $\alpha >1$, let 
$\tDalph$ denote a random variable with
 distribution characterized
by the fixed-point equation
\bea
\label{0628b}
 \tDalph \eqd U^\alpha \tDalph^{\{1\}} +(1-U)^\alpha
\tDalph^{\{2\}} +\frac{\alpha}{\alpha-1}U^\alpha +
\frac{1}{\alpha-1}(1-U)^\alpha-\frac{1}{\alpha-1},
\eea
where again
$U$ is uniform on
$(0,1)$ and independent of the other variables on the right.
Also for $\alpha>1$, let $\tFalph$ denote a
 random variable with distribution characterized
by the fixed-point equation
\bea
\label{0628d}
 \tFalph \eqd U^\alpha \tFalph +(1-U)^\alpha
\tDalph +\frac{U^\alpha}{\alpha(\alpha-1)}
+ \frac{(1-U)^\alpha}{\alpha-1} -\frac{1}{\alpha(\alpha-1)}  ,
\eea
where $U$ is uniform on $(0,1)$, $\tDalph$
has the distribution given by (\ref{0628b}), and
the $U$, $\tDalph$ and $\tFalph$ on the right are independent.
In Section \ref{secdlt} we shall see that the
random variables $\tDalph$, $\tFalph$ for
$\alpha>1$ arise as centred versions of random 
variables (denoted $\Dalph$, $\Falph$ respectively) satisfying
somewhat simpler fixed point equations.
Thus $\tDalph$ and $\tFalph$ both have mean zero;
their variances are given by eqns (\ref{0701g}) and (\ref{0701i})
below.

Let $\NN(0,s^2)$ denote the normal distribution 
with mean zero and variance $s^2$.

\begin{theorem} \label{mainth}
Suppose the weight exponent is $\alpha > 0$
%$1 \leq \alpha \leq 2$
and the partial order is $\postar$.
There exist constants $0 < t_\alpha^2 \leq s_\alpha^2$
such that,
for normal random variables 
$\Yalph \sim \mathcal{N}(0,s_\alpha^2)$ and  $\Walph
\sim \mathcal{N}(0,t_\alpha^2)$:
%\begin{itemize}
% \item[(i)]

(i)  As $n
\to \infty$,
\begin{eqnarray} n^{(\alpha-1)/2} \tLalph ( \Po^0_n )
\tod \Yalph & {\rm and } & n^{(\alpha-1)/2} \tLalph ( \X^0_n )
\tod \Walph ~~~  (0 < \alpha < 1);~~
\label{0727a}  \\
\tLone ( \Po^0 _n ) \tod \tDone^{\{1\}}
 + \tDone^{\{2\}} +Y_1   &  {\rm and } &
\tLone ( \X^0 _n ) \tod \tDone^{\{1\}}
 + \tDone^{\{2\}} +W_1 ; ~~~~ 
\label{0727b} \\
\tLalph ( \Po^0 _n ) \tod \tDalph^{\{1\}}
 + \tDalph^{\{2\}}   &  {\rm and } &
\tLalph ( \X^0 _n ) \tod \tDalph^{\{1\}}
 + \tDalph^{\{2\}} ~~~ (\alpha > 1).~~
\label{0727c}
\end{eqnarray}
Here all the random variables in the limits are independent, and
$\tDalph^{\{i\}}$, $i=1,2$ are independent copies of the random
variable $\tDalph$ defined at (\ref{0628a}) for $\alpha=1$
and (\ref{0628b}) for $\alpha >1$.

(ii)  As $n \to \infty$,
\bea
n^{(\alpha-1)/2} \tLalph ( \Po_n )
\tod \Yalph & {\rm and } & n^{(\alpha-1)/2} \tLalph ( \X_n )
\tod \Walph ~~~ (0 < \alpha < 1); ~~ 
\label{0727d} \\
\tLone ( \Po_n ) \tod \tDone^{\{1\}}
 + \tDone^{\{2\}} +Y_1  & {\rm and } &
\tLone ( \X_n ) \tod \tDone^{\{1\}}
 + \tDone^{\{2\}} +W_1 ;~~~~ 
\label{0727e} \\
\tLalph ( \Po_n ) \tod \tFalph^{\{1\}}
 + \tFalph^{\{2\}}  & {\rm and } &
\tLalph ( \X_n ) \tod \tFalph^{\{1\}}
 + \tFalph^{\{2\}} ~~~ (\alpha > 1) ~~
\label{0727f} .
\eea
 Here all the random variables in the limits
are independent, and $\tDone^{\{i\}}$, $i=1,2$,
are independent copies of $\tDone$ with
distribution defined at (\ref{0628a}),
and for $\alpha>1$, $\tFalph^{\{i\}}$, $i=1,2$, are independent
copies of 
%a random variable 
$\tFalph$ with distribution
defined at 
(\ref{0628d}).
\end{theorem}

\rems
The normal random variables $\Yalph$ or $\Walph$
arise from the edges
away from the boundary (see Section \ref{ltot}).
The non-normal variables (the $\tD$s and $\tF$s)
arise from the edges very close to the boundary, where the MDSF
is asymptotically close to the `directed linear forest'
discussed in Section \ref{secdlt}.

Theorem \ref{mainth}
 indicates a phase transition
in the character of the limit law as $\alpha$ increases.
The normal contribution (from the points away from the boundary)
dominates for $0<\alpha < 1$, while the boundary contributions dominate
for $\alpha > 1$. In the critical case $\alpha=1$, neither effect dominates
and both terms contribute significantly to the asymptotic behaviour. 

Noteworthy in the case $\alpha =1$
is the fact that by (\ref{0727b}) and (\ref{0727e}),
the limiting distribution is
the same for $\tLone(\Po_n) $ as for $\tLone(\Po_n^0)$,
and the same for $\tLone(\X_n)$ as for $\tLone(\X_n^0)$.
Note, however, that the difference $\tLone(\Po_n)-\tLone(\Po_n^0)$
is the (centred) total length of edges incident to the origin,
which is not negligible, but itself converges in distribution
(see \cite{rooted}) to a non-degenerate
random variable, namely a centred generalized Dickman random
variable with parameter $2$ (see (\ref{0629g}) below).
As an extension of Theorem \ref{mainth},
it should be possible to show that the joint distribution
of $(\tLone(\Po_n),\tLone(\Po_n^0))$ converges to that
of two coupled random variables, both having the distribution
of $\tDone$, whose difference has the centred generalized Dickman 
distribution with parameter 2. Likewise for the joint distribution
of $(\tLone(\X_n),\tLone(\X_n^0))$. 

Of particular interest
is the distribution of the variable $\tDone$ 
appearing in Theorem \ref{mainth}. 
 In Section \ref{subsecalphone},
we give a plot (Figure \ref{pdffig}) of the probability density function
of this distribution, estimated by simulation. Also,
we can use the fixed-point equation (\ref{0628a}) 
 to calculate the moments of $\tDone$ recursively. Writing
\[ f(U) := U \log U + (1-U) \log (1-U) + U, \]
and setting $m_k := \Exp[\tDone^k]$, 
we obtain
\bea
\label{0709c}
 m_k = \Exp[ (f(U))^k ] + \sum_{i=2}^k {k \choose i} \sum_{j=0}^i
{i \choose j}
\Exp [ (f(U))^{k-i} U^j (1-U)^{i-j} ]
m_j m_{i=j} .
\eea
The fact that $m_1=0$ simplifies things a little, and we can
rewrite this as
\bean
m_k = \Exp[ (f(U))^k ] + \sum_{i=1}^k {k \choose i} \left[
 m_i \Exp[ (f(U))^{k-i}(U^i+(1-U)^i)] ~^{~^{~^{~^{~^{~^{~^{~^{~}}}}}}}} \right.
~~~~
\\
\left. +
\sum_{j=2}^{i-2}
{i \choose j}
\Exp [ (f(U))^{k-i} U^j (1-U)^{i-j} ]
m_j m_{i-j} \right].
\eean
So, for example, when $k=3$ we obtain $m_3 \approx 0.15411$,
which shows $\tDone$ is not Gaussian and 
is consistent with the skewness of the plot
in Figure \ref{pdffig}.  \\

The remainder of this paper is organized as follows.
After discussion of the DLT in Section \ref{secdlt},
in Section \ref{secgeneral} we present general
limit theorems in geometric probability, which
we shall use in obtaining our main results for
the MDST.
 Theorem \ref{llnthm} is proved in Section \ref{seclln}
(this proof does not use the results of Section \ref{secdlt}).
The proof of Theorem \ref{mainth} is prepared in
Sections \ref{ltot} and \ref{bdry}, and completed in
Section \ref{totallength}. In these proofs, we repeatedly
use {\em Slutsky's theorem}
 (see e.g.~\cite{penbook})
 which says that
if $X_n \to X$ in distribution and $Y_n \to 0$ in
probability, then $X_n +Y_n \to X $ in distribution.

\section{The directed linear forest and tree}
\label{secdlt}

The directed linear forest (DLF)  and  directed linear tree (DLT)
are for us a tool
for the analysis of the limiting behaviour of
the contribution to the total weight of the random
MDSF/MDST from edges near the boundary of the unit
square.
 In the present section we derive
the properties of the DLF that we need (in particular, Theorem
\ref{dltthm}); subsequently,
 in Theorem \ref{thmbdry}, we shall see that the total 
 weight of edges from the points near
the boundaries, as $n \to \infty$,
 converges in distribution to the limit of the
total weight of the DLF.

    The  DLT is also of some intrinsic interest.
It is a one-dimensional directed analogue of the
so-called `on-line nearest neighbour graph', which is of interest
in the study of networks such as the world wide web 
(see, e.g.~\cite{boll};
and~\cite{mdp} for more on the on-line nearest neighbour graph).
Moreover, it is constructed via a fragmentation process
similar to those seen in, for example, \cite{bertoin}; the
tree provides a historical representation of the fragmentation
process.

For any finite sequence $\mathcal{T}_m = ( x_1,x_2,\ldots,x_m ) \in  (0,1]^m$,
we construct the directed linear forest (DLF) as follows.
We start with the unit interval $(0,1]$
and insert the points $x_i$ in order, one at a time, starting with $i=1$.
At the insertion of each point, we join the new point to its nearest neighbour
among those points already present that lie to the {\em left} of the
 point (provided that such a point exists). In other words,
for each point $x_i$, $i \geq 2$, we join $x_i$ by a directed edge to the point
$\max \{ x_j : 1 \leq  j < i, \; x_j < x_i \}$. If 
$ \{ x_j : 1 \leq  j < i, \; x_j < x_i \}$ is empty, we do not add
any directed edge from $x_i$.
In this way we construct a `directed linear forest', which we denote by
 $\mathrm{DLF} \left( \mathcal{T}_m \right)$.
 We denote the total
weight (under weight function with exponent $\alpha$) of 
 $\mathrm{DLF} \left( \mathcal{T}_m \right)$ by
$D^\alpha ( \mathcal{T}_m )$, that is, we set
$$
 D^\alpha \left( \mathcal{T}_m \right):= \sum_{i=2}^m 
%(x_i - \max \{x_j: 1 \leq j < i, x_j < x_i \})^\alpha{\bf 1}\{x_i \neq
%\min \{ x_j: 1 \leq j \leq i \} \}.
(x_i - \max \{ x_j: 1 \leq j < i, x_j < x_i \} )^\alpha{\bf 1}\{
\min \{x_j: 1 \leq j < i \} < x_i \}.
$$

Further, given $\mathcal{T}_m$, let
$\mathcal{T}^0_m$ be the sequence
 $( x_0,x_1,\ldots,x_m )$
where the initial term is $x_0 := 0$. Then the DLF on
$\mathcal{T}^0_m$ is constructed in the same way, where now
for each $i \geq 1$, we join $x_i$ by an edge to the point
$\max \{ x_j : 0 \leq  j < i, \; x_j < x_i \}$.  But now we see that $x_1$ will
always be joined to $x_0=0$, and $x_2$ will be joined either to $x_1$
(if $x_2 > x_1$) or to $x_0$, and so on.  In this way we construct
 a `directed linear tree' (DLT) on vertex set
 %${\cal T}_m^0$ 
 $\{ x_0,x_1,\ldots,x_m \}$
with $m$ edges.
Denote the total weight of this tree with weight exponent $\alpha $
by $D^\alpha (\mathcal{T}_m^0)$; that is, set
$$
 D^\alpha \left( \mathcal{T}_m^0 \right):= \sum_{i=1}^m 
(x_i- \max\{ x_j: 0 \leq j < i, x_j < x_i \})^\alpha.
$$
We shall be  mainly interested in the case where $\mathcal{T}_m$
is a random vector in $(0,1]^m$.
In this case, 
set
$\tD^\alpha \left( \mathcal{T}_m \right) := D^\alpha \left( \mathcal{T}_m
\right) - \Exp \left[ D^\alpha \left( \mathcal{T}_m \right)
\right]$ the centred total weight of the DLF,
and
$\tD^\alpha \left( \mathcal{T}_m^0 \right) = D^\alpha \left( \mathcal{T}_m^0
\right) - \Exp \left[ D^\alpha \left( \mathcal{T}_m^0 \right)
\right]$ the centred total weight of the DLT.
% In the case where the components
% of $\mathcal{T}_m$ are independent uniform random variables on
%$(0,1]$, we write $\mathcal{T}_m=\mathcal{U}_m$.

We take $\mathcal{T}_m$ to be a vector of
uniform  variables.
Let $(X_1,X_2,X_3,\ldots)$ be a sequence of independent
uniformly distributed 
random variables   in $(0,1]$, and for $m \in \N$ set 
 $\UU_m := ( X_1, X_2, \ldots, X_m )$.
% where all the $X_i$ are independent and uniformly distributed on $(0,1]$.
We consider $D^\alpha(\UU_m)$ and $D^\alpha(\UU_m^0)$. For these
variables, we establish asymptotic behaviour of the mean value
in Propositions \ref{dltmoms} and \ref{dlfmoms},
along with the following convergence results, which
are the principal results of this section.

For $\alpha >1$, let 
$\Dalph$ denote a random variable with
 distribution characterized
by the fixed-point equation
\bea
\label{0628bb}
 \Dalph \eqd U^\alpha \Dalph^{\{1\}} +(1-U)^\alpha
\Dalph^{\{2\}} +U^\alpha ,
\eea
where 
$U$ is uniform on
$(0,1)$ and independent of the other variables on the right.
Also for $\alpha>1$, let $\Falph$ denote a
 random variable with distribution characterized
by the fixed-point equation
\bea
\label{0628dd}
 \Falph \eqd U^\alpha \Falph +(1-U)^\alpha
\Dalph   ,
\eea
where $U$ is uniform on $(0,1)$, $\Dalph$
has the distribution given by (\ref{0628bb}), and
the $U$, $\Dalph$ and $\Falph$ on the right are independent.
The corresponding centred
random variables  $\tDalph := \Dalph -\Exp[\Dalph]$ and
$\tFalph := \Falph - \Exp [ \Falph]$
satisfy the fixed-point equations (\ref{0628b}) and
(\ref{0628d}) respectively. The solutions to \eq{0628b}
and \eq{0628d} are unique by the criterion given at
\eq{0701b}, and hence the solutions to \eq{0628bb} and \eq{0628dd}
are also unique.

\begin{theorem} \label{dltthm}
\begin{itemize}
\item[(i)]
As $m\to \infty$ we have
 $\tD^1(\UU_m^0) \inLL \tDone$ and $\tD^1(\UU_m) \inLL \tFone$
where $\tDone$ has the distribution 
given by the fixed-point equation (\ref{0628a}),
and $\tFone$  has the same distribution 
as $\tDone$. 
Also, the variance of $\tDone$ (and hence also of
$\tFone$) is $2 -\pi^2/6 \approx 0.355066$.
Finally, $\Cov(\tDone,\tFone) = (7/4)-\pi^2/6 \approx
0.105066$.
\item[(ii)] For $\alpha>1$, as $m\to \infty$ we have
$D^\alpha(\UU_m^0) \to \Dalph$, almost surely and in $L^2$,
and
 $D^\alpha(\UU_m) \inLL \Falph$, almost surely and in $L^2$,
 where the distributions of
$\Dalph$, $\Falph$ are given by the fixed-point equations 
(\ref{0628bb})
and (\ref{0628dd}) respectively.
Also, $\Exp[\Dalph]= (\alpha-1)^{-1}$ and
 $\Exp[\Falph]= (\alpha(\alpha-1))^{-1}$, 
while
 $\Var(\Dalph)$ and $\Var(\Falph)$ are
given by (\ref{0701g}) and (\ref{0701i}) respectively.
\end{itemize}
\end{theorem}
\proof 
Part (i) follows from
 Propositions \ref{dconverge},
 \ref{1020d} 
and
 \ref{ffixed}
below.
Part (ii) follows from
Propositions \ref{1119f} and  \ref{1125a} below.
We prove these results in the following sections.
 $\square$ \\

An interesting property of the DLT, which we use
in establishing fixed-point equations for limit distributions,
is its {\em self-similarity} (scaling property). In terms of the total
weight, this says that for any $t \in (0,1)$, if
$Y_1,\ldots,Y_n$ are independent and uniformly distributed on $(0,t]$,
then the distribution
of $D^\alpha(Y_1,\ldots,Y_n)$ is the same as that of
$t^\alpha D^\alpha(X_1,\ldots,X_n)$.

\subsection{The mean total weight of the DLF and DLT} \label{dltproperties}

First we consider the rooted case, i.e.~the DLT on $\UU^0_m$.
For $m=1,2,3,\ldots$
denote by $Z_m$ the random variable given by the gain in length
of the tree on the addition of one point ($X_m$) to an existing $m-1$
points in the DLT on a sequence of uniform random variables
$\UU^0_{m-1}$, i.e.~with the conventions
$D^1 (\UU_0^0)=0$ and 
$X_0 =0$,
we set
\bea
Z_m := D^1(\UU_m^0)-D^1(\UU_{m-1}^0) = 
X_m- \max\{ X_j: 0 \leq j < m, X_j < X_m \}.
\label{0728j}
\eea
Thus, with weight exponent $\alpha$, the $m$th edge to be added
has weight $Z_m^\alpha$.

\begin{lemma} \label{Zmdist} %\begin{itemize}
% \item[(i)]
(i)  $Z_m$ has distribution function $F_m$ given by
 $F_m(t)=0$ for $t<0$, $F_m(t)=1$ for $t>1$, and
$ F_m(t) = 1-(1-t)^m $  for $ 0 \leq t \leq 1$.
% \item[(ii)] 

(ii)
For $\alpha >0$, 
 $Z_m^\alpha$ has expectation and variance
\begin{equation} \label{zexpvar}
\Exp[Z_m^\alpha]=\frac{m! \Gamma(1+\alpha)}{\Gamma(1+\alpha+m)}
, ~~~ \Var[Z_m^\alpha] = \frac{m! \Gamma(1+2\alpha)}{\Gamma(1+2\alpha+m)}
- \left( \frac{m! \Gamma(1+\alpha)}{\Gamma(1+\alpha+m)}\right)^2.
%\frac{m}{(m+1)^2(m+2)}
\end{equation} 
In particular,
\bea
\label{0630b}
 \Exp[Z_m]=
\frac{1}{m+1} ; ~~ \Var[Z_m] = \frac{m}{(m+1)^2(m+2)}. 
\eea

(iii)
%\item[(iii)]
For $\alpha >0$,
as $m \to \infty$ we have
 \begin{equation} \label{asymZm}
\Exp[Z_m^\alpha]
\sim \Gamma (\alpha+1) m^{-\alpha}, \;
\Var[Z_m^\alpha] \sim \left( \Gamma(2\alpha+1) - (\Gamma(\alpha+1))^2 \right)
m^{-2\alpha} . \end{equation}

(iv)
%\item[(iv)] 
As $m \to \infty$,
$mZ_m$ converges in distribution, 
 to an exponential
 with parameter $1$.
%\end{itemize}
\end{lemma}
\proof For $0 \leq t \leq 1$ we have
\[ \Pr[Z_m > t]  =  \Pr [ 
X_{m}> t 
\textrm{ and none of }
X_1,\ldots,X_{m-1} \textrm{ lies in } (X_{m}-t,X_{m}) 
]  =  (1-t)^{m} , \]
%Then by setting $s=t^{1/\alpha}$ we obtain
%which gives us (\ref{0629e}), proving (i).
and (i) follows.
We then obtain (ii)
since for any $\alpha >0$ and for $k=1,2$,
\[
 \Exp[Z_m^{k\alpha}] = \int_0^1 \Pr [ Z_m > t^{1/(k \alpha)} ]
 \ud t = \int_0^1
( 1-t^{1/k\alpha})^m \ud t =
\frac{m! \Gamma(1+k\alpha)}{\Gamma(1+k\alpha+m)}.
\]
Then (iii) follows by Stirling's formula, which yields
$$
\Exp[Z_m^{k\alpha}] = \Gamma(1+k\alpha)
 m^{-k\alpha} (1+O(m^{-1})).
$$
For (iv), we have from (i) that, for $t \in [0,\infty)$,
and $m$ large enough so that $(t/m) \leq 1$,
 \[ \Pr[mZ_m \leq t]
= F_m \left( \frac{t}{m} \right) = 1- \left( 1 -
\frac{t}{m} \right)^m \to  1-e^{-t} , \textrm{ as } m \to \infty.
\]
But $1-e^{-t}$, $t \geq 0$
 is the exponential distribution function with
parameter 1. $\square$ \\

The following result gives the asymptotic 
behaviour of the expected total weight of the DLT.
Let $\gamma$ denote Euler's constant, so that
\bea
\label{0708b}
 \left( \sum_{i=1}^k  \frac{1}{i} \right) - \log k = \gamma +O(k^{-1}).
\eea

\begin{proposition} \label{dltmoms}
As $m \to \infty$ the expected total weight of the DLT under $\alpha$-power 
weighting on
$\mathcal{U}^0_m$ satisfies
\begin{eqnarray} 
\Exp[ D^\alpha(\UU_m^0) ] & \sim &
\frac{\Gamma(\alpha+1)}{1-\alpha} m^{1-\alpha} ~~~  (0<\alpha<1); 
\label{0720a}
\\
\Exp[ D^1(\UU_m^0) ] - \log{m}   & \to & \gamma-1; 
\label{0720b}\\
\Exp[ D^\alpha(\UU_m^0) ] & = & \frac{1}{\alpha-1}
+O(m^{1-\alpha}) ~~~ (\alpha>1). \label{0720c}
 \end{eqnarray} 
\end{proposition} 
\textbf{Proof.}
We have \[ \label{eqn100a} \Exp[D^\alpha(\UU_m^0)] = \sum_{i=1}^m
\left( \Exp [ D^\alpha(\UU_i^0) ] - \Exp[D^\alpha(\UU_{i-1}^0)]
\right) = \sum_{i=1}^m \Exp[Z^\alpha_i] .
\]
 In the case  where $\alpha=1$, $E[Z_i] = (i+1)^{-1}$ by (\ref{0630b}),
and  (\ref{0720b}) follows
by (\ref{0708b}).
 For general $\alpha>0$, $\alpha \neq 1$,
from (\ref{zexpvar}) we have that \begin{equation} \label{1117a}
\Exp[D^\alpha(\UU_m^0)] = \Gamma(1+\alpha) \sum_{i=1}^m
\frac{\Gamma(i+1)}{\Gamma(1+\alpha+i)} = \frac{1}{\alpha-1} -
\frac{\Gamma(1+\alpha) \Gamma(m+2)}{(\alpha-1) \Gamma(m+1+\alpha)}
. \end{equation} By Stirling's formula, the last term satisfies
\begin{equation}
 \label{1117b}
 -\frac{\Gamma(1+\alpha) 
\Gamma(m+2)}{(\alpha-1)
\Gamma(m+1+\alpha)} = -\frac{\Gamma(1+\alpha)}{\alpha-1}
m^{1-\alpha} (1+O(m^{-1})),
 \end{equation} 
which tends to zero as $m \to \infty$ for $\alpha >1$,
to give us (\ref{0720c}). For $\alpha <1$, we have (\ref{0720a}) from
(\ref{1117a}) and (\ref{1117b}).  $\square$ \\

Now consider the unrooted case, i.e., the directed linear
forest.  For $\mathcal{U}_m$  as above
 the total weight of the DLF is denoted $D^\alpha ( \UU_m)$,
and the centred total weight is $\tD^\alpha (\UU_m) :=
D^\alpha (\UU_m)- \Exp [D^\alpha (\UU_m)]$. We
then see that
\bea
 D^\alpha ( \UU_m^0 ) = D^\alpha ( \UU_m) + \LL^\alpha_0( \UU_m^0) , 
\label{0714b}
\eea
 where $\LL^\alpha_0( \UU_m^0)$ is the total weight of
edges incident to 0 in the DLT on $\UU_m^0$.

The following lemma says that $\LL^\alpha_0( \UU_m^0)$
converges to a random variable that has
the generalized Dickman distribution with
parameter $1/\alpha$ (see \cite{rooted}), that is,
the distribution of a random variable $X$ which satisfies
the distributional fixed-point equation
\bea
\label{0629g}
 X \eqd U^{\alpha} (1+X) ,\eea
where $U$ is uniform on $(0,1)$ and independent of the
$X$ on the right. We recall
 from Proposition 3 of \cite{rooted} that if
$X$ satisfies (\ref{0629g}) then
\bea
\label{0630j}
\Exp[X] = 1/\alpha, \textrm{ and } \Exp[ X^2] = (\alpha+2)/(2\alpha^2).
\eea
\begin{lemma} \label{dickmanlem} Let $\alpha>0$.
There is a random variable  $\LL^\alpha_0$  with
the generalized Dickman distribution with parameter
$1/\alpha$, such that
as $m \to \infty$,
we have that
$ \LL^\alpha_0( \UU_m^0) \to 
\LL^\alpha_0 $, almost surely and in $L^2$.
\end{lemma} \textbf{Proof.} 
Let $\delta_D(\UU_m^0)$ denote the degree of the origin in the
directed linear tree on $\UU^0_m$, so that
$\delta_D(\UU_m^0)$ is the number of lower records in 
the sequence $(X_1,\ldots, X_m)$. Then
\bea
\label{0629f}
 \LL_0^\alpha (\UU_m^0) = U_1^\alpha +(U_1U_2)^\alpha
 + \cdots + ( U_1 \cdots U_{\delta_D(\UU_m^0)} )^\alpha ,
\eea
where $(U_1, U_2, \ldots)$ is a certain
sequence of independent uniform random variables on $(0,1)$,
namely the ratios between successive lower records of
the sequence $(X_n)$.
The sum
$
% \LL_0^\alpha := 
U_1^\alpha +(U_1U_2)^\alpha + (U_1U_2U_3)^\alpha + \cdots  
$
 has nonnegative terms and finite expectation, so it
converges almost surely to a limit which we denote
 $\LL_0^\alpha$.
 Then $\LL_0^\alpha$
 has the generalized Dickman distribution
 with parameter $1/\alpha$ (see Proposition 2 of \cite{rooted}).

Since $\delta_D(\UU_m^0)$ tends to infinity almost surely as $m \to \infty$,
we have $\LL_0^\alpha(\UU_m^0) \to \LL_0^\alpha$ almost surely.
Also, $\Exp [ (\LL_0^\alpha)^2 ]
< \infty$, by (\ref{0630j}),
and $(\LL_0^\alpha -\LL_0^\alpha (\UU_m^0))^2 \leq (\LL_0^\alpha)^2$
for all $m$. Thus $\Exp [ (\LL_0^\alpha (\UU_m^0) -\LL_0^\alpha )^2 ] \to 0$
by the dominated convergence theorem, and so we have the $L^2$
convergence as well. $\square$

\begin{proposition} \label{dlfmoms}
As $m \to \infty$ the expected total weight of the DLF under
 $\alpha$-power weighting on
$\mathcal{U}_m$ satisfies
\begin{eqnarray} 
\Exp[ D^\alpha(\UU_m) ] & \sim &
\frac{\Gamma(\alpha+1)}{1-\alpha} m^{1-\alpha} ~~~  (0<\alpha<1); 
\label{0720d}
\\
\Exp[ D^1(\UU_m) ] -\log{m}   & \to & \gamma -2
% ~~~  
%(\alpha=1)
; \label{0720e}
\\
\Exp[ D^\alpha(\UU_m) ] & \to & \frac{1}{\alpha(\alpha-1)} ~~~
(\alpha>1). \label{0720f}
 \end{eqnarray}
\end{proposition}
\proof
By (\ref{0714b}) we have
 $  \Exp[
D^\alpha(\UU_m) ] = \Exp[ D^\alpha(\UU_m^0) ] -
 \Exp[ \LL_0^\alpha(\UU_m^0) ]. 
$ 
By
  Lemma \ref{dickmanlem} and
 (\ref{0630j}),
\[ \Exp [ \LL_0^\alpha ( \UU_m^0 )] \longrightarrow \Exp [ \LL_0^\alpha ] = 1/\alpha.\]
We then obtain (\ref{0720d}), (\ref{0720e}) and (\ref{0720f})
 from Proposition \ref{dltmoms}. $\square$

\subsection{Orthogonal increments for $\alpha =1$}
In this section we shall show (in Lemma \ref{covz}) 
that when $\alpha =1$, the variables $Z_i, i \geq 1$ are
mutually orthogonal, in the sense of having zero covariances,
which will be used later on to establish convergence of the
(centred) total length of the DLT.  
To prove this, we first need further notation.

Given  $X_1,\ldots,X_m$, let us denote the   order statistics
of $X_1,\ldots,X_m$, taken in increasing order,
 as  
\linebreak
$ X_{(1)}^m, X_{(2)}^m, \ldots, X_{(m)}^m $. Thus 
$ ( X_{(1)}^m, X_{(2)}^m, \ldots,
X_{(m)}^m ) $ is a nondecreasing sequence, forming
 a permutation of the original
$(X_1, \ldots, X_m)$.
Denote the existing
$m+1$ intervals
between points by
$I_j^{m}:= \left( X_{(j-1)}^m, X_{(j)}^m \right)$ for $j=1,2,\ldots,m+1$,
where we set $X_{(0)}^m:=0$ and $X_{(m+1)}^m:=1$. Let the widths of
these intervals (the spacings) be
\[ S_j^m := \left| I_j^m \right| = X_{(j)}^m-X_{(j-1)}^m , \]
 for $1 \leq j \leq
m+1$. Then $0 \leq S_j^m < 1$ for $1 \leq j \leq m+1$, and
$\sum_{j=1}^{m+1} S_j^m = 1$. That is, the vector $\left(
S_1^m,S_2^m, \ldots,S_{m+1}^m \right)$ belongs to the
$m$-dimensional simplex, $\Delta_m$. 
Note that only $m$ of the $S_j^m$ are
required to specify the vector.

We can arrange the spacings themselves  $(S_j^m, 1 \leq j \leq m+1)$ into
increasing order to give $S_{(1)}^m,$ $S_{(2)}^m ,$ $\ldots,$
$S_{(m+1)}^m$. Then let $\mathcal{F}_S^m$ denote the sigma field
generated by these ordered spacings, so that
\bea
\label{0629d}
 \mathcal{F}_S^m = \sigma \left(
S_{(1)}^m,\ldots,S_{(m+1)}^m \right) .
\eea
The following interpretation of $\mathcal{F}_S^m $ 
may be helpful. The set $(0,1)\setminus\{X_1,\ldots,X_m\}$
consists almost surely of $m+1$ connected components (`fragments')
of total length 1, and $\mathcal{F}_S^m $ is the 
$\sigma$-field generated by the collection of lengths
of these fragments, ignoring the order in which they appear.

By definition, the value of $Z_m$ must be one of the
(ordered) spacings $S^m_{(1)}, \ldots, S^m_{(m+1)}$.
The next result says that, given the values of these 
spacings, each of the possible values for $Z_m$
are equally likely.

\begin{lemma} \label{0222q}
For 
$m \geq 1$ we have
\bea
\label{0526a}
\Pr \left[ \left. Z_m = S^m_{(i)} \right| \FF_S^m \right] = \frac{1}{m+1}
~~~{\rm a.s.},
\textrm{ for $i=1,\ldots,m+1$}.
\eea
Hence,
\bea
\Exp \left[ Z_m \left|
\mathcal{F}_S^m \right. \right] = \frac{1}{m+1} \sum_{i=1}^{m+1} 
S_{(i)}^m 
= \frac{1}{m+1} .
\label{0526b}
\eea
\end{lemma}
\textbf{Proof.} First we note that $\left( X_{(1)}^m, \ldots,
 X_{(m)}^m \right)$
is uniformly distributed over
 \[ \left\{ ( x_1, \ldots, x_m ) : 0 \leq x_1 \leq x_2
\leq \ldots \leq x_m \leq 1 \right\} . \] Now \[ \left(
\begin{array}{c} S_1^m \\ S_2^m \\ S_3^m \\ \vdots \\ S_m^m 
\end{array} \right)
= \left( \begin{array}{cccccc} 1 & 0 & 0 & \ldots & 0 & 0 \\
-1 & 1 & 0 & \ldots & 0 & 0 \\
0 & -1 & 1 & \ldots & 0 & 0 \\
\vdots & \vdots & \vdots & \ddots & \vdots & \vdots \\
0 & 0 & 0 & \ldots & -1 & 1 \end{array} \right) \left(
\begin{array}{c} X_{(1)}^m \\ X_{(2)}^m \\ X_{(3)}^m \\ \vdots
\\ X_{(m)}^m \end{array} \right) . \]
The $m$ by $m$ matrix here has determinant 1. Hence $\left( S_1^m,
\ldots , S_m^m \right)$ is uniform over \[ \left\{ \left( x_1,
\ldots, x_m \right) : \sum_{j=1}^m x_j \leq 1; x_j \geq 0 , \;
\forall \; 1 \leq j \leq m \right\}.\] Then $\left(
S_1^m, \ldots , S_{m+1}^m \right)$ is uniform over the
$m$-dimensional simplex $\Delta_m$.
In particular, the $S_j^m$ are exchangeable.
Thus given $S_{(1)}^m, \ldots, S_{(m+1)}^m$, i.e.~$\mathcal{F}_S^m$, the
actual values of $S_1^m, \ldots, S_{m+1}^m$ are equally
likely to be any permutation of $S_{(1)}^m, \ldots, S_{(m+1)}^m$, and
given $S_1^m, \ldots, S_{m+1}^m$
the value of $Z_m$ is equally likely to be any of $S_{1}^m,
\ldots, S_{m}^m$ (but cannot be $S_{m+1}^m$).

Hence, given $S^m_{(1)},\ldots,S^m_{(m+1)}$ the probability
that $Z_m = S^m_{(i)}$ is $(1/m) \times m/(m+1) = 1/(m+1)$,
i.e.~we have (\ref{0526a}), and then 
 (\ref{0526b})
follows since $\sum_{j=1}^{m+1} S_{(j)}^m =1$.
$\square$

\begin{lemma}
\label{0526c}
Let $1 \leq m < \ell$.
Given $\mathcal{F}_S^m$,  
 $Z_\ell$ and $Z_m$ are
conditionally independent. \end{lemma}
\textbf{Proof.} Given
$\mathcal{F}_S^m$, we have $S_{(1)}^m,\ldots,S_{(m+1)}^m$,
and by (\ref{0526a}),
the (conditional)
distribution of $Z_m$ is uniform on
$\{S_{(1)}^m,\ldots,S_{(m+1)}^m\}$.
The conditional distribution of $Z_\ell$, $\ell > m$,
 given $\FF_S^m$, depends only on $S_{(1)}^m,\ldots,S_{(m+1)}^m$ and
 not which
one of them $Z_m$ happens to be. Hence $Z_m$ and $Z_\ell$ are
conditionally independent. $\square$

\begin{lemma} \label{covz} 
For $1 \leq m < \ell$,
%$\alpha=1$, 
the random variables $Z_m,Z_\ell$
satisfy
$\Cov \left[ Z_m,Z_\ell \right] = 0$.
\end{lemma}
\textbf{Proof.} From Lemmas  \ref{0526c} and 
\ref{0222q},
 \bean
 \Exp \left[ Z_m Z_\ell |
\mathcal{F}_S^m \right] & = & \Exp \left[ Z_m | \mathcal{F}_S^m
\right] \Exp \left[ Z_\ell | \mathcal{F}_S^m \right]  =
\frac{1}{m+1} \Exp \left[ Z_\ell | \mathcal{F}_S^m \right], 
\eean
and by taking expectations we obtain
\bean
\Exp \left[ Z_m Z_\ell \right] & = & \frac{1}{m+1} \Exp
\left[ Z_\ell \right] = \frac{1}{m+1} \cdot \frac{1}{\ell+1} = \Exp [ Z_m ]
\cdot \Exp [Z_\ell].
\eean
Hence the covariance of $Z_m$ and $Z_\ell$ is zero. $\square$ \\

\rems
(i) Calculations yield, for example, that
$\Exp[D^1(\UU_1^0)]=\Exp[Z_1]=1/2$, $\Exp[D^1(\UU_2^0)]=5/6$, and
$\mathrm{Var}[Z_1]=1/12$, $\mathrm{Var}[Z_2]=1/18$,
$\mathrm{Var}[D^1(\UU_2^0)]=5/36$. \\

(ii) The orthogonality structure of the $Z_m^\alpha$ is unique to the
$\alpha=1$ case. For example, it can be shown that, for $\alpha>0$, 
\[
\Exp[Z_1^\alpha] \Exp[Z_2^\alpha] =
 \frac{2}{(1+\alpha)^2(2+\alpha)} , \textrm{
and } \Exp[Z_1^\alpha Z_2^\alpha] = \frac{1}{2(1+\alpha)^2 } \left( 1+ \frac
{2 \Gamma (\alpha +2)^2}{\Gamma (2 \alpha +3)} \right).
 \] 
Then 
\[
\mathrm{Cov}[Z_1^\alpha,Z_2^\alpha] = \frac{(\alpha-2) \Gamma(2\alpha+3)
+2(\alpha+2) \Gamma (\alpha+2)^2}{2(\alpha+1)^2 (\alpha+2)
\Gamma(2\alpha +3)} , 
\]
and this quantity is zero only if $\alpha=1$; it is positive for $\alpha>1$ 
and negative for $0<\alpha<1$.

\subsection{Limit behaviour for $\alpha >1$} \label{alphagr1}

%\subsection{The rooted case: the DLT}
We now consider the limit distribution of the
total weight 
of the
DLT and DLF.
 In the present
section we consider the case of $\alpha$-power weighted
edges with $\alpha >1$;
that is, we prove part (ii) of Theorem \ref{dltthm}.
To describe the moments of the limiting distribution of
 $D^\alpha(\UU_m^0)$ and $D^\alpha(\UU_m)$, we introduce the notation
 \bea
 J(\alpha) := \int_0^1 u^\alpha (1-u)^\alpha
 \ud u = 2^{-1-2\alpha} \sqrt{\pi} \frac{ \Gamma
(\alpha+1)}{\Gamma (\alpha+3/2)} . 
%\label{Jdef})
\label{Jdef}
\eea 
We start with the rooted case ($D^\alpha(\UU_m^0)$), and subsequently 
consider the unrooted case ($D^\alpha(\UU_m)$).

\begin{proposition} \label{1119f} Let $\alpha>1$. 
Then there exists a random variable $\Dalph$
such that as $m \to \infty$ we have
$D^\alpha(\UU_m^0) \to \Dalph$ almost surely and in $L^2$.
Also, the random variable $\Dalph$ 
satisfies the distributional fixed-point equality (\ref{0628bb}).
Further,
$\Exp[ \Dalph] = 1/(\alpha-1)$ and
\bea
\label{0701g}
 \Var [ \Dalph ] =
\frac{\alpha \left( \alpha -2 + 2(2\alpha+1) J(\alpha) \right)}
{(\alpha-1)^2(2\alpha-1)}.
\eea
\end{proposition}
\textbf{Proof.} 
Let $Z_i$ be the length of the $i$th edge of the DLT,
  as defined at (\ref{0728j}).
Let $ D_\alpha := \sum_{i=1}^\infty
Z_i^\alpha . $
The sum converges almost surely since
it has non-negative terms
and, by (\ref{asymZm}), has finite expectation for $\alpha >1$.
By (\ref{asymZm}) and Cauchy-Schwarz, there
exists a constant $0<C<\infty$ such that
\bean
\Exp[ D_\alpha^2 ] = \sum_{i=1}^\infty \sum_{j=1}^\infty
\Exp [ Z_i^\alpha Z_j^\alpha ] \leq C \sum_{i=1}^\infty \sum_{j=1}^\infty
i^{-\alpha} j^{-\alpha} < \infty, \eean
since $\alpha>1$. The $L^2$ convergence
then follows from the dominated convergence
theorem.

 Taking $U= X_1 $ here, by the self-similarity of the DLT we have that 
\bea
 D^\alpha (\UU_m^0) \eqd U^\alpha
D_{\{1\}}^\alpha (\UU_{N}^0) +(1-U)^\alpha D_{\{2\}}^\alpha
(\UU_{m-1-N}^0) +U^\alpha, 
\label{0714h}
\eea
 where
 $N \sim \textrm{Bin} (m-1, U)$,
given $U$, 
and, given $U$ and $N$, $D^\alpha_{\{1\}}(\UU_N^0)$ and
 $D^\alpha_{\{2\}}(\UU_{m-1-N}^0)$ 
are independent with the distribution of $D^\alpha(\UU_N^0)$ and
 $D^\alpha(\UU_{m-1-N}^0)$, respectively. 
As $m \to \infty$, $N$ and $m-N$ both tend to
infinity almost surely, and so, by taking
$m \to \infty$ in (\ref{0714h}), we obtain 
the fixed-point equation (\ref{0628bb}).

The identity $E[\Dalph] = (\alpha -1)^{-1}$
is obtained either from  (\ref{0720c}) of
Proposition \ref{dltmoms}, or by taking expectations in (\ref{0628bb}).
 Next, if we set $\tDalph = \Dalph- \Exp[ \Dalph]$, (\ref{0628bb})
 yields (\ref{0628b}). 
 Then,
 using the definition (\ref{Jdef}) of $J(\alpha)$,
the fact that $\Exp[ \tDalph]=0$, and independence,
we obtain from (\ref{0628b})
that
\bean
\Exp [ \tDalph^2 ] & = & \frac{2\Exp [ \tDalph^2 ]}{2\alpha +1}
+ \frac{\alpha^2+1}{(\alpha-1)^2 (2\alpha+1)}
 + \frac{2\alpha J(\alpha)}{(\alpha-1)^2}-\frac{1}{(\alpha-1)^2},
\eean
and rearranging this gives (\ref{0701g}). $\square$ \\

 Recall from Lemma  \ref{dickmanlem} 
that $\LL_0^\alpha$ is the limiting weight of edges 
attached to the origin in the DLT on uniform points.
Combining this fact with Proposition \ref{1119f},
we obtain a similar result to the latter for the
unrooted case as follows:

\begin{proposition} \label{1125a} 
 Let $\alpha > 1$. There is a random variable
$\Falph$,
 satisfying the
distributional fixed-point equality (\ref{0628dd}),
such that
$ D^\alpha ( \UU_m ) 
\to \Falph$, as $n \to \infty$, almost surely and in $L^2$. 
Further, $\Exp[ \Falph] = 1/(\alpha(\alpha-1))$, and
\bea
\label{0701i}
 \Var [ \Falph ] = \frac{1}{2\alpha} \Var [ \Dalph ]
+ \frac{ \alpha + 2(2\alpha+1) J(\alpha)-2}{2\alpha^2 (\alpha-1)^2} ,
\eea
where $J(\alpha)$ is given by (\ref{Jdef}) and
$\Var[\Dalph]$ by (\ref{0701g}). 
\end{proposition}
\proof
By Lemma \ref{dickmanlem} and Proposition \ref{1119f},
there are random variables $\Dalph$ and $\LL_0^\alpha$ such
that as $m \to \infty$ we have $D^\alpha ( \UU_m^0 ) \inLL \Dalph$
and $  \LL^\alpha_0 (\UU_m^0) \inLL \LL_0^\alpha$, also
with almost sure convergence in both cases.  Hence, setting
$\Falph := \Dalph - \LL^\alpha_0 $, we have by (\ref{0714b}) that
\begin{eqnarray}
D^\alpha (\UU_m) 
=  D^\alpha (\UU_m^0) - \LL^\alpha_0 (\UU_m^0) 
\to \Falph, ~~~~{\rm a.s.~~ and~in~}L^2.
\label{0715} 
\end{eqnarray} 

Next, we show that $\Falph$ satisfies
the distributional fixed-point equality (\ref{0628dd}).
The self-similarity of the DLT implies that
 \bea
 D^\alpha (\UU_m) \eqd U^\alpha D^\alpha
(\UU_{N}) +(1-U)^\alpha D^\alpha (\UU_{m-1-N}^0), 
\label{0728k}
\eea 
 where $N \sim
\textrm{Bin} (m-1, U)$, given $U$, and
$ D^\alpha (\UU_{N})$ and $D^\alpha (\UU_{m-1-N}^0)$
are independent, given $U$ and $N $.
As $m \to \infty$, $N$ and
$m-N$ both tend to infinity almost surely,  so taking
$m\to \infty$ in (\ref{0728k}), using Proposition
\ref{1119f} and 
eqn (\ref{0715}), we obtain the fixed-point equation (\ref{0628dd}). 

The identity $\Exp[\Falph]= \alpha^{-1}(\alpha-1)^{-1}$
is obtained either by (\ref{0720f}), or by taking expectations in
 (\ref{0628dd}) and using the formula for $E[\Dalph]$  in
  Proposition \ref{1119f}. Then with $\tFalph : = \Falph - \Exp [\Falph]$,
we obtain (\ref{0628d}) from (\ref{0628dd}), and
using independence and the fact that $\Exp[ \tFalph]=
\Exp[\tDalph]=0$ we obtain
\[
\frac{2\alpha}{2\alpha+1} \Exp[ \tFalph^2]
= \frac{\Exp[ \tDalph^2]}{2\alpha+1} + \frac{2\alpha J(\alpha)-1}{\alpha^2 (\alpha-1)^2}
+\frac{\alpha^2+1}{\alpha^2(\alpha-1)^2(2\alpha+1)},\]
which yields (\ref{0701i}).  $\square$ \\

{\bf Examples.}
When $\alpha=2$ we have that $\Exp[D_2]=1$ and
$J(2)=1/30$, so that $\Var[D_2] = 2/9$.
Also, $\Exp[F_2]=1/2$ and $\Var[F_2]=7/72 \approx 0.0972$.

\subsection{Limit behaviour for $\alpha =1$}
\label{subsecalphone}

Unlike in the case $\alpha > 1$, for $\alpha =1$ the mean of the total
weight $D^1(\UU_m^0)$ diverges as $m \to \infty$
 (see Proposition \ref{dltmoms}),
so clearly there is no limiting distribution
for $D^1(\UU_m^0)$. Nevertheless, by using the
orthogonality of the increments of the sequence
$(D^1(\UU_m^0),m \geq 1)$, we are able
to show
that the {\em centred} total weight $\tD^1(\UU_m^0)$ does converge
in distribution (in fact, in $L^2$) to a limiting
random variable,
and likewise for the unrooted case; this is our next result.

Subsequently, we shall 
 characterize the distribution of the limiting
random variable (for both the rooted and unrooted cases)
% $\tD^1$ (and  also
%the limit  in the unrooted case) 
by
a fixed-point identity, and thereby complete the proof of
Theorem \ref{dltthm} (i).
%show that both
%this random variable 
%have the distribution given by (\ref{0628a}).

\begin{proposition} \label{dconverge}
(i) As $m \to \infty$,
the random variable $\tD^1(\UU_m^0)$ converges in $L^2$ to a
limiting random variable $\tDone$, with
$\Exp[\tDone]=0$ and $\Var [ \tDone ] = 2 - \pi^2/6$.
In particular,
$
 \Var \left[ D^1(\UU_m^0)
\right] \to 2 - \pi^2/6$ as $m \to \infty$. 

(ii) As $m \to \infty$, $\tD^1(\UU_m) $ converges in $L^2$ 
to the limiting random variable $\tFone : = \tDone -\LL_0^1+1$.
\end{proposition}
\proof
Adopt the convention $D^1(\UU_0^0)=0$.
By the orthogonality of the
$Z_j$ (Lemma \ref{covz}) and (\ref{0630b}),
for $0 \leq \ell < m$,
\bean
\Var \left[ \tD^1(\UU_m^0) - \tD^1 (\UU_\ell^0) \right]
& = & \Var \sum_{j=\ell+1}^m \left( Z_j - \Exp [Z_j] \right) \\
& = & \sum_{j=\ell+1}^m \frac{j}{(j+1)^2(j+2)} \longrightarrow 0
\textrm{ as } m,\ell \to \infty.
\eean
Hence $\tD_1(\UU_m^0)$ is a Cauchy sequence in $L^2$,
and so converges in $L^2$ to a limiting random variable,
which we denote $\tDone$.
Then $\Exp[\tDone] = \lim_{m \to \infty} \Exp[ \tD_1(\UU_m^0)] =0$,
and
\bean
 \Var [ \tDone] = 
\lim_{m \to \infty} \Var \left[ \tD^1(\UU_m^0) \right]
=
 \sum_{j=1}^\infty \frac{j}{(j+1)^2(j+2)}
\nonumber \\
= \sum_{j=1}^\infty \left[ \frac{2}{j+1} - \frac{2}{j+2} \right]
- \sum_{j=1}^\infty \frac{1}{(j+1)^2} = 1 - \left( \frac{\pi^2}{6}
-1 \right)
= 2-\frac{\pi^2}{6} .  
\eean
It remains to prove part (ii), the convergence for
the centred total length of the DLF
$\tD^1 (\UU_m)$. We have by (\ref{0714b}) that
\begin{eqnarray}
\tD^1 (\UU_m) & =
& \tD^1 (\UU_m^0) - \LL_0^1 (\UU_m^0)+\Exp[ \LL_0^1 (\UU_m^0)]
 %\nonumber\\ & \inLL &
\inLL
\tDone - \LL^1_0 + 1 , \nonumber \end{eqnarray}
where the convergence follows by 
Lemma \ref{dickmanlem} and part (i). 
 Thus $\tD^1 ( \UU_m)$
converges in $L^2$ as $m \to \infty$. $\square$ \\

For the next few results it is more convenient to consider the DLF
defined on a Poisson  number of points. 
Let $(X_1,X_2,\ldots)$ be a sequence of independent 
uniformly distributed random variables in $(0,1]$, and let
$(N(t),t\geq 0)$ be the counting process of a homogeneous Poisson process
of unit rate in $(0,\infty)$, independent of $(X_1,X_2,\ldots)$.
Thus $N(t)$ is a Poisson variable with parameter $t$.
As before, let $\UU_m =(X_1,\ldots,X_m)$, and 
(for this section only) let
 $\Po_t := \UU_{N(t)}$. Let
 $\Po_t^0 := \UU_{N(t)}^0$,
so that
 $\Po_t^0 = (0,X_1,X_2,\ldots,X_{N(t)})$. 

We construct the DLF and DLT on $X_1, X_2, \ldots,X_{N(t)}$ as before.
Let $\tD^1(\Po_t^0) = D^1(\Po_t^0) - \Exp \left[ D^1(\Po_t^0)
\right]$ and $\tD^1(\Po_t) = D^1(\Po_t) - \Exp \left[ D^1(\Po_t)
\right]$.
We aim to show that the limit distribution for $\tD^1(\Po_t^0)$
is the same as for $\tD^1(\UU_m^0)$, and likewise in the unrooted
case.  We shall need the following result.

\begin{lemma} \label{1120a} As $t \to \infty$, \begin{equation}
\frac{\ud}{\ud t} \Exp[D^1 (\Po_t)] = \frac{1}{t} + O(t^{-2}); 
%\textrm{ and }
{\rm ~~ and ~~ }
\frac{\ud}{\ud t} \Exp[D^1 (\Po_t^0)] = \frac{1}{t} + O(t^{-2}).
\end{equation} \end{lemma}
\textbf{Proof.} 
The point set $\{X_1,\ldots,X_{N(t)}\}$ is a homogeneous
Poisson point process in $(0,1)$, so we have
\bean
 \frac{\ud}{\ud t} \Exp[D^1 (\Po_t)] & = &
\Exp[\textrm{length of new arrival}] \\ & = & \int_0^1
\ud u \Exp[\textrm{dist.~to next pt.~to the left of $u$ in
$\mathcal{P}_t$}] \\
& = &
\int_0^1 \ud u
\int_0^u st e^{-ts} \ud s
= \frac{1}{t} + \frac{2}{t^2} \left( e^{-t} - 1 \right) +\frac{e^{-t}}{t}
\\ & = & \frac{1}{t} + O \left( t^{-2} \right).
\eean
Similarly, 
\bean
\frac{\ud}{\ud t} \Exp[D^1 (\Po_t^0)] & = &
 \int_0^1
\ud u \Exp[\textrm{dist.~to next pt.~to the left of $u$ in
$\Po_t \cup \{ 0 \}$}]
\\ & = &
\int_0^1 \ud u
\int_0^u \Pr[\textrm{dist.~to next pt.~to the left} > s] \ud s
\\ & = & \int_0^1 \ud u \int_0^u e^{-ts} \ud s
= \frac{1}{t}+\frac{e^{-t}-1}{t^2}
\\ & = &
\frac{1}{t} + O \left( t^{-2} \right). \; \square
\eean

\begin{lemma}
\label{0630c}
(i) As $t \to \infty$, $\tD^1(\Po_t^0)$
 converges in distribution
 to $\tDone$, the $L^2$ large-$m$ limit of
$\tD^1(\UU_m^0)$.

(ii) As $t\to \infty$, $\tD^1(\Po_t)$ converges in distribution to
$\tFone$, the $L^2$ large-$m$ limit of $\tD^1(\UU_m)$.
\end{lemma} \proof (i) From Proposition~\ref{dconverge}, we have
$
 \tD^1 (\UU_m^0) \inLL \tDone
$
 as $m \to \infty$.
Let $a_t := \Exp[D^1(\Po_t^0)]$ and $\mu_m :=\Exp[D^1(\UU_m^0)]$. 
Since $\mu_m = \Exp \sum_{i=1}^m Z_i = \sum_{i=1}^m (1+i)^{-1}$ by
(\ref{0630b}), for any positive integers $\ell, m$  we have
\bea 
\left| \mu_{m} - \mu_{\ell } \right| =
\sum_{j= \min ( m, \ell ) +1}^{\max ( m, \ell ) }
\frac{1}{j+1} \leq \log\left( 
\frac{ \max(m,\ell)+1 }{ \min(m,\ell)+1 }
\right) = \left| \log\left( 
\frac{ m+1 }{ \ell+1 }
\right)
\right|.
\label{0720g}
 \eea
%The number of arrivals up to time $t$, $N_t$, is Poisson with
%parameter $t$.
Note the distributional equalities
 \bea
\mathcal{L} \left( D^1(\Po_t^0) | N(t) = m \right) = \mathcal{L}
\left( D^1(\UU_m^0) \right) ; \nonumber \\
\mathcal{L} \left( D^1(\Po_t^0) - \mu_{N(t)} | N(t) = m \right) =
\mathcal{L} \left( \tD^1( \UU_m^0) \right)  \label{0720h}.
\eea

 First we aim to show that
$ a_t - \mu_{\lfloor t \rfloor} \to 0 $
as
$ t \to \infty$.
 Set $p_m(t) :=
e^{-t} \frac{t^m}{m!}$.
Then we
can write
 \bea \label{eq200} a_t -\mu_{\lfloor t \rfloor } & = & 
\sum_{m=0}^\infty
p_m(t) (\mu_m - \mu_{\lfloor t \rfloor } ) \nonumber \\ & = & 
\!\! \sum_{\left| m-\lfloor t \rfloor \right| \leq
t^{3/4}} \!\! p_m(t) (\mu_m -\mu_{\lfloor t \rfloor} )
 + \!\! \sum_{\left| m-\lfloor t \rfloor
\right| > t^{3/4}}  \!\! p_m(t) (\mu_m - \mu_{\lfloor t\rfloor}) .
 \eea
 We examine these two sums
separately. First consider the sum for $|m-\lfloor t \rfloor |
\leq t^{3/4}$. By (\ref{0720g}), we have
\bean
\sup_{m: | m - \lfloor t \rfloor | \leq t^{3/4}}
\left| \mu_m - \mu_{\lfloor t \rfloor} \right|
\leq \max \left(  
\log{ \left( \frac{\lfloor t
 \rfloor+1+t^{3/4}}{\lfloor t \rfloor+1} \right) } ,
\log \left( \frac{\lfloor t \rfloor+1}{\lfloor t \rfloor+1-t^{3/4}} \right)
\right) \nonumber \\
=O \left(t^{-1/4} \right)
\to 0 \textrm{ as } t \to \infty.
\eean
Hence the first sum in (\ref{eq200}) tends to zero as $t \to \infty$.
To estimate the second sum, observe that
\begin{eqnarray} \label{1119c} \sum_{\left| m -
\lfloor t \rfloor \right|
> t^{3/4}}  \!\!\!\! p_m(t) (\mu_m  - \mu_{\lfloor t \rfloor })
& \leq &  \sum_{\left|
m - \lfloor t \rfloor \right| > t^{3/4}} \!\!\!\! p_m(t) (m +t) 
\nonumber \\
& = & \Exp \left[ (N(t) +t) \1 {\{ |N(t)-\lfloor t \rfloor| >
t^{3/4} \} } \right] \nonumber\\ & \leq & \left( \Exp \left[
(N(t)+t)^2 \right] 
%\right)^{1/2} \left(
\cdot \Pr \left[ \left| N(t) -
\lfloor t \rfloor \right| > t^{3/4} \right] \right)^{1/2} \!\!\!. ~~~~
\end{eqnarray} 
By Chernoff bounds on the tail probabilities of a Poisson
random variable (e.g.~Lemma 1.4 of 
 \cite{penbook}), 
the expression (\ref{1119c}) is $O(t \exp(-t^2/18))$ and so tends to zero. 
Hence the second sum in (\ref{eq200}) tends to zero, and thus
 \bea
\label{0630d}
 a_t - \mu_{\lfloor t
\rfloor} 
 \to 0 \textrm{ as } t \to \infty .
\eea

 Now we show that $\tD^1(\Po_t^0) \tod \tDone$ as $t \to\infty$. We have
\begin{equation}
\label{0630g}
\tD^1(\Po_t^0) = \left( D^1(\Po_t^0) - \mu_{N(t)} \right) +
\left( \mu_{N(t)} - \mu_{\lfloor t \rfloor} \right) + \left(
\mu_{\lfloor t \rfloor} -a_t \right) .
\end{equation} 
The final bracket tends to zero,
%is $O(t^{-1/4})$, from
by (\ref{0630d}).
Also, by (\ref{0720h}) and the fact that
 $N(t) \to \infty \textrm{ a.s.}$ as $t \to \infty$, we have
 \[
D^1(\Po_t^0) - \mu_{N(t)} \tod \tDone . 
\]
Finally, using 
(\ref{0720g}),
 we have
\bean \left| \mu_{N(t)} - \mu_{\lfloor t \rfloor} \right|
& \leq & \left| \log \frac{ N(t)+1}{\lfloor t \rfloor+1}
\right| 
 \toP 0 ,
\eean
as $t \to \infty$, since $N(t)/ \lfloor t \rfloor \toP 1$.
So Slutsky's theorem
applied to (\ref{0630g}) yields
$\tD^1(\Po_t^0) \tod \tDone$ as $t \to\infty$, completing
the proof of (i)

The proof of (ii) follows in the same way as that of (i),
except that in (\ref{0720g}) the first equals sign 
 is replaced by an inequality $\leq$. This does not
affect the rest of the proof. 
%with some slight modifications.
$\square$ \\

The next two propositions complete the proof of Theorem \ref{dltthm}.

\begin{proposition} \label{1020d}
The limiting random variable $\tDone$ of Proposition
\ref{dconverge} (i) satisfies the fixed-point
equation (\ref{0628a}).
\end{proposition} \proof
For integer $n >0$,
let $T_n :=\min \{ s: N(s) \geq n\}$, the $n$th  arrival
time of the  Poisson process with counting process $N(\cdot)$.
Set $T:= T_1$, and set $U:= X_1$
(which is uniform on $(0,1)$).

By the Marking Theorem for Poisson processes \cite{kingman},
the two-dimensional point process ${\cal Q} := \{(X_n,T_n): n \geq 1\}$
is a homogeneous Poisson process of unit intensity
on $(0,1)\times (0,\infty)$. Given the value of $(U,T)$,
the restriction of ${\cal Q}$ to $(0,U]\times (T,\infty)$ and
the restriction of ${\cal Q}$ to $(U,1]\times (T,\infty)$ are
independent homogeneous Poisson processes on these regions.
Hence, by scaling properties of the Poisson process
(see the Mapping Theorem in \cite{kingman}) and of the DLT,
writing $D^1_{\{i\}}(\cdot)$, $i=1,2$ for independent
copies of $D^1 (\cdot)$, we have
\bea \label{0128a0}
D^1 (\Po_t^0)  & \eqd &  \left( U D^1_{\{1\}}
(\Po_{U(t-T)}^0)+ (1-U) D^1_{\{2\}}(\Po_{(1-U)(t-T)}^0) +U \right)
\1 \{t>T\}. ~~~ 
\eea
Let $a_s=0$ for $s\leq 0$, and $a_s=\Exp [ D^1 (\Po_s^0) ]$ for $s > 0$.
Then $\tD^1(\Po_t^0) = D^1(\Po_t^0)-a_t$, so that
 by (\ref{0128a0}),
\bea
\label{0128a}
 \tD^1 (\Po_t^0)  & \eqd &  
 \left( U \tD^1_{\{1\}} (\Po_{U(t-T)}^0)+ (1-U) \tD^1_{\{2\}}
(\Po_{(1-U)(t-T)}^0) +U \right) \1 \{ t>T \}
\nonumber\\ &  & +U \left( a_{U(t-T)}-a_t
\right) +(1-U) \left( a_{(1-U)(t-T)} -a_t \right).
\eea From Lemma \ref{1120a} we have
$ \frac{\ud a_t}{\ud t} = \frac{1}{t}+O(t^{-2}). $
Hence, if $T < t$, then
\bean
 & & a_t- a_{U(t-T)}  = \int_{U(t-T)}^t
\frac{\ud a_s}{\ud s} \ud s
 = \log{t} - \log \{U(t-T)\} +O\left( (U(t-T))^{-1} \right), 
\eean
and hence as $t \to \infty$,
\bea
 a_t-a_{U(t-T)} \to -\log U, ~~~{\rm a.s.} .
\label{0714e}
\eea
 Since $P[T<t ] $ tends to 1, 
by making $t \to \infty$ in 
 (\ref{0128a}) and using Slutsky's theorem we obtain (\ref{0628a}).
$\square$

\begin{proposition} \label{ffixed}
The limiting random variable $\tFone$ of Proposition
\ref{dconverge} (ii) satisfies the fixed-point
equation (\ref{0628a}), and so has the same distribution
as $\tDone$. Also, $\Cov(\tFone,\tDone) = (7/4)- \pi^2/6$.
\end{proposition} \textbf{Proof.} The proof follows similar lines to
that of Proposition \ref{1020d}. Once more let $a_s = \Exp [ D^1 (\Po_s^0) ]$, for
$s \geq 0$, and $a_s=0$ for $s <0$.
Let $b_s = \Exp [ D^1 (\Po_s) ]$ for $s > 0$, and $b_s=0$
for $s \leq 0$, and let $T := \min\{t: N(t) \geq 1\}$,
Then
 \begin{eqnarray} D^1(\Po_t) & \eqd & \left( U
D^1_{\{1\}} (\Po_{U(t-T)})+ (1-U) D^1_{\{2\}} (\Po_{(1-U)(t-T)}^0) \right)
\1 \{t>T\} ,
\label{0720j}
\end{eqnarray}
where $D^1_{\{1\}}(\cdot)$ and $D^1_{\{2\}}(\cdot)$ are
independent copies of $D^1(\cdot)$. 
Then $\tD^1(\Po_t) = D^1(\Po_t) - b_t$ and  
$\tD^1(\Po_t^0) = D^1(\Po_t^0) - a_t$, so that
(\ref{0720j}) yields
\begin{eqnarray}
 \tD^1 (\Po_t)  & \eqd &
 \left( U \tD^1_{\{1\}}(\Po_{U(t-T)})+ (1-U) \tD^1_{\{2\}} 
(\Po_{(1-U)(t-T)}^0)
\right) \1 \{t>T\}
\nonumber\\ &  & +U \left( b_{U(t-T)}-b_t \right) +(1-U) \left(
a_{(1-U)(t-T)} -b_t \right). 
\label{0714f}
\end{eqnarray} From
Lemma \ref{1120a} we have
$\frac{\ud b_t}{\ud t} = \frac{1}{t}+O(t^{-2})$.
Hence, by the same argument as used at (\ref{0714e}),
 \[ b_t- b_{U(t-T)}   \to
 -\log{U}  ~~~ {\rm a.s.} \]
Also,
$a_t -b_t = \Exp [\LL_0^1(\Po_t^0)]$ by (\ref{0714b}), so that
 $\lim_{t \to \infty} (a_t-b_t) = 1$, by Lemma
\ref{dickmanlem} and the fact that $\Exp[ \LL_0^1]=1$ (eqn (\ref{0630j})).
Using also (\ref{0714e}) we find that
 as $t \to \infty$,
\[ a_{(1-U)(t-T)} - b_t 
= (a_{(1-U)(t-T)} -a_t) + (a_t-b_t) 
\to 1 + \log{(1-U)} , ~~~~{\rm a.s.}  \] Taking
$t \to \infty$ in (\ref{0714f}),
and using Slutsky's theorem, we obtain
\bea
\label{0713a}
 \tFone  \eqd  U \tFone+ (1-U) \tDone +U \log {U}
+(1-U) \log{(1-U)} +(1-U). \eea
The change of variable $(1-U) \mapsto U$ then shows that $\tDone$
as defined at (\ref{0628a}) satisfies (\ref{0713a}), and so by the uniqueness
of solution, $\tFone$ has the same distribution as $\tDone$ and satisfies
(\ref{0628a}).

To obtain the covariance of $\tFone$ and $\tDone$, observe from
Proposition \ref{dconverge} (ii) that $\LL^1_0 = \tDone-\tFone + 1$,
and therefore 
 by (\ref{0630j}), we have that
\bea
1/2 = \Var [\LL^1_0] = 
\Var[\tDone]   + \Var[\tFone] - 2 \Cov(\tDone,\tFone).   
\label{0817}
\eea
Since 
$\Var[\tFone]=    \Var[\tDone] = 2 - \pi^2/6$ by
Proposition \ref{dconverge} (i), rearranging (\ref{0817})
we find that $\Cov(\tDone,\tFone) = (7/4)-\pi^2/6$.
 $\square$ \\

\begin{figure}[h!]
\begin{center}
\includegraphics[angle=270, width=0.8\textwidth]{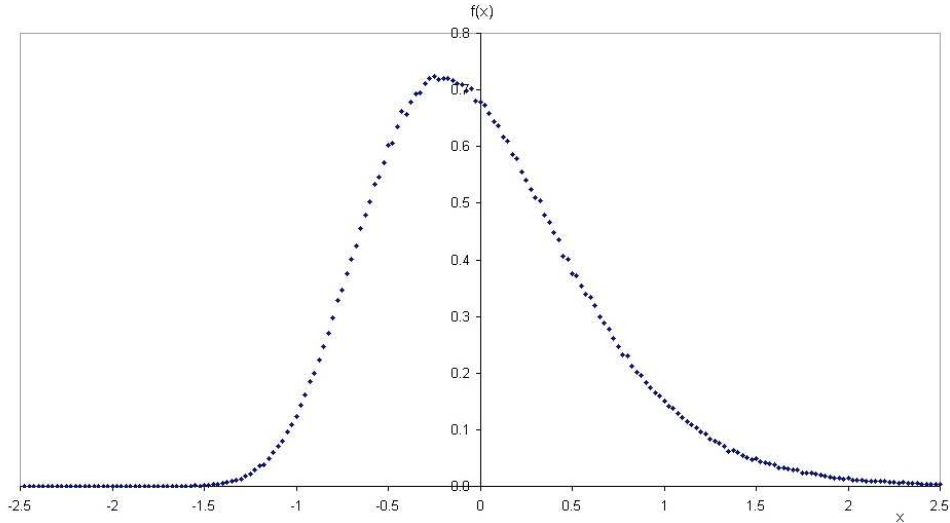}
\end{center}
\caption{Estimated probability density function for $\tDone$.}
\label{pdffig}
\end{figure}

\rem
Figure \ref{pdffig} is a plot of the estimated probability
density function of $\tDone$. This was obtained by
performing $10^6$ repeated
simulations of the DLT on a sequence of
$10^3$ uniform (simulated) random points on $(0,1]$. For each
simulation, the expected value
of $D^1(\UU_{10^3})$ (which is precisely $(1/2) + (1/3) + \cdots
(1/1001)$ by Lemma \ref{Zmdist}) was subtracted
from the total length of the simulated DLT
to give an approximate realization of $\tDone$.
The density function was then estimated from
the sample of $10^6$ approximate realizations
of $\tDone$, using a window width of $0.0025$.
The simulated
sample from which the density
estimate
for $\tDone$ was taken had sample mean $\approx -2 \times 10^{-4}$
and sample variance $\approx 0.3543$, which are reasonably
close to the expectation and variance of $\tDone$.

\section{General results in geometric probability}
\label{secgeneral}

Notions of {\em stabilizing} functionals
of point sets have recently proved to be a useful basis for
a general methodology for establishing  limit
theorems
for functionals 
of random
point sets in $\R^d$.
In particular,  Penrose and Yukich \cite{penyuk1,penyuk2}
provide general central limit theorems and laws of large numbers
for stabilizing   functionals.
One might hope to apply  these  results
 in the case of the MDSF weight. 
In fact we shall obtain our law of large numbers (Theorem
\ref{llnthm}) by application of a result from 
\cite{penyuk2}, but to obtain the central limit
theorem  for edges away from the boundary in the MDSF and MDST,
we need an extension of the general result in \cite{penyuk1}. It
is these general results that we describe in the present section.

For our general results, we use the following notation.
Let $d \geq 1$ be an integer. For $\X \subset \R^d$, 
constant $a>0$, and ${\bf y} \in
\R^d$, let ${\bf y}+a\X$ denote the transformed set $\{ {\bf y} + a\bx
: \bx \in \X\}$. Let $\diam(\X) :=\sup \{ \|\bx_1-\bx_2\|:
\bx_1,\bx_2 \in \X \}$, and let $\card(\X)$ denote the cardinality
(number of elements) of $\X$ (when finite).

For $\bx \in \R^d$ and $r>0$, let $B(\bx;r)$ denote the closed Euclidean ball
with centre $\bx$ and radius $r$, and let $Q(\bx;r)$ denote the
corresponding $l_\infty$ ball, i.e., the $d$-cube $\bx + [-r,r]^d $. For 
bounded measurable $R \subset \R^d$ let $|R|$ denote the Lebesgue measure 
of $R$, let $\partial R$ denote the topological boundary of $R$ and
 for $r>0$, set $\partial _r R := \cup_{\bx \in \partial R} 
Q(\bx;r)$, the $r$-neighbourhood of the boundary of $R$.

\subsection{A general law of large numbers}
\label{secgenlln}

Let $\xi ( \bx ; \X )$ be a measurable $\R_+$-valued function
defined for all pairs $(\bx , \X)$, where $\X \subset \R^d$ is finite
and $\bx \in \X$. Assume $\xi$ is translation invariant, that is,
 for all $\by \in \R^d$,
 $\xi( \by + \bx; \by+\X) = \xi (\bx ; \X)$.
When $\bx \notin \X$, we abbreviate the notation
$\xi(\bx;\X \cup \{\bx\})$ to $\xi(\bx;\X)$.

For our general law of large numbers, we use a notion
of stabilization defined as follows.
For any locally finite point set $\X \subset \R^d$ and any $\ell
\in \N$ define
\[ \xi^+(\X;\ell) := \sup_{k \in \N} \left( \mathrm{ess}
\sup_{\ell,k} \left\{ \xi(\0;(\X \cap B(\0;\ell)) \cup \mathcal{A}
\right\} \right) \textrm{, and}
\]
\[ \xi^-(\X;\ell) := \inf_{k \in \N} \left( \mathrm{ess}
\inf_{\ell,k} \left\{ \xi(\mathbf{0};(\X \cap B(\0;\ell)) \cup
\mathcal{A} \right\} \right) ; \] where $\mathrm{ess}
\sup_{\ell,k}$ is the essential supremum, with respect to Lebesgue
measure on $\R^{dk}$, over sets $\mathcal{A} \subset \R^d
\backslash B(\0;\ell)$ of cardinality $k$. Define the {\em limit} of
$\xi$ on $\X$ by
\[ \xi_{\infty}(\X) := \limsup_{k \to \infty}
\xi^+(\X;k)  . \] We say the functional $\xi$ \emph{stabilizes} on
$\X$ if 
\bea
 \lim_{k \to \infty} \xi^+(\X;k) = \lim_{k \to \infty}
\xi^-(\X;k) = \xi_{\infty} (\X)  .
\label{stabeq}
\eea

 For $\tau \in
(0,\infty)$, let $\H_{\tau}$ be a homogeneous Poisson
process of intensity $\tau$ on $\R^d$. 
The following general law of large numbers 
is due to Penrose and Yukich \cite{penyuk2}.
We shall use it to prove Theorem \ref{llnthm}.
\begin{lemma} \cite{penyuk2}
\label{llnpenyuk} Suppose $q=1$ or $q=2$. Suppose $\xi$ is
almost surely stabilizing on $\H_{\tau}$, with limit
$\xi_{\infty}(\H_\tau)$, for all $\tau \in (0,\infty )$.
Let $f$ be a probability density function on $\R^d$, and
let $\X_n$ be the point process consisting of $n$ independent
random $d$-vectors with common density $f$.
If $\xi$ satisfies the moments condition \begin{equation}
\label{0k715a}
 \sup_{n
\in \N} \Exp \left[ \xi \left( n^{1/d}\mathbf{X}_1;n^{1/d}
\X_n \right) ^p \right] < \infty, \end{equation}
 for some $p>q$, then
as $n \to \infty$,
 \begin{equation} \label{0k715b}
n^{-1}
\sum_{\bx \in \X_n} \xi( n^{1/d}\bx; n^{1/d} \X_n )
 \stackrel{L^q}{\longrightarrow} \int_{\R^d} \Exp \left[
\xi_{\infty} \left( \H_{f(\mathbf{x})} \right) \right]
f(\mathbf{x}) \ud \mathbf{x} ,
\end{equation} and the limit is finite. \end{lemma}
%The proof of the result is given in~\cite{penyuk2}. 

\subsection{General central limit theorems}
\label{subsecgenclt}

In the course of the proof of Theorem \ref{mainth}, we
shall use
 a modified form of a 
 general central limit theorem obtained for
functionals of geometric graphs  by
Penrose and Yukich  \cite{penyuk1}. We recall
the setup of \cite{penyuk1}. 
As in Section \ref{secgenlln},
let $\xi (\bx;\X) $ be a translation invariant real-valued functional
defined for finite $\X \subset \R^d$ and $\bx \in \X$.
 Then $\xi$ induces a translation invariant functional $H(\X;S)$
 defined on all 
finite point sets $\X \subset \R^d$ and all Borel-measurable
regions $S \subseteq \R^d$  by 
\bea
 H (\X ; S) : = \sum_{\bx \in \X \cap S} \xi (\bx ; \X) .
\label{0714}
\eea
It is this `restricted' functional that
 interests us here, while \cite{penyuk1}
is concerned rather with the global functional $H( \X; \R^d)$.
In our particular application (the length of edges  of the MDST
on random points in a square),
the global functional fails to satisfy the conditions
of the central limit theorems in \cite{penyuk1}, owing to
boundary effects.  Here we generalize the result
in \cite{penyuk1} to the
`restricted' functional $H(\X;S)$. 
It is this generalized result that  we can apply to the MDST,
  when we take  $S$ to be a region `away from
the boundary' of the square in which the random points are placed.

We use a notion of stabilization for $H$ which
is related to, but not equivalent to, the notion
of stabilization of $\xi$ used in Section
\ref{secgenlln}. Loosely speaking, $\xi$ is stabilizing
if when a point inserted at the origin into a homogeneous
Poisson process,  only nearby
 Poisson points affect the inserted point;
 for $H$ to be stabilizing we 
require also that the the inserted
point affects only nearby points.

For $B \subseteq \R^d$, let $\Delta(\X;B)$ denote the `add one
cost' of the functional $H$ on the insertion of a point at the
origin,
\[ \Delta( \X;B) := H ( \X \cup \{ \0 \};B) - H (\X;B ) . \]
Let $\Po := \H_1$ % i.e., let $\Po$ be 
(a homogeneous Poisson point process of unit intensity on $\R^d$).
Let $\QQ_n:= \Po \cap R_n$ 
(the restriction of $\Po$ to  $R_n$).
Adapting the ideas of~\cite{penyuk1}, we make the following definitions.
\begin{definition} 
\label{sstabdef}
We say the functional $H$
 is {\em strongly stabilizing} if
 there exist almost
surely finite random variables $R$ (a {\em radius of
stabilization}) and $\Delta (\infty )$ such that, with
probability 1, for any $B \supseteq B( \0 ; R )$,
\[ \Delta ( \Po \cap  B(\0;R) \cup {\cal A}; B ) = \Delta
(\infty ), ~ \forall \textrm{ finite } {\cal A} \subset \R^d \setminus
B(\0;R) .\]
\end{definition}

We say that the functional $H$ is \em polynomially bounded \em 
if, for all $B \ni \0$,
there exists a constant $\beta$ such that for all finite sets $\X
\subset \R^d$,
\bea \label{poly} \left| H(\X; B) \right| \leq \beta \left(\diam
(\X) + \card (\X) \right)^{\beta} . \eea

We say that $H$ is \textit{homogeneous of order}
$\gamma$ if for all finite $\X \subset \R^d$ and Borel $B \subseteq \R^d$,
%on which $H$ is defined,
and all $a\in \R$, $H (a\X ; aB) = a^{\gamma} H (\X ; B)$.

Let $(R_n, S_n)$, for $n=1,2,\ldots$, be a sequence of ordered pairs of
bounded Borel
subsets of $\R^d$, such that $S_n \subseteq R_n$ for all $n$.
Assume that
 for all $r>0$,
 $n^{-1} | \partial_r R_n | \to 0$ and $n^{-1} | \partial_r S_n | \to 0$
(the {\em vanishing relative boundary condition}).
Assume also that
% with the following conditions on $R_n, S_n$: that
 $|R_n| = n$ for all $n$,
and $|S_n|/n \to 1$ as $n \to \infty$;
that 
 $S_n$ tends to $\R^d$, in the sense that
 $\cup_{n \geq 1} \cap_{m \geq n} S_m = \R^d$;
 and
that there
exists a constant $\beta$ such that $\diam (R_n) \leq \beta n^\beta$
 for all $n$
(the {\em polynomial boundedness} condition on $(R_n,S_n)_{n \geq 1}$).
  Subject to these conditions, the choice of
 $(R_n, S_n)_{n \geq 1}$ is arbitrary.

Let $\bU_{1,n},\bU_{2,n},\ldots$ be i.i.d.~uniform random vectors on
$R_n$. Let
 \[ \UU_{m,n}=\{ \bU_{1,n},\ldots,\bU_{m,n} \}\] 
(a binomial point process), and for Borel $A \subseteq \R^d$ with
$0 < |A|< \infty$,
let $\UU_{m,A}$ be the binomial point process of $m$
i.i.d.~uniform random vectors on $A$.
%Let $\Po_n$ be a homogeneous Poisson
%process of intensity $1$ on $R_n$.

Let $\RR$ be the collection of all pairs $(A,B)$ with 
$A, B \subset \R^d$ of
the form $(A,B) = (\bx+R_n, \bx+ S_n ) $ with $ \bx \in \R^d$ and
$ n\in \N$.
That is, $\RR$ is the
collection of all the $(R_n,S_n)$ and their translates. 

We say that the functional $H$ satisfies the 
\em uniform bounded moments condition \em on $\RR$ if
\bea \label{ubm} \sup_{(A,B) \in \RR : \0 \in A}
\left(
\sup_{|A|/2 \leq m \leq 3|A|/2} \{ \Exp [ \Delta (
\UU_{m,A} ; B )^4 ] \}
 %\right\}
 \right)
 < \infty . \eea

We now state the general results, which extend
those of Penrose and Yukich (Theorem 2.1 
and Corollary 2.1 in~\cite{penyuk1}).

\begin{theorem}
\label{genclt}
Suppose that $H$ is strongly stabilizing,
 is polynomially bounded (\ref{poly}), 
and satisfies the uniform
bounded moments condition (\ref{ubm}) on $\RR$.
Then there exist constants
$s^2$, $t^2$, with $0 \leq t^2 \leq s^2$, such that as $n\to
\infty$,
\begin{itemize} \item[(i)] $n^{-1}
\Var\left( H \left( \QQ_n ; S_n \right) \right) \to s^2$;
\item[(ii)] $n^{-1/2}\left(H\left(\QQ_n ; S_n  \right)-
\Exp\left[ H \left(\QQ_n ; S_n  \right) \right] \right) \tod
 \NN
\left(0,s^2\right)$; \item[(iii)] $n^{-1} \Var \left( H \left(
\UU_{n,n} ; S_n  \right)\right) \to t^2$; \item[(iv)]
$n^{-1/2}\left(H\left(\UU_{n,n} ; S_n  \right) - \Exp \left[ H
\left( \UU_{n,n} ; S_n  \right)\right]\right) \tod
 \NN (0,t^2)$.
\end{itemize} Also, $s^2$ and $t^2$ are
independent of the choice of the $(R_n,S_n)$. Further, if the distribution
of $\Delta(\infty)$ is nondegenerate, then $s^2\geq t^2>0$.
\end{theorem}
%For the proof of the theorem, see~\cite{penyuk1}.

Let $R_0$ be a fixed bounded Borel subset of $\R^d$ with $|R_0|=1$ and
$| \partial R_0 | =0$. Let $(S_{0,n}, n \geq 1)$ be a sequence of
Borel sets with
$S_{0,n} \subseteq R_0$ 
such that $|S_{0,n}| \to 1$ as $n\to \infty$ and
for all $r >0$ we have
 $| \partial_{n^{-1/d} r} S_{0,n} | \to 0$ as $n \to \infty$

Let $\RR_0$ be the collection of all pairs of the form
 $(\bx + n^{1/d} R_0, \bx +n^{1/d} S_{0,n})$
 with $n \geq 1$ and $\bx \in \R^d$.
Let $\X_n$ be the binomial point process of $n$ i.i.d.~uniform
random vectors on $R_0$, and let $\Po_n$ be a homogeneous
Poisson point process of intensity $n$ on $R_0$.

\begin{corollary}
\label{cltcor}
Suppose $H$ is strongly stabilizing, satisfies the uniform
bounded moments condition on $\RR_0$, is polynomially
bounded and is homogeneous of order $\gamma$. 
Then with $s^2, t^2$ as in Theorem \ref{genclt}
we have that,
as $n \to \infty$
\begin{itemize}
  \item[(i)] $n^{(2 \gamma/d)-1} \Var \left(
H \left( \Po_n ; S_{0,n} \right) \right) \to s^2$; 
\item[(ii)]
$n^{(\gamma/d)-1/2} \left( H \left( \Po_n ; S_{0,n}  \right) - \Exp
\left[ H \left( \Po_n ; S_{0,n}  \right) \right] \right) \tod \NN
\left( 0, s^2 \right)$;
 \item[(iii)]
 $n^{(2 \gamma/d)-1} \Var \left(
H \left( \X_n ; S_{0,n}  \right) \right) \to t^2$; \item[(iv)]
$n^{(\gamma/d)-1/2} \left( H \left( \X_n ; S_{0,n}  \right) - \Exp
\left[ H \left( \X_n ; S_{0,n}  \right) \right] \right) \tod \NN
\left( 0, t^2 \right)$.
\end{itemize} \end{corollary}
\proof The corollary follows from Theorem \ref{genclt}
by taking $R_n =
n^{1/d} R_0$ and $S_n = n^{1/d} S_{0,n}$ (or suitable translates thereof),
and scaling, since $H$ is homogeneous of order $\gamma$. $\square$

\subsection{Proof of Theorem \ref{genclt}: the Poisson case}

Let $\Po$ be a Poisson process of unit intensity on $\R^d$.
We say the functional $H$ is {\em weakly stabilizing} on $\RR$ if
there is a random variable $\dinf$ such that
 \bea 
\label{weakstab}
\Delta ( \Po \cap A ; B) \toas \dinf, 
\eea 
as $(A,B) \to \R^d$ through
$\RR$, by which we mean (\ref{weakstab}) holds whenever 
$(A,B)$ is an $\RR$-valued sequence of the form
$(A_n,B_n)_{n \geq 1}$, such that
%$\cup_{n \geq 1} \cap_{m \geq n} A_n = \R^d$ and
$\cup_{n \geq 1} \cap_{m \geq n} B_m = \R^d$. Note that
 strong stabilization of $H$ implies
weak stabilization of $H$.

We say the functional $H$ satisfies the {\em Poisson bounded
moments condition} on $\RR$ if
\bea \label{Pbmc} \sup_{(A,B) \in \RR:
\0 \in A} \left\{ \Exp \left[ \Delta ( \Po \cap A ; B)^4 \right]
\right\} < \infty . \eea

\begin{theorem}
\label{0217c}
 Suppose that $H$ is weakly stabilizing on $\RR$
(\ref{weakstab}) and satisfies (\ref{Pbmc}). Then there exists
$s^2 \geq 0$ such that as $n \to \infty$, $n^{-1} \Var [ H
(\QQ_n; S_n)] \to s^2$ and 
$n^{-1/2} ( H( \QQ_n; S_n) - \Exp [ H(\QQ_n; S_n)])
\tod \NN (0,s^2)$.
\end{theorem}
Before proving Theorem \ref{0217c}, we require
further definitions and a lemma.
Let $\Po'$ be an independent copy of the Poisson process
 $\Po$. For $\bx \in \Z^d$, set
\[ \Po'' (\bx) = \left( \Po \setminus Q (\bx; 1/2) \right)
\cup \left( \Po' \cap Q( \bx; 1/2 ) \right). \]
Then given a translation invariant functional $H$
 on point sets in $\R^d$, define
\[ \Delta_{\bx} (A ; B) := H ( \Po''(\bx) \cap A ; B) -
 H ( \Po \cap A ; B) ; \]
this is the change in $H(\Po \cap A;B)$ when the Poisson
points in $Q(\bx; 1/2)$ are resampled.

\begin{lemma} Suppose $H$ is weakly stabilizing on $\RR$. Then for all $\bx \in \Z^d$,
there is a random variable $\Delta_{\bx} (\infty)$ such that for all $\bx \in \Z^d$,
\bea
\label{0219g}
  \Delta_{\bx} (A;B) \toas \Delta_{\bx} (\infty),
\eea
as $(A,B) \to \R^d$ through $\RR$. Moreover, if $H$ satisfies
(\ref{Pbmc}), then
\bea
\label{0219a}
 \sup_{ (A,B) \in \RR, \bx \in \Z^d} \Exp 
\left[ ( \Delta_{\bx} (A;B) )^4 \right] < \infty.
\eea
\end{lemma}
\proof
Set $C_0 = Q(\0; 1/2)$. By translation invariance, we need only consider
the case $\bx = \0$, and thus it suffices to prove that the variables
 $H ( \Po \cap A;B)
- H (\Po \cap A \setminus C_0 ; B)$ converge almost surely as $(A,B) \to \R^d$ through $\RR$.

The number $N$ of points of $\Po$ in $C_0$ is Poisson with parameter $1$. 
Let $\bV_1, \bV_2, \ldots, \bV_N$ be the points of $\Po \cap C_0$,
 taken in an order chosen uniformly
at random from the $N!$ possibilities. Then, provided $C_0 \subseteq A$,
\[ H ( \Po \cap A ; B) - H ( \Po \cap A \setminus C_0 ;B)
= \sum_{i=0}^{N-1} \delta_i (A;B), \]
where
\[ \delta_i (A;B)
:= H (( \Po \cap A \setminus C_0 ) \cup \{ \bV_1, \ldots, \bV_{i+1} \} 
;B) 
- H ( ( \Po \cap A \setminus C_0 ) 
\cup \{ \bV_1, \ldots, \bV_i \} 
;B) . \]
Since $N$ is a.s.~finite, it suffices
 to prove that each $\delta_i (A;B)$ converges almost surely
as $(A,B) \to \R^d$ through $\RR$. Let $\bU$ be a uniform random vector on
 $C_0$, independent of $\Po$. The distribution of the translated point 
process $- \bV_{i+1} + \{ \bV_1, \ldots, \bV_i \} \cup
( \Po \setminus C_0 )$ is the same as the conditional distribution of $\Po$ given that the
number of points in $- \bU + C_0 $ is
 equal to $i$, an event of strictly positive
probability. By assumption, this satisfies weak stabilization, which proves (\ref{0219g}).

Next we prove (\ref{0219a}). If $Q( \bx; 1/2) \cap A =
\emptyset$ then $\Delta_\bx (A;B)$ is
zero with probability 1.
By translation invariance, it suffices to consider the $\bx = \0$ case,
that is, to prove
\bea \label{0219h}
\sup_{(A,B) \in \RR : C_0 \cap A \neq \emptyset}
\Exp \left[ \left( \Delta_{\0} (A;B) \right)^4 \right] < \infty.
\eea
The proof of this now follows the proof of (3.4) of \cite{penyuk1}, but with
$\delta_i(A)$ replaced by $\delta_i(A;B)$ everywhere. 
$\square$ \\

\noindent
\textbf{Proof of Theorem \ref{0217c}.} Here we can assume, without
loss of generality, that $\QQ_n = \Po \cap R_n$. For $\bx \in
\Z^d$, let $\FF_\bx$ denote the $\sigma$-field generated by the
points of $\Po$ in $\cup_{\by \in \Z^d : \by \leq \bx} Q(\by;
1/2)$, where the order in the union is the lexicographic order on
$\Z^d$.

Let $R'_n$ be the set of points
 $\bx \in \Z^d$ such that $Q(\bx ;1/2) \cap R_n \neq \emptyset$.
 Let $k_n = \card ( R'_n )$. Then we
have that
\[ R_n \subseteq \bigcup_{\bx \in R'_n} Q ( \bx ; 1/2 ) \subseteq
R_n \cup \partial_1 (R_n) , \] so that
\[ | R_n | \leq k_n \leq | R_n | + | \partial_1 (R_n) | . \]
The vanishing relative boundary condition then implies that $k_n
 / n \to 1$ as $n \to \infty$.

Define the filtration $( \GG_0, \GG_1, \ldots, \GG_{k_n} )$ as
follows: let $\GG_0$ be the trivial $\sigma$-field, label the
elements of $R'_n$ in lexicographic order as $\bx_1, \ldots,
\bx_{k_n}$ and let $\GG_i = \FF_{\bx_i}$ for $1 \leq i \leq k_n$.
Then $H (\QQ_n ; S_n) - \Exp [ H (\QQ_n ; S_n ) ] =
\sum_{i=1}^{k_n} D_i$, where we set
\bea
 D_i = \Exp [ H (\QQ_n ; S_n) | \GG_i ] -
\Exp [ H (\QQ_n ; S_n) | \GG_{i-1} ] = \Exp [ -\Delta_{\bx_i} (
R_n ; S_n ) | \FF_{\bx_i} ] . 
\label{0819}
\eea
 By orthogonality of martingale
differences, $\Var [ H (\QQ_n ; S_n ) ] = \Exp \sum_{i=1}^{k_n}
D_i^2 $. By this fact, along with a CLT for martingale differences
(Theorem 2.3 of \cite{mcleish} or Theorem 2.10 of \cite{penbook}),
it suffices to prove the
conditions
\bea
\label{0217d} \sup_{n \geq 1} \Exp \left[ \max_{ 1
\leq i \leq k_n} \left\{ k_n^{-1/2} |D_i |\right\} ^2 \right] <
\infty,
\\ \label{0217e}
k_n^{-1/2} \max_{1 \leq i \leq
k_n} |D_i| \toP 0, \eea and for some $s^2 \geq 0$,
\bea
\label{0217f}
k_n^{-1} \sum_{i=1}^{k_n} D_i^2 \inL  s^2 .
\eea

Using \eq{0219a}, and the representation \eq{0819} for $D_i$,
we can verify \eq{0217d} and \eq{0217e} in just the same manner
as for the equivalent estimates (3.7) and (3.8) in \cite{penyuk1}.

We now prove (\ref{0217f}). By 
 (\ref{0219g}), for each $\bx \in \Z^d$ the variables $\Delta_{\bx}
(A;B)$ converge almost surely to a limit, denoted $\Delta_{\bx} (\infty)$,
as $(A,B) \to \R^d$ through $\RR$. For $\bx \in \Z^d$ and $(A,B) \in \RR$, let
\[ F_{\bx} (A;B) = \Exp [ \Delta_{\bx} (A;B) | \FF_{\bx} ]; ~~~
F_{\bx} = \Exp [ \Delta_{\bx} (\infty) | \FF_{\bx} ]. \]
Then $(F_{\bx} , \bx \in \Z^d)$ is a stationary family of random variables.
Set $s^2 = \Exp [ F^2_{\0}]$. We claim that the ergodic theorem
 implies
\bea \label{0322p}
 k_n^{-1} \sum_{\bx \in R'_n} F_{\bx}^2 \inL s^2 . \eea
The proof of this follows, with minor modifications,
 the proof of the corresponding result
(3.10) in \cite{penyuk1}.

We need to show that $F_{\bx}(R_n;S_n)^2$ approximates to $F_{\bx}^2$. We consider
$\bx$ at the origin $\0$. For any $(A,B) \in \RR$, by Cauchy-Schwarz,
\bea
 \Exp [ | F_{\0} (A;B)^2 - F_\0^2 | ] \leq \left(
\Exp[ (F_\0(A;B) + F_\0)^2 ] \right)^{1/2}
\left(
\Exp[ (F_\0(A;B) - F_\0)^2 ] \right)^{1/2}. 
\label{0801}
\eea
By the definition of $F_\0$ and the conditional Jensen inequality,
\bean
\Exp [ (F_\0(A;B) + F_\0 )^2] & = & \Exp \left[
\left( \Exp[ \Delta_\0 (A;B) + \Delta_\0 (\infty) | \FF_\0 ] \right)^2 \right] \\
& \leq & \Exp [ \Exp [ (\Delta_\0 (A;B) + \Delta_\0 (\infty))^2 | \FF_\0 ]]
\\
= \Exp [ (\Delta_\0 (A;B) + \Delta_\0 (\infty) )^2] , \eean
which is uniformly bounded by
%the alternative stabilization
(\ref{0219g})
and
(\ref{0219a}). Similarly,
\bea
\label{0325b}
\Exp [ (F_\0(A;B) - F_\0 )^2]  \leq
 \Exp [ (\Delta_\0 (A;B) - \Delta_\0 (\infty) )^2] ,
\eea
which is also uniformly bounded by (\ref{0219g}) and (\ref{0219a}).
For any $\RR$-valued sequence $(A_n,B_n)_{n \geq 1}$
with $\cup_{n \geq 1} \cap_{m \geq n} B_n = \R^d$,
the sequence $(\Delta_\0 (A_n ;B_n) - \Delta_\0 (\infty))^2$
tends to 0 almost surely by (\ref{0219g}), and is uniformly integrable by (\ref{0219a}), and therefore 
%(see \cite{durrett}, Chapter 4, Theorem 5.2) 
the expression (\ref{0325b})
tends to zero so that by (\ref{0801}),
 $\Exp[ |F_\0 (A_n;B_n)^2 - F_\0^2|] \to 0$.

Returning to the given sequence $(R_n,S_n)$, let $\eps >0$. By the vanishing relative boundary condition, we can choose $K_n$
so that $\lim_{n \to\infty} K_n= \infty$ and $| \partial_{K_n} S_n | \leq
\eps n$ for all $n$. 
Let $S'_n$ be the set of $\bx \in \Z^d$ such that $Q_{1/2}(\bx)$
has non-empty intersection with $S_n \setminus \partial_{K_n}(S_n)$.
Using the conclusion of the previous paragraph and translation invariance,
 it is not hard to deduce that
\bea
 \lim_{n \to \infty} \sup_{\bx \in S_n'} \Exp [ |F_{\bx} (R_n;S_n)^2
-F_{\bx}^2 | ] = 0 . 
\label{0727g}
\eea
Also, since we assume $|S_n | \sim n$ we have 
 $\card(S'_n) \geq |S_n| - \eps n \geq (1-2 \eps)n$ for large enough $n$. 
Using this with (\ref{0727g}),
 the uniform boundedness of $\Exp [ |F_{\bx} (R_n;S_n)^2
-F_{\bx}^2 | ]$ and the fact that $\eps$ can be taken arbitrarily small
in the above argument, it is routine to deduce that
\[ 
k_n^{-1} \sum_{\bx \in R'_n} ( F_{\bx} (R_n;S_n)^2 - F^2_{\bx} ) \inL 0,
\]
and therefore (\ref{0322p}) remains true with $F_{\bx}$ replaced by 
$F_{\bx}(R_n;S_n)$;
that is, (\ref{0217f}) holds and the proof of Theorem 
\ref{0217c} is complete. $\square$

\subsection{Proof of Theorem \ref{genclt}: the non-Poisson case}

In this section we complete the proof of Theorem \ref{genclt}. 
The first step
is to show that the conditions of Theorem \ref{genclt} imply those
of Theorem \ref{0217c}, as follows.

\begin{lemma}
\label{0326a}
If $H$ satisfies the uniform bounded moments condition (\ref{ubm}) and is
polynomially bounded, then $H$
 satisfies the Poisson bounded moments condition
(\ref{Pbmc}). \end{lemma}
\proof
The proof follows, with minor modifications, that of
Lemma 4.1 of \cite{penyuk1}.
$\square$ \\

It follows from Lemma \ref{0326a} that if $H$ satisfies the 
conditions of Theorem \ref
{genclt}, then Theorem \ref{0217c} applies and we have the Poisson parts of
Theorem \ref{genclt}. To de-Poissonize these limits we follow \cite{penyuk1}.
Define
\[ R_{m,n} := H ( \UU_{m+1,n} ; B) - H ( \UU_{m,n} ; B) .\]
We use the following coupling lemma.

\begin{lemma}
\label{0326b}
Suppose $H$ is strongly stabilizing. Let $\eps >0$. Then there
exists $\delta>0$ and $n_0 \geq 1$ such that for all $n \geq n_0$
and all $m, m' \in [(1-\delta)n, (1+\delta)n]$ with $m < m'$, there exists
a coupled family of variables $D, D', R, R'$ with the following
properties:
\begin{itemize}
\item[(i)] $D$ and $D'$ each have the same distribution as
$\Delta(\infty)$; \item[(ii)] $D$ and $D'$ are independent;
\item[(iii)] $(R,R')$ have the same joint distribution as
$(R_{m,n}, R_{m',n})$; \item[(iv)] $\Pr [ \{ D \neq R \} \cup \{
D' \neq R' \} ] < \eps$.
\end{itemize}
\end{lemma}
\proof
Since we assume $|S_n|/|R_n| \to 1$, the probability that
a random $d$-vector uniformly distributed over $R_n$ lies
in $S_n$ tends to 1 as $n \to \infty$. Using this fact 
the proof follows, with some minor modifications,
 that of the corresponding result in \cite{penyuk1},
Lemma 4.2. $\square$

\begin{lemma}
\label{0331q} Suppose $H$ is strongly stabilizing and satisfies
the uniform bounded moments condition (\ref{ubm}). Let $(h(n))_{n
\geq 1}$ be a sequence with $n^{-1} h(n) \to 0$ as $n \to \infty$.
Then \bea \label{0404p} & & \lim_{n \to \infty} \sup_{ |n-m| \leq
h(n)} \left|
\Exp R_{m,n} - \Exp \Delta ( \infty) \right| = 0; \\
\label{0404q} &  & \lim_{n \to \infty} \sup_{ n-h(n) \leq m < m'
\leq n+h(n)} \left|
\Exp R_{m,n}R_{m',n} - (\Exp \Delta ( \infty))^2 \right| = 0; \\
\label{0404r} &  & \lim_{n \to \infty} \sup_{ |n-m| \leq h(n)}
\Exp R_{m,n}^2 < \infty . \eea \end{lemma}

\proof The proof follows that of Lemma 4.3 of \cite{penyuk1}. 
$\square$ \\

\noindent \textbf{Proof of Theorem \ref{genclt}} Theorem
\ref{genclt} now follows in the same way as Theorem 2.1
in \cite{penyuk1},
replacing $H ( \; \cdot \;)$ with $H ( \; \cdot \; ; S_n)$. 
$\square$

\section{Proof of Theorem \ref{llnthm}: Laws of large numbers} \label{seclln}

We now derive our law of large
numbers for the total weight of the random MDSF on the unit
square.
We consider the general partial order
$\stackrel{\theta,\phi}{\preccurlyeq}$, for $0 \leq \theta < 2\pi$
and $ 0 < \phi \leq \pi$ or $\phi =2 \pi$.
 Recall that $\by \potp \bx$ if $\by \in
C_{\theta,\phi}(\bx)$, where $C_{\theta,\phi}(\bx)$ is the cone
formed by the rays at $\theta$ and $\theta+\phi$ measured
anticlockwise from the upwards vertical.

We consider the random point set 
$\X_n$,
the
binomial point process of $n$ independent uniformly distributed
points on $(0,1]^2$. However, the result 
 (\ref{0728e}) also holds 
(with virtually the same proof) if the points of $\X_n$
are uniformly distributed on an arbitrary convex set
in $\R^2 $ of unit area. If the points are 
distributed in 
$\R^2$ with a density function $f$ that has
convex support and is bounded away from 0 and infinity
on its support, then (\ref{0728e}) holds with a factor
of $\int_{\R^2} f(\bx)^{(2-\alpha)/2} \ud \bx$ introduced into the
right hand side (cf. eqn (2.9) of \cite{penyuk2}).

For the general partial order given by $\theta,\phi$ we apply 
Lemma \ref{llnpenyuk} to obtain a law of
large numbers for $\LL^\alpha ( \X_n )$. As a special case, we
thus obtain a law of large numbers under the partial order
$\preccurlyeq^*$  given
by $\theta=\phi=\pi/2$.
This method enables us to evaluate the limit explicitly, unlike
methods based on the subadditivity of the functional 
which may also be applicable here (see the remark at the end of this section).

In applying Lemma
\ref{llnpenyuk} to the MDSF functional, we 
take
the dimension $d$ in the lemma to be $2$,
and take   $f(\mathbf{x})$ (the underlying
probability  density function in the lemma)
to be 1 for 
$\mathbf{x} \in (0,1]^2$ and zero elsewhere.
We take $\xi(\mathbf{x} ; \X)$ to be
$d(\bx;\X)^\alpha$, where
$d(\bx;\X)$ is
 the distance from point $\mathbf{x}$ to its directed nearest
neighbour in $\X$ under $\potp$, if such a neighbour exists, or
zero otherwise. 
Thus in our case
\bea
 \xi(\mathbf{x};\X) = \left( d(\mathbf{x};\X) \right)^\alpha ~~~ 
{\rm with}~~~
d (\bx ; \X) := \min
\left\{ \| \bx - \by \| : \by \in \X \setminus \{ \bx \}, \by
\preccurlyeq \bx \right\} 
\label{0802}
\eea
with the convention that $\min \{ \} = 0$.
We need to show this choice of $\xi$ satisfies the conditions
of Lemma \ref{llnpenyuk}. As before, $\H_\tau$ denotes a homogeneous
Poisson process on $\R^d$ of intensity $\tau$, now with $d=2$.

\begin{lemma} \label{stabil} 
Let $\tau >0$.  Then $\xi$ is almost surely stabilizing on
$\mathcal{H}_{\tau}$, in the sense of (\ref{stabeq}),
with limit  
$\xi_\infty(\H_\tau) = (d(\0;\H_\tau))^\alpha$.
\end{lemma}
 \proof 
Let $R$ be the (random) distance from $\0$ to its directed nearest
neighbour in $\H_\tau$, i.e.~$R = d( \0 ; \H_\tau )$. 
%Then,
Since $\phi>0$ and $\tau>0$, we have $0< R < \infty$  almost surely.
 But then
for any $\ell >R$, we have $\xi(\0; (\H_\tau \cap B(\0;\ell)) \cup
{\cal A}) = R^\alpha$, for any finite ${\cal A} \subset \R^d \setminus
B(\0;\ell)$. Thus $\xi$  stabilizes on $\H_\tau$ with limit
$\xi_\infty(\H_\tau) =R^\alpha$.  $\square$ \\

Before proving that our choice of $\xi$ satisfies the
 moments condition for Lemma \ref{llnpenyuk},
we give a geometrical lemma.  For $B \subseteq \R^2$
with $B$ bounded, and for $\bx \in B$, 
write $\dist(\bx;\partial B)$ for $\sup\{r: B(\bx;r) \subseteq B\}$,
and for $s >0$, define the region 
\bea
\label{Atpdef} 
A_{\theta,\phi}(\mathbf{x},s;B) :=
B( \mathbf{x}; s ) \cap B \cap C_{\theta,\phi}(\bx). 
\eea
\begin{lemma}
\label{lem0727}
Let $B$ be a  convex bounded set in $\R^2$, and  let $\bx \in B$.
If $A_{\theta,\phi} (\bx,s;B) \cap \partial B(\bx;s) \neq \emptyset$,
and $s > \dist(\bx, \partial B)$, then 
$$
|A_{\theta,\phi}(\bx,s;B)| \geq   
s
\sin (\phi /2)  
 \dist(x,\partial B)  /2.
$$
\end{lemma}
\proof
The condition
 $A_{\theta,\phi} (\bx,s;B) \cap \partial B(\bx;s) \neq \emptyset$
says that there exists $\by \in B \cap C_{\theta,\phi}(\bx,s)$
with $\|\by - \bx\| = s$.
 The line segment $\bx \by$ is contained in
the cone $C_{\theta,\phi}(\bx)$; take a half-line ${\bf h}$ starting
from $\bx$, at an angle $\phi/2$ to the line segment $\bx \by$
and such that ${\bf h}$ is also contained in $ C_{\theta,\phi}(\bx)$.
Let $\bz$ be the point in $\bf h$ at a distance $\dist(\bx,\partial B)$
from $\bx$. Then the interior of the triangle $\bx \by \bz$ is entirely
contained in  $A_{\theta,\phi}(\bx,s)$, and has area
$s \sin (\phi /2)  \dist(x,\partial B)/2$. $\square$

\begin{lemma} \label{lem0k715b} Suppose 
$\alpha >0$.
Then  $\xi$ given by (\ref{0802}) satisfies the
moments condition (\ref{0k715a}) for any $p \in (1/\alpha,2 /\alpha]$.
\end{lemma} \proof 
Setting $R_n :=(0,n^{1/2}]^2$, we have
\bea
 \Exp \left[ \xi \left( n^{1/2}\mathbf{X}_1;n^{1/2} \X_n
\right)^p \right]  =  
\int_{R_n} \Exp \left[ \left( \xi
(\mathbf{x};n^{1/2} \X_{n-1}) \right)^p \right] \frac{\ud \mathbf{x}}{n}
%\nonumber
 . 
\label{0728}
\eea 
For $x \in R_n$ set $m(\bx) := \dist(\bx, \partial R_n)$.
Let us divide $R_n$ into three regions 
\bean
R_n(1) & : = & \{\bx \in R_n: m(\bx) \leq n^{-1/2} \}; ~~~~
R_n(2)  : =  \{\bx \in R_n:  m(\bx) > 1 \};
 \\
R_n(3) &  : = & \{\bx \in R_n: n^{-1/2} < m(\bx) \leq 1 \}.
\eean
%with $m(\bx) \leq n^{-1/2}$.
 For all $\bx \in R_n$, we have
$\xi(\bx;n^{1/2} \X_{n-1}) \leq (2n)^{\alpha/2}$, and hence,
since $R_n(1)$ has area at most 4, we can
bound the contribution to (\ref{0728}) from $\bx \in R_n(1)$ by 
\bea
\label{0728a}
\int_{\bx \in R_n(1)} \Exp \left[ \left( \xi
(\mathbf{x};n^{1/2}\X_{n-1}) \right)^p \right] \frac{\ud \mathbf{x}}{n}
\leq  4 n^{-1} (2n)^{p\alpha /2} = 2^{2+ p\alpha/2}  n^{(p \alpha -2)/2}, 
\eea
which is bounded provided $p \alpha \leq 2$.
 
Now, for $\mathbf{x} \in R_n$, with $A_{\theta,\phi}(\cdot)$ defined
at (\ref{Atpdef}),  
we have
\bea 
\Pr \left[  d
(\mathbf{x}; n^{1/2} \X_{n-1})  > s \right] & \leq & 
\Pr \left[ n^{1/2} \X_{n-1}
 \cap A_{\theta,\phi}(\bx,s;R_n) = \emptyset \right]
\nonumber \\
 & = & 
\left(1 - \frac{|A_{\theta,\phi}(\bx,s;R_n)| }{n} \right)^{n-1}
\nonumber \\
& \leq &  \exp( 1 - |A_{\theta,\phi}(\bx,s;R_n)| ),
\label{0728b} \eea 
since $|A_{\theta,\phi}(\bx,s;R_n)|\leq n$.
For $\bx \in R_n$ and $s>m(\bx)$, by Lemma \ref{lem0727} we have
 $$
\left| A_{\theta,\phi}(\mathbf{x},s;R_n)
 \right| \geq \sin(\phi/2 ) s m(\bx)/2 ~~~{\rm if} ~~ 
A_{\theta,\phi} (\bx,s;R_n) \cap \partial B(\bx;s) 
\neq \emptyset,
$$
and also
$$
\Pr [ d(\bx;n^{1/2} \X_{n-1}) > s ] = 0 ~~~{\rm if}~~~
A_{\theta,\phi} (\bx,s;R_n) \cap \partial B(\bx;s) 
= \emptyset.
$$
 For $s \leq m(\bx)$, 
we have that 
$
 \left|
A_{\theta,\phi}(\mathbf{x},s;R_n) \right| =
\frac{\phi}{2} s^2 \geq \sin (\phi/2)s^2. 
$ 
Combining these observations and (\ref{0728b}),
 we obtain for all $\bx \in R_n$  and $s >0$ that 
\bean 
\Pr \left[  d
(\mathbf{x}; n^{1/2} \X_{n-1})  > s \right] & \leq & 
\exp \left( 1 - \sin (\phi/2) s \min (s, m(\bx) )/2 
\right), ~~~ \bx \in R_n.
\eean
Setting $c =(1/2) \sin(\phi/2)$, we therefore have for $\mathbf{x} \in R_n$
that
\bea
  \Exp \left[  \xi (\mathbf{x};n^{1/2} \X_{n-1})^p
\right] 
%\!\! 
& = & 
%\!\! 
\int_0^\infty  
\Pr  \left[  \xi (\mathbf{x}; n^{1/2} \X_{n-1})^p > r \right]
 \ud r 
\nonumber \\
& = & \int_0^\infty \Pr \left[ d(\bx; n^{1/2} \X_{n-1})
> r^{1/(\alpha p)} \right] \ud r
\nonumber \\
%\!\! 
& \leq &
% \!\! 
\int_0^{m(\bx)^{\alpha p}} \ud r 
 \exp {  \left(1 - c r^{2/(\alpha p)} \right) } 
\nonumber \\
&& +
 \int_{m(\bx)^{\alpha p }}^\infty \!\! \ud r
\exp { \left(1 - c m(\bx)r^{1/(\alpha p)} \right) }
\nonumber \\
& = & O(1) + \int_{m(\bx)^2}^\infty 
e^{1-cu} \alpha p u^{\alpha p-1} m(\bx)^{-p \alpha} \ud u
\nonumber \\
& = & O(1) + O(m(\bx)^{-\alpha p}).
\label{0728d}
\eea 
For $\bx \in R_n(2)$, this bound is $O(1)$, and the area
of $R_n(2)$ is less than $n$, so that the contribution
to (\ref{0728}) from $R_n(2)$ satisfies
\bea
\label{0728c}
\limsup_{n \to \infty}
\int_{R_n(2)} \Exp \left[ \left( \xi
(\mathbf{x};n^{1/2} \X_{n-1}) \right)^p \right] 
\frac{\ud \mathbf{x}}{n}
%\ud \bx 
< \infty.
\eea
Finally, by (\ref{0728d}), there is a constant $c'$ such 
that if $\alpha p >1$, 
the contribution to (\ref{0728}) from $R_n(3)$ satisfies
% is bounded by a constant times
\bean
\int_{R_n(3)} \Exp \left[ \left( \xi
(\mathbf{x};n^{1/2} \X_{n-1}) \right)^p \right] \frac{\ud \mathbf{x}}{n}
&\leq & 
  c' n^{-1/2} \int_{y=n^{-1/2}}^1 y^{-\alpha p} \ud y 
\\
& \leq & \left( \frac{c' n^{-1/2}}{\alpha p-1} \right) n^{(\alpha p -1)/2}
\eean
which is bounded provided  $\alpha p \leq 2$.
Combined with the bounds in (\ref{0728a}) and (\ref{0728c}),
this shows that the expression (\ref{0728}) 
is uniformly bounded, provided $1 < \alpha p \leq 2$.
$\square$ \\

Following notation from Section \ref{subsecgenclt},
for $k \in \N$, and for $a<b$ and $c<d$
 let $\UU_{k,(a,b] \times (c,d]}$ denote the point process 
consisting of
$k$ independent  random vectors uniformly distributed
on the rectangle $(a,b] \times (c,d]$.
Before proceeding further,
we recall 
that
if $M(\X)$ denotes the number of minimal elements (under the ordering
$\postar$) of a point 
set $\X \subset \R^2 $, then 
\bea
\Exp [ M(\UU_{k,(a,b]\times (c,d]}) ] =
\Exp [ M(\X_k) ]
 = 1 + (1/2) + \cdots + (1/k) \leq 1 + \log k. 
\label{harmonicbd}
\eea
The first equality in (\ref{harmonicbd})  comes
from some obvious scaling which shows that the
distribution of 
$ M(\UU_{k,(a,b]\times (c,d]}) $ does not depend on $a,b,c,d$.
For  the second equality 
in (\ref{harmonicbd}),
see \cite{BNS} or the proof of 
Theorem 1.1(a) of \cite{bhattroy2002}.  \\

\noindent \textbf{Proof of Theorem \ref{llnthm}.} 
Suppose $\alpha <2$,
and set $f(\cdot)$ to be the indicator of the unit square
$(0,1]^2$. 
 By Lemmas
\ref{stabil} and \ref{lem0k715b}, our functional $\xi$,
given at (\ref{0802}), satisfies the
conditions of Lemma \ref{llnpenyuk} with $p= 2/\alpha$
and $q=1$, with this choice of $f$. 
 So by Lemma \ref{llnpenyuk},
we have that
\bea
n^{(\alpha/2)-1} \LL^\alpha(\X_n) = 
  n^{-1} 
\sum_{\bx \in \X_n} \xi (n^{1/2}\mathbf{x};n^{1/2}\X_n) 
 \nonumber \\
\inL
\int_{\R^2} \Exp \left[ \xi_\infty ( \H_{f(\mathbf{x})} ) \right]
f( \mathbf{x} ) \ud \mathbf{x}  = \Exp \xi_\infty(\H_1).
\label{0728f}
 \eea
Since the disk sector $C_{\theta,\phi}(\mathbf{x}) \cap B(\bx;r)$
has area $(\phi/2) r^2$, by Lemma \ref{stabil} we have
 \bean \Pr\left[ \xi_{\infty} (
\H_1 ) >s \right] & = & \Pr \left[ \H_1 \cap
C_{\theta,\phi}(\mathbf{0}) \cap B(\0;s^{1/\alpha}) = \emptyset \right] =
\exp{\left(-(\phi/2) s^{2/\alpha}\right)}.
\eean
Hence, the limit in (\ref{0728f}) is
 \[
  \Exp \left[ \xi_{\infty} ( \H_1 ) \right]  =
\int_0^{\infty} \Pr \left[ \xi_\infty \left( \H_1 \right)
> s \right] \ud s
= \alpha 2^{(\alpha-2)/2} \phi^{-\alpha/2} \Gamma( \alpha/2 ), 
\]
and this gives us (\ref{0728e}). Finally, in the case where $\potp=\postar$,
(\ref{0728e}) remains true when $\X_n$ is replaced by $\X_n^0$, since
\bea
\label{0806a}
 \Exp [n^{(\alpha/2)-1} | \LL^\alpha (\X_n^0) - \LL^\alpha (\X_n)| ]
\leq 2^{\alpha/2} n^{(\alpha/2)-1} \Exp [ M(\X_n)] ,
\eea
where $M(\X_n)$ denotes the number of minimal elements of $\X_n$. By (\ref{harmonicbd}),
$\Exp[ M(\X_n)] \leq 1+\log n$, and hence the right hand side of (\ref{0806a})
tends to 0 as $n \to \infty$ for $0<\alpha<2$. This gives us (\ref{0728e})
with $\X_n^0$ under $\postar$.   
 $\square$ \\

\noindent \textbf{Remark.} A law of large numbers for Euclidean
functionals of many random geometric structures can be treated by
the boundary functional approach of Yukich~\cite{yukbook}. It
can be shown that the MDSF satisfies some, but possibly not all, of the
appropriate conditions that would allow this approach to be
successful. The MDSF functional is subadditive, its corresponding
boundary functional is superadditive, and the functional and its
boundary functional are sufficiently `close in mean'. However, it
is not clear that the functional is `smooth', since
the degree of the graph is not bounded.

\section{Central limit theorem away from the boundary} 
\label{ltot}

While it should be possible to adapt the 
argument of the present section to more general partial orders,
from now on we take the partial order $\preccurlyeq$ 
on $\R^2$ to be $\postar$.
For each $n$, define the region
 $S_{0,n} := (n^{\eps -1/2},1]^2$, 
where $\eps \in (0,1/2)$ is a small
 constant to be chosen later.  
In this section, 
 we use the general central limit theorems
of Section \ref{subsecgenclt} to
 demonstrate a central limit theorem for the contribution to the
 total weight
of the MDSF, under $\postar$, from edges away from the
 boundary, that is from points in the
region $S_{0,n}$.

Given $\alpha >0$,
consider the MDSF total weight functional $H = \LL^\alpha$ on
point sets in $\R^2$. 
For $\bx \in \X$, let the directed nearest
neighbour distance $d( \bx ; \X)$ and
the corresponding $\alpha$-weighted functional
 $\xi(\bx;\X)$
be  given by (\ref{0802}), where now we take $\preccurlyeq$ to be
$\postar$.
For $R \subseteq \R^2$, set
\bea
 \LL^\alpha ( \X;R ) = \sum_{\bx \in \X\cap R} \xi (\bx ; \X) ,
\label{0714a}
\eea
and set $\LL^\alpha(\X):= \LL^\alpha(\X;\R^2)$.

Let $\X_n$ be the binomial point process of $n$ i.i.d.~uniform
random vectors on $(0,1]^2$, and let $\Po_n$ be the homogeneous
Poisson process of intensity $n$ on $(0,1]^2$. The main result of
this section is the following.

\begin{theorem} \label{CLT}
Suppose that $\alpha>0$ and the partial order is $\postar$.
 Then there exist constants $0<t_\alpha \leq
s_\alpha$, not depending on the choice of $\eps$,
such that, as $n \to \infty$,
\begin{itemize}  \item[(i)] $n^{\alpha-1} \Var \left[
\LL^\alpha \left( \X_n ; S_{0,n} \right) \right] \to t_\alpha^2$;
\item[(ii)]
$n^{(\alpha-1)/2} \tLalph
 \left( \X_n ; S_{0,n} \right) 
 \tod \NN \left( 0, t_\alpha^2 \right)$;
\item[(iii)] $n^{\alpha-1}
\Var \left[ \LL^\alpha \left( \Po_n ; S_{0,n} \right) \right]
\to s_\alpha^2$;
\item[(iv)] $n^{(\alpha-1)/2} 
\tLalph \left( \Po_n ; S_{0,n} \right) 
\tod \NN
\left( 0, s_\alpha^2 \right)$.
\end{itemize} \end{theorem}

The following corollary states that 
Theorem \ref{CLT} remains true in the rooted
cases too, i.e.~with $\X_n$ replaced by $\X^0_n$ and
$\Po_n$ replaced by $\Po^0_n$.

\begin{corollary} \label{0804c}
Suppose that $\alpha>0$ and the partial order is $\postar$.
 Then, with $t_\alpha$, $s_\alpha$
as given in Theorem \ref{CLT}, we have that as $n \to \infty$,
\begin{itemize} \item[(i)] $n^{\alpha-1} \Var \left[
\LL^\alpha \left( \X_n^0 ; S_{0,n} \right) \right] \to t_\alpha^2$;
\item[(ii)]
$n^{(\alpha-1)/2} \tLalph
 \left( \X_n^0 ; S_{0,n} \right) \tod \NN \left( 0, t_\alpha^2 \right)$;
\item[(iii)] $n^{\alpha-1}
\Var \left[ \LL^\alpha \left( \Po_n^0 ; S_{0,n} \right) \right]
\to s_\alpha^2$;
\item[(iv)] $n^{(\alpha-1)/2} 
\tLalph \left( \Po_n^0 ; S_{0,n} \right) 
\tod \NN
\left( 0, s_\alpha^2 \right)$.
\end{itemize} \end{corollary}
\proof
For each region $R \subseteq [0,1]^2$ and
point set ${\cal S} \subset [0,1]^2$
with $\0 \in {\cal S}$, let $\LL_0^\alpha ( {\cal S} ; R)$ denote
the total weight of the edges incident to $\0$
in the MDST on ${\cal S}$ from points
in $R$.  Then
 $\Lalph ( \Po_n^0 ; S_{0,n} ) $ equals
$\Lalph ( \Po_n ; S_{0,n} ) + \LL_0^\alpha ( \Po_n^0 ; S_{0,n} )$,
so that
\bea
 \Var [   \Lalph ( \Po_n^0 ; S_{0,n} ) ]  - 
 \Var [   \Lalph ( \Po_n ; S_{0,n} ) ]
= 2 \Cov [ \Lalph ( \Po_n ; S_{0,n} ) , \LL_0^\alpha ( \Po_n^0 ; S_{0,n} ) ]
\nonumber \\
+
 \Var[ \LL_0^\alpha ( \Po_n^0 ; S_{0,n} ) ]. 
\label{0817a}
\eea

 Let $N_n$ denote the number of points of $\Po_n$, and 
let $E_n$ denote the event
that at least one point of $\Po_n \cap S_{0,n}$ is
joined to $\0$ in the MDST on $\Po_n^0$. Then 
\bean
%\label{0804i}
\Pr[E_n] \leq \Pr \left[ (0,n^{\eps-1/2}]^2 \cap \Po_n =
\emptyset \right] = \exp (-n^{2\eps}) ,
\eean
and $\LL_0^\alpha(\Po_n^0; S_{0,n}) \leq 2^{\alpha/2} N_n \1_{E_n} $.
Thus by the Cauchy-Schwarz inequality, 
%from (\ref{0804j}) and (\ref{0804i}) 
for some finite constant $C$
we have
\bea
 \Var \left[ \LL_0^\alpha ( \Po_n^0 ; S_{0,n} ) \right]
 \leq \Exp \left[ 
%\left|
 \LL_0^\alpha ( \Po_n^0 ; S_{0,n} )^2 
%\Lalph ( \Po_n ; S_{0,n} ) \right|
 \right]
\leq 
%2^{\alpha}
C 
%n^{\alpha+3}
n^2  \exp (-n^{2\eps}/2) ,
\label{0817b}
\eea
and combining this with (\ref{0817a}),
Theorem \ref{CLT} (iii) and
 the Cauchy-Schwarz inequality 
shows that
$$
n^{\alpha-1} 
( \Var
\left[   \Lalph ( \Po_n^0 ; S_{0,n} ) \right]  - 
 \Var \left[   \Lalph ( \Po_n ; S_{0,n} ) \right] )
\to 0,
$$
so that from Theorem \ref{CLT} (iii) we obtain
the corresponding rooted result (iii).
Also, since (\ref{0817b}) implies $n^{\alpha-1} 
 \Var \left[ \LL_0^\alpha ( \Po_n^0 ; S_{0,n} ) \right]$ tends to zero,
 from Theorem \ref{CLT} (iv) and Slutsky's theorem we obtain
 the corresponding
rooted result (iv).

The binomial results (i) and (ii) follow in the same manner as above,
 with slight
modifications. $\square$ \\

To prove Theorem \ref{CLT}, we demonstrate that our functional
 $\LL^\alpha$ satisfies
suitable versions of the conditions of Theorem \ref{genclt} and
Corollary \ref{cltcor}. First, we see that $\LL^\alpha$ is polynomially
bounded (see (\ref{poly})), since
\[ \LL^\alpha ( \X ; B) \leq \left( \diam ( \X ) \right)^\alpha
\card ( \X ) .\] Also, $\LL^\alpha$ is homogeneous of order $\alpha$.

\begin{lemma} \label{stab}
$\LL^\alpha$ is strongly stabilizing, in the sense of Definition
\ref{sstabdef}.
\end{lemma}
\proof
 To prove stabilization it is sufficient to show that there
exists an almost surely finite random variable $R$, the radius of
stabilization, such that the add one cost is unaffected by changes
in the configuration at a distance greater than $R$ from the added
point. We show that there exists such an $R$.

\begin{figure}[h]
\centering
\input{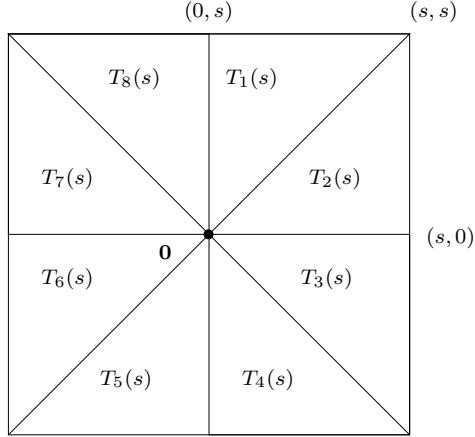}
\caption{The triangles $T_1(s), \ldots, T_8(s)$, $s>0$. }
\label{fig9}
\end{figure}

For $s>0$ construct eight disjoint triangles
$T_j (s), 1 \leq j \leq 8$, by splitting the square $Q(\0;s)$
into eight triangles via drawing
in the diagonals of the square and the $x$ and $y$ axes.
Label the triangle with vertices $(0,0),(0,s),(s,s)$ as
$T_1 (s)$ and then label increasingly in a clockwise
manner. See Figure \ref{fig9}. Note that $T_j (t) \subset T_j (s) $ for $t<s$.
Let the random variable $S$ be the minimum $s$ such that the
triangles $T_j (s), 1 \leq j \leq 8$, each contain at least one
point of
$\Po$. Then $S$ is almost surely finite.

We claim that $R=3S$ is a radius of stabilization for $\LL^\alpha
$, that is any points at distance $d \geq 3S$ from the origin have
no impact on the set of added or removed edges when a point is
inserted at the origin.

First, $\0$ can have no point
at a distance of at least $3S$ away as its directed nearest neighbour,
since there will be points in $T_5$ and $T_6$ within a distance of at most
$\sqrt{2} S$ of $\0$.

We now need to show that no point at a distance at least $3S$ from $\0$
can have the origin as its directed nearest neighbour. Clearly, for the
partial order $\postar$, we need only consider
points in the region $(0,\infty)^2$.

Consider a point $(x,y)$ in the first quadrant, such that $\| (x,y)\| \geq 3S$.
Consider the disk sector
\[ D_{(x,y)} := B \left( (x,y), \|(x,y)\| \right) \cap \left\{ {\bf w }: {\bf w}
\postar (x,y) \right\} .\]
We aim to show that given any $(x,y)$ of the above form,
at least one of the $T_j(S)$, $j=1,\ldots,8$, is contained in $D_{(x,y)}$, which
implies that the origin cannot be the directed nearest neighbour of $(x,y)$.
To demonstrate this, we show that given such an $(x,y)$, $D_{(x,y)}$ contains all
three vertices of at least one of the $T_j(S)$.

First suppose $x>S$, $y>S$. Then we have that
$T_1(S)$ and $T_2(S)$ are in $D_{(x,y)}$, since we have, for
example,
\bean \| (x,y) -\0 \|^2
- \| (x,y) - (0,S) \|^2
& = & \left( x^2 + y^2 \right) - \left( x^2 +(y-S)^2 \right)
\\ & = & S (2y-S) > 0 .
\eean
By symmetry, the only other situation we need consider is when $0<x \leq S$. Then $y^2
\geq 9S^2 - x^2 \geq 8S^2$, so $y \geq 2\sqrt{2} S$. Then we have that
$T_8(S)$ is in $D_{(x,y)}$,
since
\bean
\| (x,y) -\0 \|^2
- \| (x,y) - (-S,S) \|^2
& = & \left( x^2 + y^2 \right) - \left( (x+S)^2 +(y-S)^2 \right)
\\ & = & 2S(y-x-S) \geq 4S^2 (\sqrt{2} -1) > 0. \; 
\eean
This completes the proof. 
$\square$

\begin{lemma} \label{nondeg}
The distribution of $\Delta (\infty )$ is non-degenerate.
\end{lemma}
\proof
 We demonstrate the existence of two
 configurations that occur with strictly positive probability
 and give rise to different values for $\Delta(\infty)$.
Note  that
adding a point at the origin causes some new edges to be formed
(namely those incident to the origin), and the possible deletion
of some edges (namely the edges from points which have the origin as
their directed nearest neighbour after its insertion).

Let $\eta >0$, with $\eta < 1/3$. Later we shall impose further conditions
on $\eta$.
Again we refer to the construction in Figure \ref{fig9}. Let $E_1$
denote the event that for
each $i$, $1\leq i \leq 8$,
there is a single point of $\Po$, denoted ${\bf W}_i$, in
each of $T_i (\eta)$,
and that there are no other points in $[-1,1]^2$.
Suppose that $E_1$ occurs.
Then, on addition of the origin, the only edges that
can possibly  be removed are those from ${\bf W}_1$ and from ${\bf W}_2$
(see the proof of Lemma \ref{stab}). These removed edges
have length at most $\eta \sqrt{8}$, and hence
\bea
\Delta \geq - 2( \eta \sqrt{8})^\alpha := \delta_1, ~~~~{\rm on}~~ E_1.
\label{0629a}
\eea

\begin{figure}[h]
\centering
\input{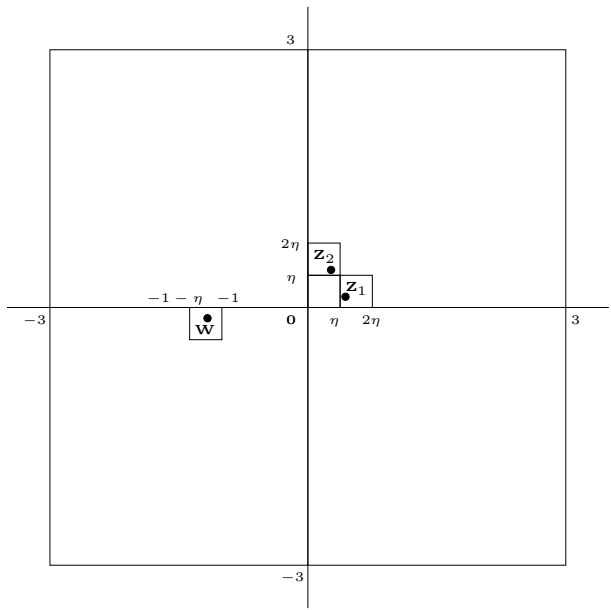}
\caption{A possible configuration for event $E_2$. }
\label{fige2}
\end{figure}

Now let $E_2$ denote the event that there is a single point of $\Po$,
denoted ${\bf Z}_1$,
 in the square $(\eta,2\eta)\times(0,\eta)$,
a single point denoted ${\bf Z}_2$ in the square
 $(0,\eta)\times(\eta,2 \eta)$,
 a single point denoted ${\bf W}$ in the square $(-1-\eta,-1) \times (-\eta,0)$, and no other point in $[-3,3]^2$. See Figure \ref{fige2}.

%(A pic might help here. Maybe not
%needed)
%\Comment{pic?}
Suppose that $E_2$ occurs. Now, on  addition of the origin,
an edge of length at most $1+2\eta$ is added from the origin
to ${\bf W}$. On the other hand, for $i=1,2$ the edge
from ${\bf Z}_i$ to ${\bf W}$ (of length at least 1) is replaced
by an edge from ${\bf Z}_i$ to the origin (of length at most
$3 \eta$). It is also possible that some other edges from points
outside $[-3,3]^2$ are replaced by shorter edges from these
points to the origin. Combining the effect of all these
additions and
replacements of edges, we find that
\bea
\Delta \leq (1+ 2 \eta)^\alpha + 2((3 \eta)^\alpha - 1) := \delta_2,
~~~~{\rm on} ~~ E_2.
\eea
Given $\alpha$,
by taking $\eta$ small enough we can arrange that $\delta_1 > -1/4$
and $\delta_2 < -3/4$. With such a choice of $\eta $, events
$E_1$ and $E_2$ both have strictly positive probability
which shows that the distribution of $\Delta$ is non-degenerate.
$\square $ \\

For the next lemma, we set $R_0 := (0,1]^2$, 
recalling that $S_{0,n}:= (n^{\eps-1/2},1]^2$ throughout this section,
and let
$\RR_0$ be as defined just before Corollary \ref{cltcor}.

\begin{lemma} \label{ubmlem} $\LL^\alpha$
satisfies the uniform bounded moments condition (\ref{ubm}) on $\RR_0$.
\end{lemma}
\proof Choose some $(A,B) \in \RR_0$ such that
 $\0 \in A$, i.e., such that for some $n \in \N$ the set
$A$ is a translate of 
$(0,n^{1/2}]^2$
% suitably translated so that it
 containing the origin and $B$
is the corresponding translate of $n^{1/2}S_{0,n}=(n^\eps,n^{1/2}]^2$.
Note that $|A|=n$,
and choose $m \in [n/2,3n/2]$. 

Denote the $m$ independent random vectors on $A$ comprising
$\UU_{m,A}$ by $\bV_1, \ldots ,\bV_m$.
For contributions to
$\Delta(\UU_{m,A};B)$
 we are only interested in edges
from points in the region $B$ away from the boundary of $A$, although the
origin can be inserted anywhere in $A$. Contributions to
%$\LL^\alpha ( \; \cdot \; ; B)$
$\Delta(\UU_{m,A};B)$
 come from the edges that are added
or deleted on the addition of $\0$.
We split $\Delta (\UU_{m,A};B)$ into two parts: the positive
contribution from added edges, $\Delta^+ (\UU_{m,A};B)$, and the
negative contribution, $\Delta^- (\UU_{m,A};B)$, from removed
edges.

 By construction of the MDSF, the added edges are those that
 have $\0$ as an end-point after it has been inserted.
Thus an upper bound on $\Delta^+ (\UU_{m,A} ;B)$ is $L_{\max}^{\alpha}
\delta(\0) + L_0^\alpha$, where $L_{\max}$ is the length of the
longest edge  from a point of $\UU_{m,A} \cap B$ to $\0$,  
and $\delta(\0)$ is the number of such edges
(or zero if no such edge exists), 
and $L_0$ is the length of the edge from $\0$,
or zero if no such edge exists.

For $\bw \in A $ and $ \bx \in B$, %  in $\R^2 $,
 with $\bw \postar \bx$, define the region
$$
 R(\bw,\bx) :=
\{ {\bf y} \in A : {\bf y} \postar
\bx, \|{\bf y}-\bx\| \leq \|\bw - \bx\| \}  .
$$
Since points in $B$ are distant at least 1 from the lower 
or left boundary of $A$,
by Lemma \ref{lem0727} there exists a constant
%(depending on $\eps$)
$0<C<\infty$ such that
\bea
| R(\bw, \bx )| \geq C \| \bx  - \bw \|,
{\rm ~for~ all~} \bw \in A, \bx \in B ~{\rm with}~
\bw \postar \bx
~{\rm and}~ \| \bx -\bw \| \geq 1.
\label{0719d}
\eea
Suppose there is a point at $\bx$ with $\0 \postar \bx$.
 Then, the probability of
the event $E(\bx)$ that $\bx$ is joined to  the origin in the MDSF
on $\UU_{m,A} \cup \{\0\}$ is
\bea
 \Pr[E(\bx)] & = & \Pr[ R(\0,\bx) \textrm{ empty} ]
= \left( 1-\frac{|R(\0,\bx)|}{|A|}
\right)^{m-1} \!\!\!\! 
\nonumber \\
& \leq &
% \left( 1-\frac{|R(\0,\bx)|}{n}
%\right)^{n/2-1} \!\!\!\! \leq \exp {(-|R(\0,\bx)|/3)} , 
 \exp \left( (1-m) \left(\frac{|R(\0,\bx)|}{n}\right) \right) 
% \!\!\!\! 
\leq \exp {(1-|R(\0,\bx)|/2)} , 
\label{0719a}
\eea
since $m \geq n/2$ and $|R(\0,\bx)| \leq n$.

We have that $ L_{\max}^\alpha \delta(\0) \leq
\max_{i=1,\ldots,m} W_i , $ where \[ W_i =
\|\bV_i\|^\alpha \; \card( B( \0; \| \bV_i\|) \cap \UU_{m,A} \cap \{
{\bf y} : \0 \postar {\bf y} \}
 ) \; {\bf
1} { \{ \bV_i \textrm{ joined to } \0 \textrm{ and } \bV_i \in B\} } .
\]
Let $N(\bx)$ denote the number of points of $\UU_{m-1,A}$ in
$B(\0; \|\bx\|) \cap  \{ {\bf y} : \0 \preccurlyeq {\bf y} \}$. Then
we obtain
\bean
 \Exp[ L_{\max}^{4 \alpha} \delta(\0)^4] \leq \Exp
\sum_{i=1}^m W_i^4 = m \int_{B} \|\bx\|^{4\alpha} \Exp[
(N(\bx)+1)^{4} {\bf 1} { \{ E(\bx) \}} ] \frac{\ud \bx}{|A|} .
\eean
 By the Cauchy-Schwarz inequality
and the fact that $m \leq 3|A|/2$ by assumption, 
\bea
 \Exp
[L_{\max}^{4\alpha} \delta(\0)^4] \leq \frac{3}{2} \int_{B}
\|\bx\|^{4\alpha} (\Exp [(N(\bx) +1)^{8} ])^{1/2} \Pr
[E(\bx)]^{1/2} \ud \bx .
\label{0719b}
\eea
The mean of $N(\bx)$ is bounded by a constant times $\|\bx\|^2$ so
$\Exp [(N(\bx) +1)^{8}] = O(\max(\|\bx\|^{16},1))$.
This follows from the binomial moment generating function for
$\mathrm{Bin}(n,p)$, from which we have
for $\beta >0$ that $\Exp[X^\beta] \leq
k_1 ( \Exp[X])^\beta$ if $pn>1$ and $\Exp[X^\beta] \leq k_2
 \Exp[X]$ if $pn<1$, for some constants $k_1,k_2 >0$.

Combined with (\ref{0719d}), (\ref{0719a}) and (\ref{0719b}), 
this shows that $\Exp [L_{\max}^{4\alpha} \delta(\0)^4] $ is
bounded by a constant times
\[
\int_{\bx \in B: \|\bx\| \geq 1} \|\bx\|^{4 \alpha + 8}
 \exp{ \left( -C \|\bx\|/4 \right)} \ud
\bx + \int_{\bx \in B: \|\bx\|\leq 1} \|\bx\|^{4 \alpha}\ud \bx,
\]
 which is
bounded by a constant that does not depend on
the choice of $(A,B)$.

We need to consider $L_0$ only when $\0 \in B$.
For $\bx \in \R^2$ with  $\bx \postar \0$, 
let $E'(\bx)$ denote the event that $R(\bx,\0)$ is empty
(i.e., contains no point of $\UU_{m-1,A}$).
By (\ref{0719d}) and (\ref{0719a}), 
for $\0 \in B$ we have
\bean
\Exp [ L_0^{4\alpha} ]
%& = & \int_0^\infty \Pr[ L_0 > r^{1/(4\alpha)} ] \ud r
%\\
 & \leq & m
\int_{\bx \in A: \bx \postar \0 }
 \|\bx\|^{4 \alpha}
P[E'(\bx)]\frac{\ud \bx}{|A|}
%\int_0^{1/2} \int_0^{1/2}  \| x \|^4 \exp{ (-\pi \| x \|^2/12 ) } \ud x \ud y
\\
& \leq & \frac{3}{2}
\left[
 \int_{ \bx \in A: \bx \postar \0,
\|\bx\| \geq 1 }
\!\!\!\!
 \|\bx\|^{4 \alpha} \exp(1-C\|\bx\|/2) \ud \bx
  + \int_{ \bx \in A: \bx \postar \0,
\|\bx\| \leq 1 } \!\!\!\! \|\bx\|^{4 \alpha}  \ud \bx
\right]
\eean
which is bounded by a constant.
 Thus $\Delta^+(\UU_{m,A};B)$ has bounded
fourth moment.

Now consider the set of deleted edges.
As at (\ref{0802}), 
let $d(\bx;\UU_{m,A})$ denote the distance from $\bx$ to its
directed nearest neighbour in $\UU_{m,A}$, or zero
if no such point exists. Again use $E(\bx)$ for
the event that $\bx$ becomes joined to $\0$ on the
addition of the origin, and let $E''(\bV_i):= E(\bV_i) \cap
 \{\bV_i \in B\}$. 
%Let $\bU_{1},\ldots,\bU_{n}$
%be the independent random $d$-vectors uniformly
%distributed on $A$. %comprising
Then
\bea
\Exp[ \Delta^- (\UU_{m,A};B)^4 ]  = 
\sum_{i=1}^m \sum_{j=1}^m \sum_{k=1}^m \sum_{\ell=1}^m
\Exp[ d(\bV_i;\UU_{m,A})^\alpha
 d(\bV_j;\UU_{m,A})^\alpha
\nonumber \\
 \times d(\bV_k;\UU_{m,A})^\alpha
 d(\bV_\ell;\UU_{m,A})^\alpha
{\bf 1}\{
E''(\bV_i) \cap
E''(\bV_j) \cap
E''(\bV_k) \cap
E''(\bV_\ell)\}].  \label{0728h}
\eea
For $i,j,k,\ell$ distinct, the $(i,j,k,\ell)$th term of
%  expectation in this sum
(\ref{0728h}) is bounded by
\bea
\int_{B} \int_{B} \int_B \int_B 
\frac{\ud  \bw}{n} 
\frac{\ud  \bx}{n}
\frac{\ud  \by}{n}
\frac{\ud  \bz}{n}
 \Exp[ d_{m-4}(\bw)^\alpha
 d_{m-4}(\bx)^\alpha
 d_{m-4}(\by)^\alpha
 d_{m-4}(\bz)^\alpha
\nonumber \\
\times
{\bf 1}\{
E_{m-4}(\bw) \cap
E_{m-4}(\bx) \cap
E_{m-4}(\by) \cap
E_{m-4}(\bz)\}], 
\label{0728i}
\eea 
where $d_{m-4}(\bx):=d(\bx, 
\UU_{m-4,A}\cup\{\bw,\bx,\by,\bz\})$ (using the notation of (\ref{0802})),
 and $E_{m-4}(\bx)$ is the event
that $\0$ is the directed nearest neighbour of $\bx$ in
the set $\UU_{m-4,A} \cup \{\0,\bx\}$.

Let $I_{m-4}(\bx)$ denote the indicator variable of the
event that $\bx$ is a minimal element of
$\UU_{m-4,A} \cup \{\bx\} $. 
An upper bound for $d_{m-4}(\bx)$ is provided by $d(\bx;\UU_{m-4,A} \cup \bx)$
except when this is zero, so that
\bea
d_{m-4}(\bx)^{8\alpha} & \leq & d(\bx; \UU_{m-4,A}\cup\{\bx\})^{8\alpha} 
%\nonumber \\
%& &
 + d(\bx;\{\bw,\bx,\by,\bz\})^{8\alpha}
I_{m-4}(\bx).
%{\bf 1}\{\bx \textrm{ minimal in }
%\UU_{m-4,A} \cup \{\bx\} \}. 
\label{0817d}
\eea
For $\bx \in B$, it can be shown,
by a similar argument to the one used above for
$L_0$,  that there is a constant $C'$ such that  
\bea
\Exp[(d(\bx;\UU_{m-4,A}\cup\{\bx\}))^{8 \alpha}] <C'.
\label{0817c}
\eea 
Moreover, if $\bw \in A$ with $\bw \preccurlyeq \bx$
and $\|\bx -\bw\| =t >0$, then 
by a similar argument to that at (\ref{0719a}), and (\ref{0719d}), we have
that
$$
%P[ \bx \textrm{ minimal in } \UU_{m-4,A} \cup \{\bx\} ]
\Exp[  I_{m-4}(\bx) ]
\leq \exp(4 - |R(\bw,\bx)|/2) \leq \exp(4 -C t/2), ~~~t \geq 1,
$$  
and hence, uniformly over $A,B$ and $\{\bw,\bx,\by,\bz\}\subset A$ with
$\bx \in B$, we have
$$
\Exp[ d(\bx;\{\bw,\bx,\by,\bz\})^{8\alpha}
%{\bf 1}\{\bx \textrm{ minimal in }
I_{m-4}(\bx)
%\UU_{m-4,A} \cup \{\bx\} \}
] \leq \max\left\{\sup_{t \geq 1}
\left( t^{8 \alpha} \exp(4-Ct/2) \right),1\right\}.
%< \infty. 
$$
Combining this with \eq{0817c}, we see from \eq{0817d}
that $\Exp[d_{m-4}(\bx)^{8 \alpha}]$ is bounded by a constant. 
Also, by a similar argument to (\ref{0719a}) and (\ref{0719d}), 
it can be shown that
 $P[E_{m-4}(\bx)] \leq \exp(4 -C\|\bx\|/2)$ for $\|\bx\| \geq 1$.
 Therefore, by H\"older's inequality,
 the expression (\ref{0728i}) is bounded by a constant times
$$
n^{-4} \int \int \int \int 
\ud \bw 
\ud \bx 
\ud \by 
\ud \bz 
\exp ( - C ( \|\bw\| + 
\| \bx\| +
\| \by\| +
\| \bz\|
)/16 )
$$
and therefore is $O(n^{-4})$. Since the number of
distinct $(i,j,k,\ell)$ in the summation
(\ref{0728h}) is bounded by $m^4$, and hence 
by $(3/2)^4 n^4$, this shows
that the contribution to (\ref{0728h}) from
$i,j,k,\ell$ distinct is uniformly bounded.

Likewise, the number of terms $(i,j,k,\ell)$
with only three distinct values  (e.g., $i=j$ with $i,k,\ell$ distinct)
is $O(n^3)$. Such a term is bounded by an expression like
(\ref{0728i}) but now with a triple integral, which
by a similar argument is $O(n^{-3})$. Hence
the contribution to (\ref{0728h})
of these terms is also bounded.
Similarly, the contribution to (\ref{0728h})
from $(i,j,k,\ell)$
with two distinct values has $O(n^2)$ terms which
are $O(n^{-2})$, and so is bounded.
Likewise  the contribution to (\ref{0728h}) from 
terms with $i=j=k= \ell$ is  bounded. Thus
the expression (\ref{0728h}) is uniformly bounded.

Hence $\Delta (\mathcal{U}_{m,A};B) $ has bounded fourth moments,
uniformly in
$A,B,m$.  $\square$ \\

\noindent
 \textbf{Proof of Theorem \ref{CLT}.}
By Lemmas \ref{stab}, \ref{nondeg}, \ref{ubmlem}
 and the fact that $\LL^\alpha$ is homogeneous of order $\alpha$,
we can apply
Corollary \ref{cltcor}, taking $R_0:= (0,1]^2$ and $S_{0,n} := 
(n^{\eps-\1/2},1]^2$, to obtain
 Theorem \ref{CLT}.  $\square$ \\

\rem An alternative method for proving central limit theorems in
geometrical probability is based on dependency graphs.
%a dependency graph approach
%\cite{baldirinott1989} 
%as 
Such a method was employed by Avram and Bertsimas
\cite{avrambertsimas1993} to give central limit theorems for
nearest neighbour graphs and other random geometrical structures.
A general version of this method is provided by \cite{pynorm}.
By a similar argument to \cite{avrambertsimas1993}, one can show
that, under $\postar$, the total weight (for $\alpha
> 2/3$)
 of edges in the MDST from
points in the region $(\varepsilon_n,1)^2$ (for $\varepsilon_n$
given below) satisfies a central limit theorem, where
\[ \varepsilon_n = \left( \left \lfloor 
\sqrt{ \frac{n}{c \log{n}} } \right \rfloor
\right)^{-1}. \] 
Such an approach can be suitably adapted to show
that a central limit theorem also holds under the more general
partial order specified by $\theta, \phi$, in the region
$(\varepsilon_n, 1-\varepsilon_n)^2$. The benefit of this method
is that it readily yields rates of convergence bounds for the CLT.
The martingale method employed has the advantage of yielding the
convergence of the variance.

\section{The edges near the boundary} \label{bdry}
Next in our analysis of the MDST on random points in the unit square,
we consider  the length of the edges close to the boundary of
the square.  The limiting structure of the MDSF and MDST
 near the boundaries is
described by the directed linear forest model
discussed in Section \ref{secdlt}.
%the preceding  section. 

Initially we consider the `rooted' case where we insert a point at
the origin. Later we analyse the multiple sink (or `unrooted')
case, where we do not insert a point at the origin, in a similar
way.

Fix $\sigma \in (1/2,2/3)$.
Let $B_n$ denote the L-shaped boundary region 
$(0,1]^2 \setminus (n^{-\sigma},1]^2$.
Recall from (\ref{0714a}) that
 $\LL^\alpha (\X ; R)$ denotes the contribution to the
total weight of the MDST on $\X$ from %those
 edges starting at points of $\X \cap R$.
When $\X$ is a random point set, set $\tilde \LL^\alpha (\X; R)
:= 
 \LL^\alpha (\X; R)
- \Exp \LL^\alpha (\X; R)$.

\begin{theorem} \label{thmbdry}
Suppose the partial order is $\postar$.
Then 
%for $\alpha \geq 1$, 
as $n \to \infty$ we have
\bea
 \tLalph ( \Po^0 _n ; B_n ) \tod \tDalph^{\{1\}}
+ \tDalph^{\{2\}} ~~~~~~(\alpha \geq 1);
\label{bdry1} \\
 \tLalph ( \X^0_n ; B_n ) \tod \tDalph^{\{1\}}
+ \tDalph^{\{2\}} ~~~~~~(\alpha \geq 1),
\label{bdry1X}
\eea
 where
$\tDalph^{\{1\}}$, $\tDalph^{\{2\}}$ are
independent random variables with the
distribution of $\tDalph$ given by
the fixed-point equation (\ref{0628a}) for $\alpha=1$
and by (\ref{0628b}) for $\alpha>1$.
Also, as $n \to \infty$,
\bea
 \tLalph  ( \Po_n ; B_n ) \tod
 \tFalph^{\{1\}} + \tFalph^{\{2\}}
~~~~~~(\alpha \geq 1);
\label{bdry2}
\\
 \tLalph  ( \X_n ; B_n ) \tod \tFalph^{\{1\}} +
\tFalph^{\{2\}}
~~~~~~(\alpha \geq 1),
\label{bdry2X}
\eea
 where $\tFalph^{\{1\}}$,
$\tFalph^{\{2\}}$ are independent random variables
with the same distribution as $\tDone$ for $\alpha=1$
and with the distribution given by the fixed-point equation
 (\ref{0628d}) for $\alpha>1$.
Also, as $n\to \infty$,
\bea \label{0214h}
n^{(\alpha-1)/2} \LL^\alpha ( \Po_n ; B_n ) \inL 0 
~~~~~~(0 < \alpha < 1);
%\eea
%and
%\bea
\\
 \label{0802a}
n^{(\alpha-1)/2} \LL^\alpha ( \Po_n^0 ; B_n ) \inL 0 
~~~~~~(0 < \alpha < 1).
\eea
%Further, all the results hold with the Poisson process
%$\Po_n$ replaced by the binomial point process $\X_n$ of $n$
%independent uniformly distributed random points on $(0,1]^2$.
\end{theorem}

The idea behind the proof of Theorem \ref{thmbdry} is
to show that the MDSF near each of the two boundaries is close to
a DLF system defined on a sequence of uniform random variables
coupled to the points of the MDSF. To do this,
%In order to proceed with the proof of Theorem \ref{thmbdry}, we
we produce two explicit sequences of random variables on which we
construct the DLF coupled to
$\Po_n$, a Poisson process of intensity $n$ on $(0,1]^2$,
on which the MDSF is constructed.

\begin{figure}[h]
\centering
\input{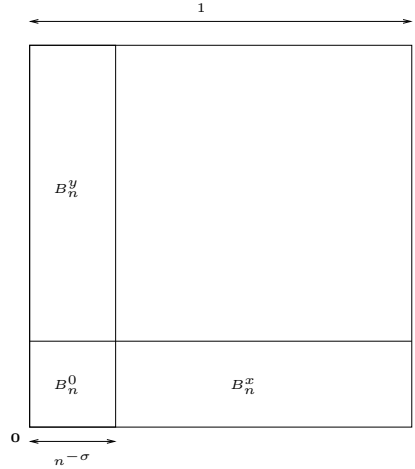}
\caption{The boundary regions}
\label{fig7}
\end{figure}

 Let $B_n^x$ be the rectangle
 $(n^{-\sigma},1] \times (0, n^{-\sigma}]$, let $B_n^y$ be
 the rectangle $(0,n^{-\sigma}]
\times (n^{-\sigma},1]$, and let $B_n^0$ be the square $(0,n^{-\sigma}]^2$;
see Figure \ref{fig7}. Then $B_n = B_n^0 \cup B_n^x \cup B_n^y$.
Define the point processes
\bea
\label{0701f}
 \VV^x_n := \Po_n \cap ( B_n^x \cup B_n^0 ),
~~~
\VV^y_n := \Po_n \cap ( B_n^y \cup B_n^0 ), \textrm{ and } \VV^0_n
:= \Po_n \cap B_n^0.
\eea
 Let $N^x_n := \card( \VV^x_n )$,  $N^y_n :=
\card( \VV^y_n)$ and $N^0_n := \card( \VV^0_n)$. List
$\mathcal{V}^x_n$ in order of increasing $y$-coordinate as
$\bX^x_i$, $i=1,2,\ldots,N^x_n$. In coordinates, set $\bX^x_i =
(X^x_i,Y^x_i)$ for each $i$. Similarly, list $\VV^y_n$ in order of
increasing $x$-coordinate as $\bX^y_i = (X^y_i, Y^y_i)$,
$i=1,\ldots, N^y_n$.  Set $\UU^x_n = ( X^x_i , i=1,2,\ldots,N^x_n
)$ and $\UU^y_n = ( Y^y_i , i=1,2,\ldots,N^y_n )$. Then $\UU^x_n$
and $\UU^y_n$ are sequences of uniform random variables in
$(0,1]$, on which we may construct a DLF. Also, we write $\UU^{x,0}_n$
for the sequence $(0, X^x_1, X^x_2, \ldots, X^x_{N^x_n})$, and
$\UU^{y,0}_n$
for the sequence $(0, Y^y_1, Y^y_2, \ldots, Y^y_{N^y_n})$.

With the total DLF/DLT weight functional $D^\alpha(\cdot)$ 
defined in Section \ref{secdlt} for  random finite sequences in $(0,1)$,
the DLF weight $D^\alpha(\UU^x_n)$ is coupled in a natural way
to the MDSF contribution $\LL^\alpha(\VV_n^x)$, 
and likewise for $D^\alpha(\UU_n^y) $ and
 $\LL^\alpha(\VV_n^y)$,
 for $D^\alpha(\UU_n^{x,0})$ and
 $\LL^\alpha(\VV_n^x \cup \{\0\})$, and
for $D^\alpha(\UU_n^{y,0})$ and   
 $\LL^\alpha(\VV_n^y \cup \{\0\})$.

\begin{lemma} \label{1108c}
For any $\alpha \geq 1$, as $n \to \infty$,
\begin{equation} \label{0209c} \LL ^\alpha (\VV^x_n) - D^\alpha
(\UU^x_n) \inLL 0, \textrm{ and }  \LL^\alpha (\VV^y_n) -
D^\alpha (\UU^y_n) \inLL 0;
\end{equation} 
%and
\begin{equation} \label{0209d}  \LL ^\alpha (\VV^x_n \cup \{ \0 \} )
- D^\alpha
(\UU^{x,0}_n ) \inLL 0, \textrm{ and } 
\LL^\alpha (\VV^y_n \cup \{ \0 \}) -
D^\alpha (\UU^{y,0}_n  ) \inLL 0.
\end{equation}
Further, for $0 < \alpha <1$, as $n\to \infty$,
\bea \label{0214a}
\Exp \left[ \left| \LL^\alpha (\VV^x_n) - D^\alpha (\UU^x_n ) \right|^2
\right] = O \left( n^{2-2\sigma-2\alpha \sigma} \right) ,\eea
and the corresponding result holds for $\VV^y_n$ and $\UU^y_n$, and
for the rooted cases (with the addition of the origin).
\end{lemma} \textbf{Proof.}
We approximate the MDSF in the region $B_n$ by two
DLFs, coupled to the MDSF.
Consider $\VV^x_n$; the argument for $\VV^y_n$ is entirely
analogous.

We have the set of points $\VV^x_n = \{(X^x_i,Y^x_i),
i=1,\ldots,N^x_n\}$. We construct the MDSF on these points,
 and construct the DLF on the $x$-coordinates, $\UU^x_n =
 (X_i^x, i=1,\ldots, N^x_n)$.
Consider any point $(X_i^x,Y_i^x)$. For any single point, either
an edge exists from that point in both constructions, or in
neither. Suppose an edge exists, that is suppose $X_i^x$ is joined
to a point $X^x_{D(i)}$, $D(i)<i$ in the DLF model, and $(X^x_i,Y^x_i)$ to
a point $(X^x_{N(i)},Y^x_{N(i)})$ in the MDST (we do not necessarily have
$N(i)=D(i)$). 
By construction, we know that $| X^x_i-X^x_{D(i)} | \leq
| X^x_i-X^x_{N(i)} |$, since $N(i)<i$ by the order of our points. It then
follows that
 \[ \|(X^x_i,Y^x_i)-(X^x_{N(i)},Y^x_{N(i)})\|^\alpha \geq
|X^x_i-X^x_{N(i)} |^\alpha \geq |X^x_i-X^x_{D(i)} |^\alpha,
\]
%That is, for every point the directed nearest neighbour distance
%must exceed the edge length in the linear scheme,
and so we have established that, for all $\alpha>0$,
\[
%\label{1122b}
D^\alpha ( \UU^x_n ) \leq \LL^\alpha
( \VV^x_n ) ; \textrm{ and } D^\alpha ( \UU^{x,0}_n  )
\leq \LL^\alpha (\VV^x_n \cup \{ \0 \}) . \]
Now, by the construction of the MDST, we have that
\bea \label{0209a}
\| (X^x_i,Y^x_i)-(X^x_{N(i)},Y^x_{N(i)})
\| \leq \| (X^x_i,Y^x_i)-(X^x_{D(i)},Y^x_{D(i)}) \|.
\eea
If $(x,y) \in (0,1]^2$ then $\| (x,y) \| \leq x + y$, and by the
Mean Value Theorem for the function $t \mapsto t^\alpha$, for $\alpha
\geq 1$,
\[ \| (x,y) \|^\alpha - x^\alpha \leq (x+y)^\alpha - 
x^\alpha \leq \alpha 2^{\alpha-1} y ~~~ (\alpha \geq 1). \]
Hence, for $\alpha \geq 1$,
\bea \label{0209b}
 \| (X_i^x, Y_i^x) - (X_{D(i)}^x, Y_{D(i)}^x) 
\|^\alpha - ( X_i^x - X_{D(i)}^x)
^\alpha \leq \alpha 2^{\alpha -1}
(Y_i^x - Y_{D(i)}^x) .
\eea
Then (\ref{0209a}) and (\ref{0209b}) yield, for $\alpha \geq 1$,
\[ \| ( X^x_i, Y^x_i ) - (X^x_{N(i)}, Y^x_{N(i)} ) \|^\alpha
- (X_i^x - X^x_{D(i)} )^\alpha
\leq \alpha 2^{\alpha-1} ( Y^x_i - Y^x_{D(i)} ) . \]
Hence, for $\alpha \geq 1$,
\[ 0 \leq \LL^\alpha ( \VV_n^x ) -
 D^\alpha ( \UU_n ^x ) \leq \alpha 2^{\alpha-1}
\sum_{i=1}^{N_n^x} (Y_i^x - Y_{D(i)}^x). \]
Thus, for $\alpha \geq 1$,
\begin{eqnarray}
\label{1202e} & & 0  \leq \LL^\alpha ( \VV^x_n ) - D^\alpha ( \UU^x_n )
\leq \alpha 2^{\alpha-1} N^x_n n^{-\sigma} ; \nonumber\\
 & \textrm{and } &
0 \leq
\LL^\alpha ( \VV^x_n \cup \{ \0 \} ) - D^\alpha ( \UU^{x,0}_n  )
\leq \alpha 2^{\alpha-1} N^x_n n^{-\sigma} .
\end{eqnarray}
We have $N^x_n \sim \textrm{Po}\left(n^{1-\sigma}\right)$, so that
since $\sigma > 1/2$, we have
$$
\Exp [
(\LL^\alpha ( \VV^x_n \cup \{ \0 \} ) - D^\alpha ( \UU^{x,0}_n  ) )^2 ]
\leq \alpha^2 2^{2 \alpha -2}  n^{-2 \sigma} \Exp[(N_n^x)^2] \to 0,
~~~ \alpha \geq 1.
$$
 An entirely analogous
argument leads to the same statement for $\UU^y_n$ and $\VV^y_n$, and we
obtain (\ref{0209c}), and (\ref{0209d}) in identical fashion.

We now consider $0 < \alpha <1$.
By the concavity of the function $t \mapsto t^\alpha$
for $\alpha <1$,
we have for $x >0, y>0$ that
 \[
 \| (x,y) \|^\alpha -x^\alpha \leq
(x+y)^\alpha -x^\alpha \leq y^\alpha ~~~ (0 <\alpha <1) .\]
Then, by a similar argument to (\ref{1202e}) in
the $\alpha \geq 1$ case, we obtain
\[ 0 \leq \LL^\alpha (\VV^x_n ) - D^\alpha (\UU^x_n ) \leq N_n^x n^{-\alpha
\sigma} . \]
Then (\ref{0214a}) follows since $N_n^x \sim 
\textrm{Po}\left(n^{1-\sigma}\right)$, and the rooted case is similar.
$\square$

\begin{lemma} \label{1108d} 
%For $\alpha > 1$,
Suppose $\tDone$ has distribution given by (\ref{0628a}),
$\tDalph$, $\alpha >1$, has distribution given by (\ref{0628b}), and 
$\tFalph$, $\alpha >1$, has distribution given by (\ref{0628d}).
Then
as $n \to \infty$,
\bea
\tilde \LL^1 (\VV^x_n \cup \{ \0 \} ) \tod \tDone,
 %\textrm{ and } 
{\rm ~~~and ~~~}
 \tilde \LL^1 (\VV^x_n) \tod \tDone;
\label{0715a}
\eea
\bea
\tilde
\LL^\alpha (\VV^x_n \cup \{ \0 \} ) \tod \tDalph,
{\rm ~~~and ~~~}
% \textrm{ and }
 \tilde \LL^\alpha (\VV^x_n) \tod \tFalph ~~~(\alpha >1).
\label{0715b}
\eea
Moreover, (\ref{0715a}) and (\ref{0715b}) also hold
with $\VV^x_n$ replaced by $\VV^y_n$.
%where $\tF^1$ has distribution given by (\ref{0628a}) and
\end{lemma}
\proof
As usual we present the argument for $\VV_n^x$ only, since
the result for $\VV_n^y$ follows in the same manner.
First consider the $\alpha>1$ case. We have the distributional
equality
\[
 \LL \left( \left. D^\alpha ( \UU_n^{x,0} ) \right| N_n^x =m
\right) = \LL \left( D^\alpha (\UU_m^0) \right) ; ~~~~~
 \LL \left( \left. D^\alpha ( \UU_n^{x} ) \right| N_n^x =m
\right) = \LL \left( D^\alpha (\UU_m) \right) .
\]
But $N_n^x$ is Poisson with mean $n^{1-\sigma}$, and so
tends to infinity almost surely. Thus by
% Proposition \ref{1119f},
Theorem \ref{dltthm} (ii),
$D^\alpha (\UU_n^{x,0}) \tod \Dalph$
and
$D^\alpha (\UU_n^{x}) \tod \Falph$
 as $n \to \infty$, and so
by Lemma \ref{1108c} and Slutsky's theorem,
 we obtain 
\bea
 \LL^\alpha (\VV_n^x \cup \{ \0 \}) \tod \Dalph 
{\rm ~~~ and ~~~}
 \LL^\alpha (\VV_n^x ) \tod \Falph 
{\rm ~~ as ~~} n \to \infty.
\label{0803} 
\eea
Also, $\Exp [D^\alpha(\UU_n^{x,0})] \to (\alpha -1)^{-1} $
by (\ref{0720c}), so by Lemma \ref{1108c} and Proposition
\ref{1119f}, 
$\Exp[\LL^\alpha
(\VV_n^x \cup \{\0\}) ] \to (\alpha -1)^{-1} =
\Exp[\Dalph]$. 
Similarly, by (\ref{0720f}), 
 Lemma \ref{1108c} and
 Proposition \ref{1125a},
$\Exp[\LL^\alpha
(\VV_n^x ) ] \to (\alpha(\alpha -1))^{-1} =
\Exp[\Falph]$. 
Hence, (\ref{0803}) still holds with the centred variables,
i.e., \eq{0715b} holds.

Now suppose $\alpha=1$.
Since $N_n^x$ is Poisson with parameter $n^{1-\sigma}$,
Lemma \ref{0630c} (i), with $t=n^{1-\sigma}$,
then shows that
$
\tD^1 (\UU^{x,0}_n) \tod \tDone$
as $n \to \infty$.
Slutsky's theorem with Lemma \ref{1108c}
then implies that
$\tLL^1 (\VV^x_n \cup \{ \0 \}) \tod \tDone$.
In the same way we obtain
$\tLL^1 (\VV^x_n ) \tod \tDone$,
this time using part (ii) instead of part (i) of
 Lemma \ref{0630c}, along with Proposition \ref{ffixed}.
$\square$ \\

Note that
$D^\alpha (\UU^x_n)$ and $D^\alpha (\UU^y_n)$ are not independent.
 To
deal with this, we
define
\[ \tilde \VV^x_n := \Po_n \cap B_x^n, \textrm{ and }
\tilde \VV^y_n := \Po_n \cap B_y^n . \]
Also, recall the definition of $\VV^0_n$ at (\ref{0701f}).
Let $\tilde N^x_n :=
\card( \tilde \VV^x_n )$ and $\tilde N^y_n := \card( \tilde
\VV^y_n)$. 
Since $B_n^x$ and $B_n^y$ are
disjoint, $\LL^\alpha (\tilde \VV^x_n)$ and $\LL^\alpha (\tilde
\VV^y_n)$ are independent,
% and $D^\alpha (\tilde \UU^x_n)$ and 
%$D^\alpha (\tilde \UU^y_n)$ are also independent, 
by the spatial
independence property of the Poisson process $\Po_n$.

\begin{lemma} \label{1108e} Suppose $\alpha > 0$. Then:

(i)  As $n \to \infty$,
\bea
 \LL^\alpha(\VV^x_n) - \LL^\alpha (\tilde \VV^x_n) 
\inL 0 , \textrm{ and } \LL^\alpha (\VV^y_n) - \LL^\alpha (\tilde
\VV^y_n) \inL 0; 
\label{0803a}
\eea
\bea
 \LL^\alpha(\VV^x_n\cup\{\0\}) - \LL^\alpha (\tilde
\VV^x_n\cup\{\0\}) \inL 0 , \textrm{ and } \LL^\alpha (\VV^y_n\cup\{\0\}) - \LL^\alpha (\tilde
\VV^y_n\cup\{\0\}) \inL 0. ~ 
\label{0803b}
\eea
% \item[(ii)]

(ii)  As $n \to \infty$, we have
$ \LL^\alpha (\VV^0_n) \inL 0 $,
 and
$ \LL^\alpha (\VV^0_n\cup \{\0\}) \inL 0  $.
%\end{itemize}
\end{lemma}
\textbf{Proof.}
We first prove (i). We give only the argument
for $\VV_n^x$; that for $\VV_n^y$ is analogous.
Set $\Delta := 
 \LL^\alpha(\VV^x_n) - \LL^\alpha (\tilde \VV^x_n)$. 
Let $\beta = (\sigma +(1/2))/2$.
Then $1/2 < \beta < \sigma$.

 Assume without loss of generality that $\Po_n$ is the restriction
 to $(0,1]^2$ of a homogeneous Poisson process $\H_n$ of intensity $n$ on 
$\R^2$.  Let $\bX^-=(X^-,Y^-)$
be the point of $\H_n \cap ( (0,n^{-\beta}] \times (0,\infty))$
with minimal $y$-coordinate. 
Then $X^-$ is uniform on $(0,n^{-\beta}]$.
Let $E_n$ be the event
that $X^- > 3 n^{-\sigma}$; then  
$\Pr[E_n^c] =3 n^{\beta - \sigma} $ for $n$ large enough.

Let $\Delta_1$ be the the contribution to $\Delta$ from edges
starting at points in $(0,n^{-\beta}]\times(0,n^{-\sigma}]$.
Then the absolute value of 
$\Delta_1$ is bounded by the product of $(\sqrt{2}n^{-\beta})^\alpha$
 and the number of points of $\Po_n $ in $ (0,n^{-\beta}] \times
(0,n^{-\sigma}]$. Hence,  for any $\alpha >0$,
\bea
 \Exp \left[  | \Delta_1 | \right]
& \leq & (\sqrt{2}n^{-\beta})^\alpha \Exp
\left[ \card \left(
\Po_n \cap ((0,n^{-\beta}] \times (0,n^{-\sigma}] ) \right)
\right]
\nonumber \\
& = & 2^{\alpha/2} n^{1-\beta-\sigma-\alpha \beta} \to 0. 
\label{0802b}
\eea
 %since $\sigma>\beta >1/2$.

Let $\Delta_2 := \Delta - \Delta_1$, the contribution
to $\Delta$ from edges starting at points in
$(n^{-\beta},1]\times(0,n^{-\sigma}]$.
Then by the triangle inequality,
if $E_n$ occurs then these edges are unaffected by points in
$B_n^0$, so that
 $\Delta_2$ is zero if $E_n$ occurs.
Also, 
only minimal elements
of $\Po_n \cap  (n^{-\beta},1]\times(0,n^{-\sigma}]$ can
possibly have their
 directed nearest neighbour in 
$(0,n^{-\sigma}] \times (0,n^{-\sigma}]$; hence,  
if $M_n$ denotes 
the number of such minimal
elements then
 $|\Delta_2|$ is bounded  by  $2^{\alpha/2}M_n$.
Hence, using (\ref{harmonicbd}), we obtain
$$
\Exp[ | \Delta_2 | ] \leq 2^{\alpha/2}P[E_n^c] \Exp[M_n] = O(n^{\beta -\sigma}
\log n )  
$$
which tends to zero. Combined with (\ref{0802b}), this  gives
us (\ref{0803a}). The same argument gives us (\ref{0803b}). 
%part (i). 

For (ii), note that
\[ \Exp \left[
 \LL^\alpha (\VV^0_n) \right] \leq (\sqrt{2}n^{-\sigma})^\alpha
\Exp[ N_n^0] = 2^{\alpha/2} n^{1-2\sigma-\sigma\alpha}
 \to 0 , \textrm{ as } n\to
\infty, \]
for any $\alpha >0$. Thus
 $\LL^\alpha (\VV^0_n) \inL 0$,
and similarly 
 $\LL^\alpha (\VV^0_n \cup \{\0\}) \inL 0$.
$\square$ \\

In proving our next lemma (and again later on)
we use the following elementary fact.
If $N(n)$ is Poisson with parameter $n$, then
as $n \to \infty$ we have
\bea
\Exp [|N(n)-n|\log \max(N(n),n)] = O(n^{1/2} \log n).
\label{logmax}
\eea
To see this, set $Y_n :=
|N(n)-n|\log \max(N(n),n)$.
Then $Y_n\1_{\{N(n) \leq 2n\} } \leq |N(n)-n |\log (2 n) $,
and the expectation of this is $O(n^{1/2} \log n)$ by Jensen's
inequality since $\Var(N(n)) = n$. On the other hand,
the Cauchy-Schwarz inequality shows that
 $\Exp[Y_n\1_{\{N(n) > 2n\} } ] \to 0$, and (\ref{logmax}) follows.

We now state a lemma for coupling $\X_n $ and $\Po_n$.
 The $\alpha \geq 1$ part will be
used in the proof of Theorem \ref{thmbdry}. The $0<\alpha<1$
part will be needed later, in the proof of Theorem \ref{mainth}.
As in Section \ref{ltot}, let $S_{0,n}$ denote
the `inner' region $(n^{\eps-1/2},1]^2$, with $\eps \in (0,1/2)$ a constant.
The boundary region $B_n$ is disjoint from $ S_{0,n} $;
let $C_n$ denote the intermediate region $(0,1]^2 \setminus
( B_n \cup S_{0,n} )$, so that $B_n \cup C_n = (0,1]^2 \setminus S_{0,n}$. 

\begin{lemma}
\label{lem0803} There exists a coupling of $\X_n$
and $\Po_n$ such that:  
\begin{itemize}
\item[(i)]
For $0< \alpha <1$, provided $\eps < (1-\alpha)/2$, 
we have that as $n \to \infty$,
\bea
\label{0804f}
n^{(\alpha -1)/2} \Exp [|\LL^\alpha(\X_n;B_n \cup C_n) - \LL^\alpha(\Po_n;
B_n \cup C_n) |]   \to 0 
\eea
and
\bea
\label{0804g}
n^{(\alpha -1)/2} \Exp [|\LL^\alpha(\X_n^0;B_n \cup C_n) - \LL^\alpha(\Po^0_n;
B_n \cup C_n) |]   \to 0. 
\eea
\item[(ii)]
For $\alpha \geq 1$, we have that
as $n \to \infty$,
\bea
\label{0803d}
\Exp [|\LL^\alpha(\X_n;B_n ) - \LL^\alpha(\Po_n;
B_n ) |]   \to 0 
\eea
and
\bea
\label{0803e}
\Exp [|\LL^\alpha(\X_n^0;B_n) - \LL^\alpha(\Po^0_n;
B_n ) |]   \to 0. 
\eea
\end{itemize}
\end{lemma}
\proof 
We couple $\X_n$ and $\Po_n$ in the following standard way.
 Let $\bX_1, \bX_2, \bX_3, \ldots$
be independent uniform random vectors on $(0,1]^2$, and
 let $N(n) \sim {\rm Po}(n)$
be independent of $(\bX_1,\bX_2,\ldots)$. For $m \in \N$
(and in particular for $m =n$) set
$\X_m := \{ \bX_1, \ldots, \bX_m \}$;
set $\Po_n := \{ \bX_1, \ldots, \bX_{N(n)} \}$.

For each $m \in \N$, let $Y_m$ denote the
in-degree of vertex $\bX_m$ in the MDST on $\X_{m}$. 
Suppose $\bX_m=\bx$. Then an upper bound for
$Y_m$ is provided by the number of minimal
elements of the restriction of $\X_{m-1}$
to the rectangle $\{\by \in (0,1]^2: \bx \postar \by\}$.
Hence, conditional on $\bX_m =\bx$ and on 
there being $k$ points of $\X_{m-1}$ in this rectangle,
the expected value of $Y_m$ is bounded
by the expected number of minimal elements
in a random uniform sample of $k$ points in this
rectangle, and hence
(see (\ref{harmonicbd})) by
$1 + \log k$.
Hence, given the value of $\bX_m$, the 
conditional expectation of $Y_m$ is bounded by
 $1+ \log m$. 

First we prove the statements in part (i) ($0<\alpha <1$). Suppose
$\eps < (1-\alpha)/2$. Then
\bea
|\LL^\alpha(\X_m;B_n\cup C_n) - \LL^\alpha(\X_{m-1};B_n \cup C_n)
 |
\leq
2^{\alpha/2}( Y_m +1)
 {\bf 1}\{\bX_m \in B_n \cup C_n\}. ~
\label{0803c1}
\eea
 Since $B_n\cup C_n$ has area 
$2 n^{\eps -1/2} - n^{2\eps -1} $,
we obtain
$$
\Exp[ (Y_m +1) {\bf 1}\{\bX_m \in B_n\cup C_n\} ] \leq
 (2 + \log m )2 n^{\eps -1/2} .
$$
Hence, by (\ref{0803c1}) there is a constant $C$ such that 
\bean
n^{(\alpha -1)/2} \Exp[ (| \LL^\alpha(\Po_n;B_n \cup C_n) -
\LL^\alpha(\X_n;B_n \cup C_n) | )| N(n) ] 
\\
\leq C 
|N(n)-n|
\log(\max(N(n),n) ) n^{(\alpha + 2 \eps -2)/2 },  
\eean
and since we assume $\alpha + 2\eps  < 1 $,
by \eq{logmax}
the expected value of the right hand side tends to zero
as $n \to \infty$, and we obtain (\ref{0804f}). Likewise in the rooted case
(\ref{0804g}).

Now we prove part (ii). For $\alpha \geq 1$, we have 
\bea
|\LL^\alpha(\X_m;B_n) - \LL^\alpha(\X_{m-1};B_n)
 |
\leq
2^{\alpha/2}( Y_m +1)
 {\bf 1}\{\bX_m \in B_n \}.
\label{0803c2}
\eea
Since $B_n$ has area 
$2 n^{-\sigma} - n^{-2\sigma}$,
by (\ref{0803c2}) there is a constant $C$ such that 
\bean
\Exp[ (| \LL^\alpha(\Po_n;B_n) - \LL^\alpha(\X_n;B_n) | )| N(n) ] 
%\\
\leq C |N(n)-n| \log(\max(N(n),n) ) n^{-\sigma },  
\eean
and since $\sigma> 1/2$, by \eq{logmax}
the expected value of the right hand side tends to zero
as $n \to \infty$, and we obtain (\ref{0803d}). 
We get (\ref{0803e}) similarly.
$\square$ \\

\noindent \textbf{Proof of Theorem \ref{thmbdry}.} Suppose
$\alpha \geq 1$. We have that
\bean
\tilde \LL^\alpha ( \tilde \VV_n^x ) = \tilde \LL^\alpha (
\VV_n^x) + ( \tilde \LL^\alpha ( \tilde \VV_n^x ) - \tilde \LL^\alpha (
\VV_n^x) ).\eean
The final bracket converges to zero in probability,
by Lemma \ref{1108e} (i). Thus by Lemma \ref{1108d} and Slutsky's theorem,
we obtain $\tilde \LL^\alpha ( \tilde \VV_n^x) \tod \tFalph$
(where we have $\tFone \eqd \tDone$).
Now
\bean
\tilde \LL^\alpha ( \VV_n^x) +
\tilde \LL^\alpha ( \VV_n^y) =
\tilde \LL^\alpha ( \tilde \VV_n^x)
+
\tilde \LL^\alpha ( \tilde \VV_n^y)
+ (\tilde \LL^\alpha ( \VV_n^x)-
\tilde \LL^\alpha ( \tilde \VV_n^x))
+(\tilde \LL^\alpha ( \VV_n^y)-
\tilde \LL^\alpha ( \tilde \VV_n^y)).
\eean
The last two brackets converge to zero in probability,
by Lemma \ref{1108e} (i). Then the independence
of $\tilde \LL^\alpha (\VV^x_n)$ and $\tilde \LL^\alpha (\VV^y_n)$
and another
application of Slutsky's theorem 
%and Lemma \ref{1108d}
 yield
\[ \tilde \LL^\alpha (\VV^x_n) + \tilde \LL^\alpha (\VV^y_n)
\tod \tFalph^{\{1\}} + \tFalph^{\{2\}} ,
\]
where $\tFalph^{\{1\}}$ and $\tFalph^{\{2\}}$
are independent copies of $\tFalph$.
Similarly,
\[ \tilde \LL^\alpha (\VV^x_n \cup \{ \0 \} ) +
\tilde \LL^\alpha (\VV^y_n \cup \{ \0 \} ) \tod \tDalph^{\{1\}}
+ \tDalph^{\{2\}} . \]
Finally, since $\tilde \LL^\alpha (\Po_n ; B_n)
= \tilde \LL^\alpha (\VV^x_n)+\tilde \LL^\alpha (\VV^y_n)
-\tilde \LL^\alpha (\VV^0_n)$
(with a similar statement including the origin)
Lemma \ref{1108e} (ii) and Slutsky's
theorem complete the proof of (\ref{bdry1}) and (\ref{bdry2}).
%in the Poisson case for $\alpha \geq 1$.

To deduce (\ref{bdry1X}) and (\ref{bdry2X}), 
assume without loss of generality
that $\X_n$ and $\Po_n$ are coupled
  in the manner of Lemma \ref{lem0803}. 
 Then $\tLL^\alpha(\Po_n;B_n) -  \tLL^\alpha(\X_n;B_n) $ tends
to zero in probability by (\ref{0803d}), and 
 $\tLL^\alpha(\Po_n^0;B_n) -  \tLL^\alpha(\X_n^0;B_n) $ tends
to zero in probability by (\ref{0803e}). Hence by Slutsky's theorem,
the convergence  results (\ref{bdry1}) and (\ref{bdry2})
 carry through to the binomial point process case,
i.e., (\ref{bdry1X}) and (\ref{bdry2X}) hold.

Now suppose
 $0 < \alpha <1$.
Then (\ref{0214a}) gives us
\bea \Exp \left[ \left| n^{(\alpha-1)/2} \left( \LL^\alpha (\VV^x_n ) -
D^\alpha (\UU^x_n) \right) \right|^2 \right] = O \left(
n^{(\alpha +1)(1 -2\sigma)} \right), 
\label{0802c}
\eea
which tends to 0 as $n\to \infty$,
  since $\sigma>1/2$. Likewise for the rooted case,
\bea
 \Exp \left[ \left| n^{(\alpha-1)/2} \left( \LL^\alpha (\VV^x_n \cup\{\0\}
 ) -
D^\alpha (\UU^{x,0}_n) \right) \right|^2 \right] = O \left(
n^{(\alpha +1)(1 -2\sigma)} \right), 
\label{0802d}
\eea
%Also, for $0<\alpha<1$,
%from
By Proposition \ref{dlfmoms} 
we have
 \[
\Exp[ n^{(\alpha-1)/2} D^\alpha (\UU^x_n) ] = O(n^{(\alpha-1)/2}
\Exp[(N_n^x)^{1-\alpha} ] ) = O(n^{(\alpha -1)(\sigma -1/2)}) \to 0,
%\inL 0 ,
 \]
%and the fact that
%$\Exp[(N^x_n)^{1-\alpha}] = O( n^{(1-\alpha)(1-\sigma)})$,
and combined with (\ref{0802c}) this 
 completes the proof of (\ref{0214h}). Similarly, 
by Proposition \ref{dltmoms}, 
 \[
\Exp[ n^{(\alpha-1)/2} D^\alpha (\UU^{x,0}_n) ] = O(n^{(\alpha-1)/2}
\Exp[(N_n^x)^{1-\alpha} ] ) = O(n^{(\alpha -1)(\sigma -1/2)}) \to 0,
 \]
and  combined with (\ref{0802d}) this gives us (\ref{0802a}).  $\square$

\section{Proof of Theorem \ref{mainth}} \label{totallength}

Let  $ \sigma \in (1/2, 2/3)$.
Let $\eps > 0 $ with
\bea
\eps < \min (1/2, (1-\sigma)/3,
(3 - 4 \sigma) /10, (2 - 3 \sigma)/8).
\label{epsdef}
\eea
In addition,  if $0<\alpha <1$,
we impose the further condition
that $\eps < (1-\alpha)/2$. As in Section \ref{ltot}, 
 denote by $S_{0,n}$ the region
 $(n^{\eps-1/2},1]^2$.
% (which we denoted $S_{0,n}$
As in Section \ref{bdry}, let $B_n$
denote the region $(0,1]^2 \setminus (n^{-\sigma},1]^2$,
and let $C_n$ denote $(0,1]^2 \setminus
( B_n \cup S_{0,n} )$. 

We know from Sections \ref{ltot} and \ref{bdry}
 that, for large $n$, the weight of
edges starting in $S_{0,n}$ satisfies
 a central limit theorem, and the weight
of edges starting in
 $B_n$ can be approximated by the directed linear forest.
We shall show in Lemmas \ref{varClem} and \ref{lem0817}
 that (with a suitable
scaling factor for $\alpha<1$) the contribution
to the total weight from points in $C_n$ 
has variance converging to 
 zero. 
To complete the proof of Theorem \ref{mainth} in the Poisson case,  
we shall show that the lengths from $B_n$ and $S_{0,n}$
are asymptotically independent by virtue of the fact that
the configuration of points in $C_n$ is (with probability
approaching one) sufficient to ensure that the configuration
of points in $B_n$ has no effect on the edges from points in $S_{0,n}$.
To extend the result to the binomial point process case,
we shall use a de-Poissonization argument related to that used 
in \cite{penyuk1}.

\begin{figure}[h]
\centering
\input{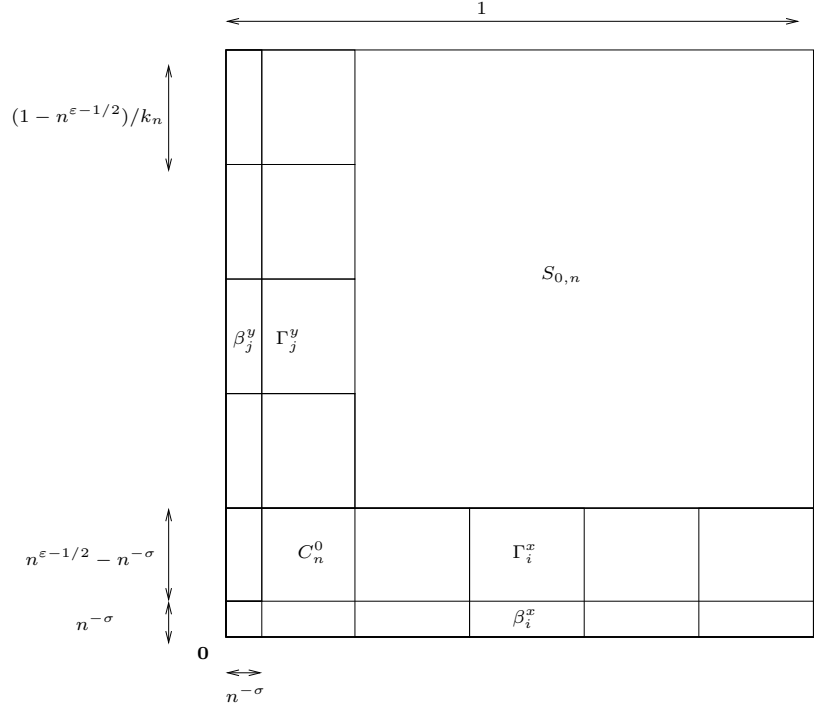}
\caption{The regions of $[0,1]^2$.} 
\label{fig1}
\end{figure}

First consider the region $C_n$. We naturally divide this
into three regions. Let
\bean
C_n^x := (n^{\eps -1/2},1] \times
(n^{-\sigma},n^{\eps-1/2}], ~~
C_n^y := (n^{-\sigma},n^{\eps-1/2}] \times (n^{\eps-1/2},1], 
\\
C_n^0 := (n^{-\sigma},n^{\eps-1/2}]^2.
\eean
Also, as in Section \ref{bdry}, let
\bean
B_n^x := (n^{-\sigma},1] \times (0,n^{-\sigma}], ~~
B_n^y := (0,n^{-\sigma}] \times (n^{-\sigma},1], ~~
B_n^0 := (0,n^{-\sigma}]^2.
\eean
We 
 divide the $C_n$ and $B_n$ into rectangular cells as follows
(see Figure \ref{fig1}.)
We leave $C_n^0$ undivided. We set
\bea
 k_n := \lfloor n^{1-\sigma -2 \eps} \rfloor
\label{kndef}
\eea
and
divide $C_n^x$ lengthways into $k_n$ cells.
For each cell, 
\bea
{\rm width} = (1 - n^{\eps- 1/2})/k_n \sim n^{2 \eps + \sigma -1 }
; ~~~~~
{\rm height} = n^{\eps - 1/2} - n^{-\sigma} \sim n^{\eps -1/2} 
.
\label{0729}
\eea
Label these
cells $\Gamma_i^x$ for $i=1,2,\ldots,k_n$ from left to right.
For each cell $\Gamma^x_i$, define the adjoining
cell of $B_n^x$, formed by extending the
vertical edges of $\Gamma^x_i$, to be $\beta^x_i$.
The cells $\beta^x_i$ then have width
$(1-n^{\eps-1/2})/k_n \sim n^{2 \eps + \sigma -1}$
and height $n^{-\sigma}$.

In a
similar way we divide $C_n^y$ into $k_n$ cells $\Gamma_i^y$ of
 height $(1-n^{\eps-1/2})/k_n$
 and
width $n^{\eps-1/2}-n^{-\sigma}$,
 and
divide $B_n^y$ into the corresponding cells $\beta^y_i$, $i=1,\ldots,k_n$.
%\Comment{see caption}

For $i=2,\ldots,k_n$, let $E_{x,i}$ denote the event that the cell
$\beta^x_{i-1}$
contains at least one point of $\Po_n$, and
 let $E_{y,i}$ denote the event that
$\beta^y_{i-1}$ contains at least one point of $\Po_n$.

\begin{lemma}
\label{neighbours} 
 For $n$ sufficiently large,
and for $1 \leq j <i \leq k_n$ with $i -j >3$,
 if $E_{x,i}$ (respectively $E_{y,i}$) occurs then
no point in the cell $\Gamma^x_i$ 
(respectively  $\Gamma^y_i$)
has a directed nearest neighbour in the cell $\Gamma^x_j$
or $\beta_j^x$
($\Gamma^y_j$ or $\beta_j^y$).
\end{lemma} 
\textbf{Proof.}
Consider a point $X$, say, in cell $\Gamma_i^x$ in $C_n^x$. 
Given $E_{x,i}$, we know that there is a point, $Y$ say, in the
cell $\beta^x_{i-1}$ to the left of the $\beta^x_i$ cell
immediately below $\Gamma_i^x$,
 such that $Y \postar X$, but the difference
in $x$-coordinates between $X$ and $Y$ is no more than twice the
width of a cell. So, by the triangle inequality, we have 
\bea
\label{0712a} \| X - Y \| \leq 2 (1-n^{\eps-1/2})/k_n +
n^{\eps-1/2} \sim 2 n^{2 \eps + \sigma -1  }, 
\eea 
since $\sigma > 1/2$.
Now, consider a point $Z$ in a cell
$\Gamma^x_j$ or $\beta_j^x$ 
with $j \leq i-4$. In this case, the difference in
$x$-coordinates between $X$ and $Z$ is at least the width of 3
cells, so that 
\bea
 \label{0712b}
 \| X-Z \| \geq 3 (1-n^{\eps-1/2})/k_n \sim 3 n^{2 \eps + \sigma -1 }. 
\eea
Comparing (\ref{0712a}) and (\ref{0712b}), we see that
$X$ is not connected to $Z$, which completes the proof.  
$\square$ \\

Recall from (\ref{0714a})
% Section \ref{con}
 that for a point set $\SS \subset \R^2$
and a region $R \subseteq \R^2$,
 $\LL^\alpha(\SS;R)$ denotes
the total weight of edges of the MDSF on $\SS$ which originate
in the region $R$.

\begin{lemma}
\label{varClem}
As $n \to \infty$, we have that
\bea
\label{0525aa}
 \Var [ \LL^\alpha( \Po_n ; C_n) ] \to 0
~~{\rm and} ~~
 \Var [ \LL^\alpha( \Po_n^0 ; C_n) ] \to 0
~~~~~ (\alpha \geq 1);
\eea
\bea
\label{0525bb}
 \Var [ n^{(\alpha-1)/2}
\LL^\alpha( \Po_n ; C_n) ] \to 0 
~~~~~ (0 < \alpha < 1);
\eea
\bea
\label{0525bbb}
 \Var [ n^{(\alpha-1)/2}
\LL^\alpha( \Po_n^0 ; C_n) ] \to 0 
~~~~~ (0 < \alpha < 1).
\eea
%[MP: do we need $\X_n$ version of these last two???] \Comment{MP}
\end{lemma}
\proof For ease of notation, write $X_i = \LL^\alpha (\Po_n ;
\Gamma^x_i)$ and $Y_i = \LL^\alpha (\Po_n ; \Gamma^y_i)$, for
$i=1,2,\ldots,k_n$. Also let $Z = \LL^\alpha (\Po_n ; C_n^0 )$.
Then 
\bea 
\label{0712v}
 \Var[ \LL^\alpha( \Po_n ; C_n)
] & = & \Var \left[ Z + \sum_{i=1}^{k_n} X_i + \sum_{i=1}^{k_n}
Y_i 
 \right].
\eea
Let $N^x_i$, $N^y_i$, $N_0$, respectively,
 denote the number of points of $\Po_n$ in $\Gamma_i^x$, $\Gamma_i^y$,
$C_n^0$, respectively.
Then by (\ref{0729}), 
$N^x_i$ 
is Poisson with parameter asymptotic
to $n^{3 \eps +\sigma-1/2}$,
while $N^x_1 + N^y_1 + N_0$ 
is Poisson with parameter asymptotic
to $2 n^{3 \eps +\sigma-1/2}$; hence as $n \to \infty$
and we have
 \bea
 \label{0712g}
 \Exp[(N_i^x)^2] \sim n^{6\eps + 2\sigma-1}, ~~~~
 \Exp[(N_1^x + N_1^y + N_0)^2] \sim 4 n^{6\eps + 2\sigma-1} .
  \eea 
Edges from points in $\Gamma_1^x \cap \Gamma_1^y \cap C_n^0$
are of length at most $2 n^{2\eps +\sigma -1}$, and hence,
 \bea
 \Var [ X_1 + Y_1 + Z] & \leq  & (2n^{2 \eps + \sigma -1} )^{2 \alpha}
 \Exp [ (N_1^x + N_1^y + N_0)^2]
\nonumber \\
  & \sim & 2^{2 + 2 \alpha}
n^{6 \eps + 2 \sigma -1 + 2 \alpha (2 \eps + \sigma -1) }.
%\to 0. 
\label{0802e}
 \eea
For $\alpha \geq 1$, since $\eps$ is small (\ref{epsdef}), 
 the expression (\ref{0802e}) is $O(n^{10 \eps + 4 \sigma -3})$
and in fact  tends to zero, so that 
\bea
\Var(X_1 + Y_1 + Z) \to 0 ~~~(\alpha \geq 1).
\label{0712x}
\eea 

By Lemma \ref{neighbours} and (\ref{0712a}), given
$E_{x,i}$, an edge from a point of $\Gamma^x_i$ can be of length
no more than 
%$2/k_n + n^{\eps -1/2}$, which is bounded by
 $3 n^{2 \eps + \sigma -1}$.
 Thus using (\ref{0712g}) we have 
\bea
\label{0712d} \Var [ X_i \1 \{E_{x,i} \} ]  & \leq &  \Exp [ X_i^2 
{\bf 1}\{E_{x,i}\}] \leq 
(3 n^{2\eps + \sigma -1})^{2\alpha} 
\Exp [ (N^x_i)^2 ]
\nonumber\\ & = &
O(n^{6 \eps + 2 \sigma-1  + 2 \alpha(2 \eps + \sigma-1) }) .
 \eea
Next, observe that ${\rm Cov} [ X_i \1 \{ E_{x,i} \} , X_j \1 \{ E_{x,j}
\}] =0$ for $i-j >3$, since  by Lemma \ref{neighbours},
$X_i {\bf 1} \{E_{x,i}\}$ is determined by the restriction of
$\Po_n$ to the union of the regions $\Gamma_\ell^x  \cup \beta_\ell^x,
i-3 \leq \ell \leq i$.
 Thus 
by (\ref{kndef}), Cauchy-Schwarz and (\ref{0712d}),
we obtain
 \bea 
\label{0712e}
 \Var \left[ \sum_{i=2}^{k_n} X_i \1
\{ E_{x,i} \} \right] 
& = & \sum_{i=2}^{k_n} \Var
[X_i \1 \{ E_{x,i} \}] 
\nonumber \\ && 
+ \sum_{i=2}^{k_n} \sum_{j:1 \leq |j-i|
\leq 3} \!\!\!\! {\rm Cov} [X_i \1 \{ E_{x,i} \} , X_j \1 \{ E_{x,j} \} ]
\nonumber\\ 
& = &
O(n^{4 \eps +  \sigma  + 2 \alpha(2 \eps + \sigma-1) }) .
\eea 
For $\alpha \geq 1$,
the bound in 
 (\ref{0712e}) tends to zero 
 as $n \to \infty$, 
since $1/2 < \sigma < 2/3$
and $\eps$ is small (\ref{epsdef}).

By (\ref{kndef}), the cells $\beta^x_i$,
 $i=1,\ldots,k_n$, have width asymptotic to $n^{2 \eps + \sigma -1}$
and height $n^{-\sigma}$, so the mean number
of  points of $\Po_n$  in one of these cells is
asymptotic to $n^{2 \eps}$; hence 
%the mean number ofwe have
%\[ \frac{1}{k_n} = n^{\sigma-1+\eta} + 
%O \left( n^{2\sigma-2+2\eta} \right),
%\] and
for any cell $\beta^x_i$ or $\beta^y_i$,
 $i=1,\ldots,k_n$, the probability that   the cell
 contains no point of $\Po_n$ is given by
$ \exp \{-n^{2\eps}(1+o(1))\}$.  
Hence for $n$ large enough, and $i=2,\ldots,k_n$,  we have
$
\Pr[E^c_{x,i}] 
%= \Pr[E^c_{y,i}]
 \leq \exp (-n^\eps),
$
and hence by (\ref{0712g}),
 \bea \label{0712f}
 \Var [ X_i \1 \{ E^c_{x,i} \} ] \leq \Exp [
X_i^2 | E^c_{x,i} ] \Pr[ E^c_{x,i} ] 
& \leq &
2^\alpha \Exp[ (N^x_i)^2 ] \Pr [
E^c_{x,i}]
\nonumber \\
&  = & O( n^{6 \eps + 2\sigma-1} \exp (-n^\eps)). ~ 
\eea
Hence by Cauchy-Schwarz
we have 
\bea \label{0712j} \Var
\left[ \sum_{i=2}^{k_n} X_i \1 \{ E^c_{x,i} \} \right] & = &
\sum_{i=2}^{k_n} \Var[ X_i \1 \{ E^c_{x,i} \} ] + \sum_{i \neq j}
{\rm Cov} [ X_i \1 \{ E^c_{x,i} \}, X_j \1 \{ E^c_{x,j} \}] \nonumber\\
& = & O \left( k_n^2 n^{6 \eps + 2\sigma-1} \exp (-n^\eps)
\right)  \to 0, 
\eea 
as $n \to \infty$. 
Then by
 (\ref{0712e}), (\ref{0712j}), and the analogous
estimates for $Y_i$, along with the Cauchy-Schwarz inequality, we
obtain for $\alpha \geq 1$ that
 \bea \label{0712z} \Var \left[ \sum_{i=2}^{k_n} X_i \1 \{
E_{x,i} \}+\sum_{i=2}^{k_n} Y_i \1 \{ E_{y,i} \}+\sum_{i=2}^{k_n}
X_i \1 \{ E^c_{x,i} \}+\sum_{i=2}^{k_n} Y_i \1 \{ E^c_{y,i} \}
\right] \to 0, \eea as $n \to \infty$. 
By (\ref{0712v}) with
 (\ref{0712x}),
 (\ref{0712z}),
and Cauchy-Schwarz again, we obtain the first part of (\ref{0525aa}). 
The argument for $\Po_n^0$ is the same as for $\Po_n$,
so we have (\ref{0525aa}).

Now suppose  $0<\alpha<1$.
We obtain (\ref{0525bb}) and (\ref{0525bbb}) in a similar
 way to (\ref{0525aa}),
since (\ref{0802e}) implies that
$$
\Var(n^{(\alpha-1)/2}(X_1 + Y_1 + Z)) = O(n^{6 \eps + 2 \sigma -2
+ \alpha(4 \eps + 2 \sigma -1) })
$$
 and 
(\ref{0712e}) implies 
$$
\Var \left( n^{(\alpha-1)/2}\sum_{i=2}^{k_n} X_i{\bf 1}\{E_{x,i}\}
\right)
 = O(n^{4 \eps +  \sigma -1
+ \alpha(4 \eps + 2 \sigma -1) }),
$$
 and both of these bounds tend to zero when $0 < \alpha < 1$,
$1/2 < \sigma < 2/3$, and $\eps $ is small (\ref{epsdef}).
$\square$ \\

To prove those parts of Theorem \ref{mainth}
which refer to the binomial process $\X_n$, we need
further results comparing the processes $\X_n$ and $\Po_n$
when they are coupled as in Lemma \ref{lem0803}.

\begin{lemma}
\label{lem0817}
Suppose $\alpha \geq 1$. With $\X_n$ and $\Po_n$
coupled as in Lemma \ref{lem0803}, we have that as $n \to \infty$
\bea
\label{0525aaa}
 \LL^\alpha( \X_n ; C_n) - \LL^\alpha(\Po_n; C_n) \inL 0
~~{\rm and} ~~
 \LL^\alpha( \X_n^0 ; C_n) - \LL^\alpha(\Po_n^0; C_n) \inL 0.
%~~~ (\alpha \geq 1).
\eea
\end{lemma}
\proof
Let  $\Po_n$ and $\X_m$ ($m \in \N$) be coupled as 
described in Lemma \ref{lem0803}. 
Given $n$, for $m \in \N$
define the event
$$
E_{m,n} :=
\cap_{1 \leq i \leq k_n}(\{ \X_{m-1} \cap \beta_i^x \neq \emptyset\}
\cap  \{\X_{m-1} \cap \beta_i^y \neq \emptyset\}),
$$
with the sub-cells $\beta_i^x$ and $\beta_i^y$ 
of $B_n$ as defined near the start of Section \ref{totallength}.
Then by similar arguments to those for $\Pr[E_{x,i}^c]$ above,
we have
$$
\Pr[E_{m,n}^c] = O (n^{1-\sigma-2\eps}
\exp(-n^\eps/2) ), ~~~~m \geq n/2 +1.
$$
As in the proof of Lemma \ref{lem0803},
let $Y_m$ denote the in-degree
of vertex $\bX_m$ in the MDST on $\X_m$.
Then
\[ \left| \LL^\alpha (\X_m ; C_n) -
\LL^\alpha (\X_{m-1} ; C_n ) \right|
\leq  (Y_m+1) \1\{\bX_m \in C_n \} 
\left( (3n^{2\eps +\sigma -1})^{\alpha}
+ 2^{\alpha/2} \1\{E_{m,n}^c\} \right).\]
Thus, given $N(n)$,
\bean
%\label{0807a}
 \left| \LL^\alpha (\X_n ; C_n) -
\LL^\alpha (\Po_n ; C_n) \right|
& \leq  \sum_{m=\min (N(n),n)}
^{\max (N(n),n)}
 & (Y_m+1) \1\{\bX_m \in C_n \} 
\\
& & \times \left( 3^\alpha n^{\alpha(2\eps+\sigma-1)}
+2^{\alpha/2} \1\{E_{m,n}^c\} \right).
\eean
Since $C_n$ has area less than $2n^{\eps-1/2}$, 
by \eq{harmonicbd} there exists a constant $C$ such that,
for $n$ 
 sufficiently large 
and $N(n) \geq n/2+1$,
\bea
\label{0807b}
\Exp \left[ \left.
\left( \left| \LL^\alpha (\X_n ; C_n) -
\LL^\alpha (\Po_n ; C_n) \right|
\right) \right| N(n) \right] \leq
2^{\alpha/2} n  \1_{\{N(n) < n/2 +1\}}
\nonumber \\ + 
 C |N(n)-n| \log( \max (N(n),n))   
n^{\alpha(2\eps+\sigma-1)+\eps-1/2} \1_{\{N(n) \geq n/2 +1\}}.
\eea
By tail bounds for the Poisson distribution, we have
$nP[N(n) < n/2 +1] \to 0$ as $n \to \infty$, and hence,
taking expectations in (\ref{0807b}) and
using (\ref{logmax}),
we obtain 
\[ \Exp \left[ 
 \left| \LL^\alpha (\X_n ; C_n) -
\LL^\alpha (\Po_n ; C_n) \right|
 \right] = O(  n^{\alpha(2\eps+\sigma-1)+\eps}\log n) + o(1),\]
which tends to zero since $\alpha \geq 1$,
$1/2< \sigma <2/3$
and $\eps$ is small (see
(\ref{epsdef})). So
%, finally, by Cauchy-Schwarz 
we obtain
the unrooted part of (\ref{0525aaa}). The argument is the same
in the rooted case. $\square$ \\

%We shall need one final lemma on coupling of $\Po_n$ and $\X_n$.
\begin{lemma}
\label{depolem}
Suppose $\X_n$ and $\Po_n$ are coupled 
as described in Lemma \ref{lem0803}, with $N(n):= \card(\Po_n)$.
Let $\Delta(\infty)$ be given by Definition \ref{sstabdef} with
$H = \LL^1$, and
set $\alpha_1 := \Exp[\Delta(\infty)]$.
Then as $n \to \infty$
 we have
\bea
\label{eqdepo1}
  \LL^1 (\Po_n;S_{0,n}) - \LL^1 (\X_n;S_{0,n})
- n^{-1/2}\alpha_1 (N(n)-n)  \inLL 0; \\
\label{eqdepo2}
  \LL^1 (\Po_n^0;S_{0,n}) - \LL^1 (\X_n^0;S_{0,n})
- n^{-1/2}\alpha_1 (N(n)-n)  \inLL 0.
\eea
\end{lemma}
{\em Proof.}
The proof of the first part (\ref{eqdepo1}) follows that of eqn (4.5)
of \cite{penyuk1}, using our Lemma \ref{0331q} and the fact that
%for our present choice of
%$\xi$, given by (\ref{0802})  with $\alpha =1$ and with
%the ordering $\postar$ 
the functional $\LL^1$ is homogeneous of order 1,
is strongly stabilizing
by Lemma \ref{stab},
and satisfies the moments condition \eq{ubm}
by Lemma \ref{ubmlem}.

As shown in the proof of Corollary \ref{0804c} 
(see in particular eqn \eq{0817b}), we have that
$  \LL^1 (\Po_n^0;S_{0,n}) - \LL^1 (\Po_n;S_{0,n})$
converges to zero in $L^2$ and
$  \LL^1 (\X_n^0;S_{0,n}) - \LL^1 (\X_n;S_{0,n})$
converges to zero in $L^2$. 
Therefore the second part (\ref{eqdepo2}) follows from 
(\ref{eqdepo1}).  $\square$\\

We are now in a position to prove Theorem \ref{mainth}.
We divide the proof into two cases: $\alpha \neq 1$ and $\alpha=1$. 
In the latter case,
to prove the result for 
 the Poisson process $\Po_n$, we
need to show that
$\LL^1 ( \Po_n ; B_n)$ and $\LL^1 ( \Po_n ; S_{0,n})$ are
asymptotically independent; likewise for $\Po_n^0$.
 We shall then obtain the results for the binomial process
$\X_n$ and for $\X_n^0$
from those for $\Po_n$ and $\Po_n^0$ via the
 coupling described in Lemma \ref{lem0803}.  \\

\noindent \textbf{Proof of Theorem \ref{mainth} for $\alpha \neq 1$.}
First suppose $0<\alpha <1$. For the Poisson case, we have
\bea
 n^{(\alpha-1)/2} \tLalph (\Po_n) = n^{(\alpha-1)/2}
\tLalph (\Po_n;S_{0,n}) + n^{(\alpha-1)/2} \tLalph (\Po_n;B_n)
\nonumber \\
+ n^{(\alpha-1)/2}
\tLalph (\Po_n;C_n).
\label{0804a}
\eea
The first term in the
right hand side of (\ref{0804a})
converges in distribution to $\NN(0,s_\alpha^2)$
by Theorem \ref{CLT} (iv), and the other two terms converge
in probability to 0 by eqns (\ref{0214h}) and (\ref{0525bb}).
Thus Slutsky's theorem yields the
first (Poisson) part of 
 (\ref{0727d}).
To obtain the second (binomial) part of (\ref{0727d}),
we use the coupling of Lemma \ref{lem0803}.
We write
\bea
n^{(\alpha -1)/2}\tLL^\alpha(\X_n) =
n^{(\alpha -1)/2}
\tLL^\alpha(\X_n;S_{0,n} ) 
+ n^{(\alpha -1)/2} (\tLL^\alpha (\Po_n; B_n \cup C_n) ) 
\nonumber \\
+ n^{(\alpha -1)/2} (
\tLL^\alpha (\X_n; B_n \cup C_n)
-
\tLL^\alpha (\Po_n; B_n \cup C_n)
 ) . ~
\label{0803f}
\eea
The first term in the right side of  (\ref{0803f}) is asymptotically
$\NN(0,t_\alpha^2)$ by Theorem \ref{CLT} (ii). The second term tends to
zero in probability by 
 (\ref{0214h}) and (\ref{0525bb}).
The third term tends to zero in probability by 
(\ref{0804f}). Thus we have 
 the binomial case of (\ref{0727d}).

The rooted case (\ref{0727a}) is similar. Now,
for the first (Poisson) part  of
 (\ref{0727a}), we use Corollary \ref{0804c} (iv) with
(\ref{0802a}) and (\ref{0525bbb}), and Slutsky's theorem. The
second part of
(\ref{0727a}) follows from the analogous statement to (\ref{0803f})
with the addition of the origin, using Corollary \ref{0804c} (ii)
with (\ref{0802a}), (\ref{0525bbb}), (\ref{0804g}), and
Slutsky's theorem again.

Next, suppose $\alpha >1$. 
%Now suppose $\alpha>1$. 
We have 
\bea
\label{0804b}
\tLalph (\Po_n) = 
\tLalph (\Po_n;S_{0,n}) + 
\tLalph (\Po_n;C_n)
+
\tLalph (\Po_n;B_n).\eea
The first term in the right hand side
converges to 0 in probability, by Theorem \ref{CLT} (iii).
The second term also converges to 0 in probability,
by the first
part of (\ref{0525aa}). Then by (\ref{bdry2})
and Slutsky's theorem, we obtain
the first (Poisson) part of (\ref{0727f}).
To obtain the rooted version,
i.e. the first part of (\ref{0727c}),
we replace $\Po_n$ by $\Po_n^0$
in (\ref{0804b}), and
 combine (\ref{bdry1}) 
  with Corollary \ref{0804c} (iii)  
and the second part of (\ref{0525aa}),
and apply Slutsky's theorem again.

To obtain the binomial versions of the results (\ref{0727c}) and
(\ref{0727f}), we again
make use of the coupling described in Lemma \ref{lem0803}.
We have 
\bea
\label{0804h}
\tLalph (\X_n)  =   
\tLalph (\X_n;S_{0,n}) + 
\tLalph (\X_n;C_n)
+
\tLalph (\X_n;B_n). 
\eea
The first term in the right hand
side converges in probability to zero by
Theorem \ref{CLT} (i). 
The second term converges in probability to zero by 
the first part of (\ref{0525aa})  and the first part of (\ref{0525aaa}).
The third part converges in distribution to
 $\tFalph^{\{1\}} + \tFalph^{\{2\}}$ by
by  (\ref{bdry2X}). Hence, 
 Slutsky's theorem
yields the binomial part of (\ref{0727f}). 

Similarly, by replacing $\Po_n$ by $\Po^0_n$ and $\X_n$
by $\X^0_n$ in (\ref{0804h}), and using Corollary \ref{0804c}
(i), the second part of (\ref{0525aa})
%, (\ref{0803e}), the second part
and of (\ref{0525aaa}),
(\ref{bdry1X}) and Slutsky's theorem,
we obtain the binomial part of (\ref{0727c}).
This completes the proof for $\alpha \neq 1$. \\

\noindent \textbf{Proof of Theorem \ref{mainth} for $\alpha = 1$:
 the Poisson case.}
We now prove the first part of (\ref{0727b}) and the first 
part of (\ref{0727e}).
Given $n$, set $q_n := 4 \lfloor n^{\eps   + \sigma -1/2} \rfloor$. 
Split each
cell $\Gamma^x_i$ of $C_n^x$ into $  4 q_n $ 
 rectangular sub-cells, by splitting
the horizontal edge into $q_n$ segments and the vertical edge into
4 segments by a rectangular grid. Similarly, split each
cell $\Gamma_i^y$ by splitting the vertical edge into $q_n$ segments
and the horizontal edge into 4 segments. 
Finally, add a single square sub-cell  in the top right-hand
corner of $C_n^0$, of side $(1/4)n^{\eps-1/2 }$, and denote this
``the corner sub-cell''.

The total number of all such sub-cells is $1+8 k_nq_n \sim 32 n^{(1/2)-\eps}$.
Each of the sub-cells has width asymptotic to $(1/4)n^{\eps-1/2}$
and height asymptotic to $(1/4)n^{\eps-1/2}$, and so
 the area of each cell is asymptotic to
$(1/16)n^{2 \eps -1}$.  So for large $n$,
for each of these sub-cells,  the probability
that it contains no point of $\Po_n$
is bounded by $\exp (-n^{\eps})$.

Let $E_n$ be the event that
each of the  sub-cells
 described above contains at least one point of $\Po_n$. Then
\bea
\label{0807d}
 \Pr[ E_n^c ] = 
%\exp(-n^\eps) +
 O \left( n^{(1/2)-\eps}
\exp (-n^{\eps}) \right)
  \to 0 . 
\eea
Suppose $\bx$ lies on the lower boundary of  $S_{0,n}$.
Consider the rectangular sub-cell of $\Gamma_i^x$
lying just to the left of the sub-cell directly below
$\bx$ (or the corner sub-cell if that lies just to the left of
the sub-cell directly below $\bx$). 
 All points $\by$ in this sub-cell
satisfy $\by \preccurlyeq^* \bx$, and
 for large $n$, satisfy  $\|\by -\bx\| <
(3/4)n^{\eps-1/2}$, whereas the nearest point to $\bx$
in $B_n$ is at a distance at least $(3/4)n^{\eps -1/2}$.
Arguing similarly for $\bx$ on the left boundary of $S_{0,n}$,
and using the triangle inequality, we see that
if $E_n$ occurs, no point in $S_{0,n}$ can be connected to 
any point in $B_n$,
provided $n$ is sufficiently large.

For simplicity of notation, set 
$X_n := \tilde \LL^1 ( \Po_n ; B_n)$
and $Y_n := \tilde \LL^1 ( \Po_n ; S_{0,n})$. 
Also, set
 $X := \tDone^{\{ 1 \}} + \tDone^{\{ 2 \}}$ 
and $Y\sim \NN(0,s_1^2)$, independent of $X$,
with $s_1$ as given in Theorem \ref{CLT}.
We know from 
Theorem \ref{thmbdry} and Theorem \ref{CLT} that
$X_n \tod X$ and $Y_n \tod Y$ as $n \to \infty$.

We need to show that
 $X_n + Y_n \tod X+ Y$, where $X$ and $Y$
are independent random variables.
% By the Cramer-Wold device, it suffices to prove that
%$X_n +Y_n \tod X + Y$ for independent $X, Y$. 
We show this by convergence of the
characteristic function,
\bea \label{0210a}
\Exp { \left[ \exp{ \left( i t (X_n +Y_n ) \right) } \right] }
\longrightarrow \Exp { \left[ \exp{ \left( i t X \right) } \right] }
\Exp { \left[ \exp{ \left( i t Y \right) } \right] }.
\eea
With $\omega$ denoting the configuration of points in
$C_n$, we have
\bean
\Exp \left[ \exp{ \left( it(X_n + Y_n ) \right) } \right]
& = & \int_{E_n} \Exp \left[ \left. e^{itX_n} e^{itY_n}
\right| \omega \right] \ud \Pr (\omega) + \Exp \left[
e^{it(X_n +Y_n)} \1_{ E_n^c } \right] \\
& = & \int_{E_n} \Exp \left[ e^{itX_n}
\right]
\Exp \left[ \left. e^{itY_n}
\right| \omega \right] \ud \Pr (\omega) + \Exp \left[
e^{it(X_n +Y_n)} \1_{ E_n^c } \right] , \eean
where we have used the fact that $X_n$ and $Y_n$ are
conditionally independent, given $\omega \in E_n$, for $n$ sufficiently large,
and that $X_n$ is independent of the configuration in $C_n$. Then
$\Exp [ e^{it(X_n+Y_n)} \1_{ E_n^c } ] \to 0$ 
as $n \to \infty$, since $P[E_n^c]\to 0$.
So 
\[
\Exp \left[ \exp{ \left( it(X_n + Y_n ) \right) } \right]
- \Exp \left[ e^{itX_n} \right] \Exp \left[ e^{itY_n} \1_{E_n} \right] \to 0,
\]
and we obtain (\ref{0210a})
since $\Exp [ e^{itY_n} \1_{E_n} ] = \Exp
[ e^{itY_n} ] - \Exp [ e^{itY_n} \1_{E_n^c} ]$,
$\Exp [ e^{itY_n} \1_{E_n^c} ]\to 0$,
$\Exp[ e^{itX_n}] \to \Exp[
e^{itX}]$, and $\Exp[ e^{itY_n}] \to \Exp[
e^{itY}]$ as $n \to \infty$.

We can now prove the first (Poisson) part of (\ref{0727e}). 
We have the $\alpha=1$ case of (\ref{0804b}). 
%The above argument
%deals with the contributions from $S_{0,n}$ and $B_n$. 
The contribution
from $C_n$ converges in probability to 0 by the first part of (\ref{0525aa}).  
Slutsky's theorem and 
(\ref{0210a}) then give the first (Poisson) part of (\ref{0727e}). 
The rooted Poisson case (\ref{0727b}) follows from the rooted version of (\ref{0804b}), this time
applying the argument for (\ref{0210a})
taking $X_n := \tilde \LL^1 (\Po^0_n;B_n)$, $Y_n:= \tilde \LL^1
(\Po^0_n;S_{0,n})$ and $X$, $Y$ as before,
and then using the second part of (\ref{0525aa}) and Slutsky's theorem again. Thus 
we obtain the first (Poisson) part of (\ref{0727b}). \\

\noindent \textbf{Proof of Theorem \ref{mainth} for $\alpha = 1$:
 the binomial case.}
It remains for us to 
prove the second part of (\ref{0727b}) and the second 
part of (\ref{0727e}).
To do this, we use the coupling of Lemma \ref{lem0803}
once more. Considering first the unrooted case, 
we here set $X_n :=  \LL^1 (\X_n;B_n)$ and 
$Y_n:= \LL^1 (\X_n;S_{0,n})$.
Set $X'_n := \LL^1 (\Po_n;B_n)$ and 
$Y'_n:=  \LL^1 (\Po_n;S_{0,n})$
 (note that all these random variables are uncentred). 

Set $Y \sim \NN(0,s_1^2)$ with $s_1$ as given in Theorem \ref{CLT}. 
Set $X := \tDone^{\{1\}} + \tDone^{\{2\}}$, independent of
$Y$. Then by \eq{0210a} we have (in our new notation)
\bea
X'_n -\Exp X'_n  + Y'_n - \Exp Y'_n \tod X + Y.  
\label{0818}
\eea
By (\ref{0803d}), we have $X_n-X'_n \toP 0$ and $\Exp X_n - \Exp X'_n \to 0$. 
Also, with $\alpha_1 $ as defined in Lemma \ref{depolem},
 eqn (\ref{eqdepo1}) of that result gives us
\bea
Y'_n - Y_n - n^{-1/2}\alpha_1 (N(n)-n)  \inLL 0
\label{0818a}
\eea
so that 
 $\Exp[Y'_n] - \Exp[ Y_n] \to 0$.
Combining these observations with (\ref{0818}), 
and using Slutsky's theorem, we obtain
\bea
X_n -\Exp X_n  + Y_n - \Exp Y_n  + n^{-1/2}\alpha_1 (N(n) -n)  \tod X + Y.  
\label{0818b}
\eea
By Theorem \ref{CLT} (iii) we have $\Var(Y'_n) \to s^2_1$ as $n \to \infty$.
By
 (\ref{0818a}), and the independence of $N(n)$ and $Y_n$, we have
\bea
s^2_1 = \lim_{n \to \infty} \Var[ Y_n + n^{-1/2}\alpha_1(N(n)-n) ]
= \lim_{n \to \infty} ( \Var [Y_n] + \alpha_1^2)  
% Y_n - \Exp Y_n  + n^{-1/2}\alpha_1 (N(n) -n)  \tod  Y  
\label{0818c}
\eea
so that $\alpha_1^2 \leq s_1^2$.   
Also, $ n^{-1/2}\alpha_1( N(n)-n)$ is
 independent of $X_n + Y_n$, and asymptotically $\NN(0,\alpha_1^2)$.
%By considering
Since the $\NN(0,s^2)$  characteristic function is
$\exp(-s^2t^2/2)$, for all $t \in \R$ we obtain 
from \eq{0818b} that 
$$
\Exp[ \exp(it(X_n - \Exp X_n + Y_n  - \Exp Y_n))] \to \exp( -(s_1^2 - \alpha_1^2)t^2/2 ) \Exp [\exp(it X )] 
$$
so that 
\bea
X_n -\Exp X_n  + Y_n - \Exp Y_n   \tod X + W,  
\label{0818d}
\eea
where $W \sim \NN(0,s_1^2 -\alpha_1^2)$, and $W$ is independent of $X$.

We have the $\alpha=1$ case of 
(\ref{0804h}). 
 By the first part of
  \eq{0525aa} and the first part of (\ref{0525aaa}),
the contribution from $C_n$ tends to zero in probability.
 Hence
by (\ref{0818d}) and Slutsky's theorem,
we obtain
 the second  (binomial) part of (\ref{0727e}).

For the rooted case, 
%i.e. for  the second part of (\ref{0727b}), 
we apply the argument for (\ref{0818d}), now
taking $X_n := \LL^1 (\X^0_n;B_n)$, $Y_n:=  \LL^1
(\X^0_n;S_{0,n})$, with $X$, $Y$ and $W$ as before.
% $Y \sim \NN(0,s_1^2)$ as before, but now with $X = \tDone^{\{1\}}  
%+ \tDone^{\{2\}}.$
The rooted case of (\ref{0818}) follows from the rooted 
case of (\ref{0210a}), and now we have
 $X_n-X'_n \toP 0$ and $\Exp X_n - \Exp X'_n \to 0$ by 
 (\ref{0803e}).
In the rooted case (\ref{0818a}) still holds by
(\ref{eqdepo2}),
and then we obtain the rooted case
of (\ref{0818d}) as before. 

To obtain  the second (binomial) part of
 (\ref{0727b}), we start with
the rooted version of the $\alpha =1$ case of (\ref{0804h}). 
 By the second part of
  \eq{0525aa} and of (\ref{0525aaa}),
the contribution from $C_n$ tends to zero in probability.
 Hence by the rooted version of (\ref{0818d}) and Slutsky's theorem,
we obtain the second part of (\ref{0727b}).

 This completes 
the proof of the $\alpha=1$ case,
and hence the proof of Theorem \ref{mainth}
 is complete.
$\square$ \\

\begin{center} \textbf{Acknowledgements} \end{center}
The first author began this work while at the University of
Durham, and was also supported by the Isaac Newton Institute for
Mathematical Sciences, Cambridge. The second author was
 supported by the EPSRC.

\end{document}